\documentclass[reqno]{amsart}

\usepackage{amssymb,amsmath,amsthm,latexsym,amscd}
\usepackage{epsfig, textcomp}
\usepackage{psfrag}
\usepackage{graphicx}
\usepackage{tikz}
\usepackage{float}

\newtheorem{Theorem}{\bf Theorem}
\newtheorem{lemma}[Theorem]{\bf Lemma}
\newtheorem{proposition}[Theorem]{\bf Proposition}
\newtheorem{corollary}[Theorem]{\bf Corollary}

\newtheorem{definition}[Theorem]{\bf Definition}

\newtheorem{remark}[Theorem]{\bf Remark}
\newtheorem{theorem}[Theorem]{\bf Theorem}
 
\def\scfig #1 #2 {\resizebox{#2}{!}{\includegraphics{#1}}}

\newcommand{\be}{\begin{equation}}
\newcommand{\ee}{\end{equation}}

\def\hpic #1 #2 {\mbox{$\begin{array}[c]{l} 
\epsfig{file=#1,height=#2}\end{array}$}}
\def\wpic #1 #2 {\mbox{$\begin{array}[c]{l} 
\epsfig{file=#1,width=#2}\end{array}$}}

\begin{document}
\title{From a Kac algebra subfactor to Drinfeld double}
\author{Sandipan De}
\address{Stat-Math Unit\\Indian Statistical Institute, 8th Mile, Mysore Road\\ Bangalore-560059}
\email{sandipan$\_$vs@isibang.ac.in}


\keywords{Subfactors, Kac algebras, Planar algebras, Drinfeld double}
\subjclass[2010]{46L37, 16S40, 16T05}
\begin{abstract}
Given a finite-index and finite-depth subfactor, we define the notion of \textit{quantum double inclusion} - a certain unital inclusion of 
von Neumann algebras constructed from the given subfactor - which is closely related to that of Ocneanu's asymptotic inclusion.
We show that the quantum double inclusion when applied to the Kac algebra subfactor $R^H \subset R$ produces Drinfeld double of $H$ where $H$ is a 
finite-dimensional Kac algebra acting outerly on the hyperfinite $II_1$ factor $R$ and $R^H$ denotes the fixed-point subalgebra. 
More precisely, quantum double inclusion of $R^H \subset R$
is isomorphic to $R \subset R \rtimes D(H)^{cop}$ for some outer action of $D(H)^{cop}$ on $R$ where $D(H)$ denotes the Drinfeld double of $H$.
\end{abstract}
\maketitle
\section*{Introduction}
Ocneanu's asymptotic inclusion \cite{Oc1988} for finite-index and finite-depth inclusion of hyperfinite $II_1$ factors could be viewed as the
subfactor analogue of Drinfeld's quantum double construction. This connection has been clarified by a number of authors
including Evans-Kawahigashi \cite{EK1995}, Izumi \cite{Izumi2000, Izumi2001} and M\"{u}ger \cite{Mug2003} and the precise formulation of 
Ocneanu's analogy is that the bimodule category arising from the asymptotic inclusion is the 
 Drinfeld center of the bimodule category of the original subfactor. 
In \cite{Mug2003} M\"{u}ger remarked that - see \cite[Remark 8.7]{Mug2003} - the asymptotic subfactor of $R^H \subset R$,
where $H$ is a 
finite-dimensional Kac algebra acting outerly on the hyperfinite $II_1$ factor $R$ and $R^H$ is the fixed-point subalgebra, cannot be isomorphic
to $R^{D(H)} \subset R$ or its dual, where $D(H)$ is the Drinfeld double of $H$, since the index of
the asymptotic subfactor of $R^H \subset R$ coincides with the index of $R^H \subset R$. This gave rise to a 
natural question: Starting with $R^H \subset R$, is it possible to obtain $R^{D(H)} \subset R$ or its dual by some procedure.
This answer is perhaps known to the experts but does not appear to have been published anywhere with details and we feel it deserves to be better known.
In this article we show that
a modification of the asymptotic inclusion, which we call \textit{quantum double inclusion}, does the job. More precisely,
we show that quantum double inclusion of $R^H \subset R$ is isomorphic to $R \subset R \rtimes D(H)^{cop}$ for some outer action of $D(H)^{cop}$ on $R$ 
and this is the main result of this article. All proofs in this paper are pictorial and are a testament to the elegance of planar algebra 
techniques.

Let $N \subset M$ be a finite-index subfactor of finite-depth and let $N ( = M_0) \subset M ( = M_1) \subset M_2 \subset M_3 \subset \cdots$ be the Jones' basic
construction tower of $N \subset M$. Let $M_{\infty}$ denote the $II_1$ factor obtained as the von Neumann closure 
$(\cup_{n = 0}^{\infty} M_n)^{\prime \prime}$ in the GNS representation with respect to the trace on $\cup_{n = 0}^{\infty} M_n$. Recall that - see 
\cite{EK1993}, for instance - the inclusion $M \vee (M^{\prime} \cap M_{\infty}) \subset M_{\infty}$ is defined as the asymptotic inclusion constructed 
from $N \subset M$ where, of course, $M \vee (M^{\prime} \cap M_{\infty})$ denotes the von Neumann algebra generated by 
$M$ and $M^{\prime} \cap M_{\infty}$. We define the inclusion $N \vee (M^{\prime} \cap M_{\infty}) \subset M_{\infty}$ to be the 
\textit{quantum double inclusion} of $N \subset M$ where $N \vee (M^{\prime} \cap M_{\infty})$ denotes the von Neumann algebra generated by
$N$ and $M^{\prime} \cap M_{\infty}$.

One of the main steps towards understanding the quantum double inclusion associated to the subfactor $R^H \subset R$ is to construct a model of 
this. Given any finite-dimensional Kac algebra $H$, let $H^i$, where $i$ is any integer, denote $H$ or $H^*$ according as $i$ is odd or even.
In $\S 3$ we construct a subfactor $\mathcal{N} \subset \mathcal{M}$ where $\mathcal{N} = ((\cdots \rtimes H^{-3} \rtimes H^{-2} \rtimes 
 H^{-1}) \otimes (H^2 \rtimes H^3 \rtimes \cdots))^{\prime \prime}, \ \mathcal{M} = (\cdots \rtimes H^{-1} \rtimes H^0 \rtimes H^1 \rtimes 
 \cdots)^{\prime \prime}$ and show that $\mathcal{N} \subset \mathcal{M}$ is a model for the quantum double inclusion of $R^H \subset
 R$. In $\S 5$, we give an explicit description of the planar algebra 
 associated to the subfactor $\mathcal{N} \subset \mathcal{M}$ which turns out to be an interesting planar subalgebra of
 $^{*(2)}\!P(H^*)$ (the adjoint of the $2$-cabling of the planar algebra of $H^*$). 
 Finally, this description of the planar algebra 
 of $\mathcal{N} \subset \mathcal{M}$ is used in $\S 6$ to prove the main result, namely Theorem \ref{result}, of this article.
 The proofs all rely on explicit pictorial computations in the planar algebra of $H^*$.

We are grateful to the referee for pointing out that an alternate approach to this result
 - possibly in greater generality - could be based on the Longo-Rehren inclusion and we hope to turn to this approach in a subsequent paper.

We give below a brief section-wise description of the contents of this paper.

$\S 1$: The goal of this section is to summarise
relevant facts concerning crossed products by  Kac algebras. We begin with recalling the notion of action of a Kac algebra on a 
complex $*$-algebra and given such an action, we describe the construction of 
the crossed product algebra. We then introduce the notion of infinite iterated crossed products.
Using this we define a family, indexed by positive integers, of inclusions of (infinite-dimensional) algebras which will be used in $\S 3$ 
in order to understand the model for quantum double inclusion of $R^H \subset R$.

$\S 2$: In section $\S 2.1$, we collect together results concerning subfactor planar algebras. We begin with introducing a few important tangles we need. 
Next, we discuss two methods of constructing new planar algebras from the old, namely, cabling and adjoint and also recall two important theorems - one
concerning `generating set of tangles' and the other being the fundamental theorem due to Jones relating subfactors and subfactor planar 
algebras - which we shall use in $\S 5$ in order to describe $P^{\mathcal{N} \subset \mathcal{M}}$. In $\S 2.2$ we briefly discuss
the planar algebra associated to a Kac algebra in terms of generators and relations and identify its vector 
spaces explicitly in terms of iterated crossed products of the Kac algebra and its dual.

The material of the first 2 sections is all very well known and is meant just to establish notation for the
convenience of the reader.

$\S 3$: This section is devoted to constructing a model for the quantum double inclusion of $R^H \subset R$. The main result of this section is 
Proposition \ref{qdim} which shows that the subfactor 
$\mathcal{N} \subset \mathcal{M}$, where $\mathcal{N} = ((\cdots \rtimes H^{-3} \rtimes H^{-2} \rtimes 
 H^{-1}) \otimes (H^2 \rtimes H^3 \rtimes \cdots))^{\prime \prime}$ and $\mathcal{M} = (\cdots \rtimes H^{-1} \rtimes H^0 \rtimes H^1 \rtimes 
 \cdots)^{\prime \prime}$, is a model for the quantum double inclusion of $R^H \subset R$.
 
 $\S 4$: This section begins with $\S 4.1$ which paves the way for constructing the basic construction tower of the subfactor 
 $\mathcal{N} \subset \mathcal{M}$ by studying some finite-dimensional basic constructions associated to inclusions of finite iterated 
 crossed product algebras, the main result being Proposition \ref{basic cons}.
In $\S 4.2$ the basic construction tower of $\mathcal{N} \subset \mathcal{M}$ is explicitly constructed in Proposition \ref{towerbasic}.
In $\S 4.3$ we compute the relative commutants of the basic construction towers using Ocneanu's compactness theorem.

 $\S 5$: The penultimate $\S 5$ studies the planar algebra associated to $\mathcal{N} \subset \mathcal{M}$.
 The main result of this section is Theorem \ref{maintheorem} which describes
the subfactor planar algebra associated to $\mathcal{N} \subset \mathcal{M}$.

$\S 6$: In the final $\S 6$ we prove the main result, Theorem \ref{result}, which says that $\mathcal{N} \subset \mathcal{M}$ is 
isomorphic to $R \subset R \rtimes D(H)^{cop}$ for some outer action of $D(H)^{cop}$ on the hyperfinite $II_1$ factor $R$.

\section{Crossed product by Kac algebras}
In this section we briefly review the notion of crossed product by a Kac algebra. For a detailed exposition of this concept, the reader may consult 
\cite{Jjo2008}. We refer to $\S 4$ of \cite{KdyLndSnd2003} for the standard
facts concerning finite-dimensional Kac algebras which will be used frequently throughout this article. 
Unless otherwise specified, $H (= (H, \mu, \eta, \Delta, \epsilon, S, *))$ will denote a finite-dimensional Kac algebra and
$\delta$, the positive square root of $dim ~H$. We set $H^i = H$ or $H^*$ according as $i$ is odd or even. The letters $x, y, z$ as well as the symbols
$x^i, y^i, z^i$, $i$ being any odd integer, will always denote an element of $H$. 
The letters $f, g$ as well as the symbols $f^i, g^i$, where $i$ is any even integer, will always represent an element of 
$H^*$. For the rest of this paper, the unique non-zero idempotent integrals of $H^*$ and $H$ will be denoted by
$\phi$ and $h$ respectively and moreover, for any non-negative integer $i$, the symbols $\phi^i$ and $h^i$ will always denote a copy of 
$\phi$ and $h$ respectively. It is a fact that $\phi(h) = \frac{1}{dim ~H}$.

\begin{definition}\label{def}

By an action of $H$ on a finite-dimensional complex $*$-algebra $A$
we will mean a linear map $\alpha : H \rightarrow End(A)$
(references to endomorphisms without further qualification will be to $\mathbb{C}$-linear endomorphisms)
satisfying (i) $\alpha_1 = id_A$, (ii) $\alpha_{xy} =\alpha _x \circ \alpha _y$, (iii) $\alpha_x(1_A) = \epsilon(x) 1_A$,
(iv) $\alpha_x(ab) = \alpha_{x_1}(a)\alpha_{x_2}(b)$, and (v) $\alpha_x(a)^* = \alpha_{Sx^*}(a^*)$ for all $x,y \in H$ and $a,b \in A$. 
To clarify notation, $\alpha_x$ stands for $\alpha(x)$
and $\Delta(x)$ is denoted by $x_1 \otimes x_2$ (a simplified version of the Sweedler cooproduct notation). For simplicity, we often use the notation $x . a$ to denote $\alpha_x(a)$.
\end{definition}

Suppose that $\alpha$ is an action of $H$ on $A$. The crossed product algebra, denoted $A \rtimes_\alpha H$ (or mostly, simply as 
$A \rtimes H$, when the action is understood) is defined to be the $*$-algebra whose underlying vector space is $A \otimes H$ 
(where we denote $a \otimes x$ by $a \rtimes x$) and the multiplication is defined by 
\begin{align*}
 (a \rtimes x)(b \rtimes y) = a\alpha_{x_1}(b) \rtimes x_2y.
\end{align*}
The $*$-structure on $A \rtimes H$ is given by $(a \rtimes x)^* = \alpha_{x_1^*}(a^*)
\rtimes x_2^*$. 
This is an algebra with unit $1_A \rtimes 1_H$ and
there are natural inclusions (which are $*$-maps also) of algebras $A \subseteq A \rtimes H$ given by $a \mapsto a \rtimes 1_H$ and 
$H \subseteq A \rtimes H$ given by $x \mapsto 1_A \rtimes x$.
We draw the reader's attention to a notational abuse of which we will often be guilty. We denote elements of
a tensor product as decomposable tensors with the understanding that there is an implied omitted summation.
Thus, when we write `suppose $f \otimes x \in H^* \otimes H$', we mean `suppose $\sum_i f^i \otimes x^i \in H^* \otimes H$' 
(for some $f^i \in H^*$ and $x^i \in H$, the sum over a finite index set).

There is a natural action of $H^*$ on $H$ given by 
 $f . x = f(x_2) x_1$ for $f \in H^*, x \in H$. Similarly we have action of $H$ on $H^*$.
 If $H$ acts on $A$, then $H^*$ also acts on $A \rtimes H$ just by acting on $H$-part and ignoring the $A$-part, meaning that, 
 $f. (a \rtimes x) = a \rtimes f.x = f(x_2) \ a \rtimes x_1$ for $f$ in $H^*$ and $a \rtimes x \in A \rtimes H$ and consequently, we can construct $A \rtimes H \rtimes
 H^*$. Continuing this way, we may construct $A \rtimes H \rtimes H^* \rtimes \cdots$. 
 Note that $H^i$ and $H^j$ commute whenever $|i-j| \geq 2$.

 For integers $i \leq j$, we define $H_{[i, j]}$ to be the crossed product algebra $H^i \rtimes H^{i+1} \rtimes \cdots \rtimes H^j$. 
 If $i = j$, we will simply write $H_i$ to denote $H_{[i, i]}$ and if $i > j$, we take $H_{[i, j]}$ to be $\mathbb{C}$.
 A typical element of $H_{[i, j]}$ will be denoted by $x^i / f^i \rtimes f^{i+1} / x^{i+1} \rtimes \cdots$ ($j-i+1$ terms). 
 For instance, a typical element of $H_{[0, 3]}$
 will be denoted by $ f^0 \rtimes x^1 \rtimes f^2 \rtimes x^3$.
 The multiplication rule shows that if $p \leq i \leq j \leq q$, the natural inclusion of $H_{[i,j]}$ into $H_{[p,q]}$ is an algebra map.
 For any positive integer $l$, we use the notation $A(H)_l$ to denote the crossed product algebra
 $\underbrace{H \rtimes H^* \rtimes \cdots}_{\text{$l$ times}}$. Note that given integers $i, l$ with $l$ positive, 
 $H_{[i, i+l-1]} = A(H)_l$ or $A(H^*)_l$ according as $i$ is odd or even.

 Define the algebra $H_{(-\infty,\infty)}$ to be the `union' of all the $H_{[i,j]}$. We may suggestively write 
 $ H_{(-\infty,\infty)} = \cdots \rtimes H \rtimes H^* \rtimes H \rtimes \cdots$ and represent a typical element of 
$H_{(-\infty, \infty)}$ as $\cdots \rtimes x^{-1} \rtimes f^0 \rtimes x^1 \rtimes \cdots$. We repeat that this means that a typical element of $H_{(-\infty,\infty)}$
is in fact a finite sum of such terms. Note that in any such term all but finitely many of the $f^i$ are $\epsilon$ and all but finitely many of 
the $x^i$ are $1$. 
Next, for any integer $m$, we define a subalgebra of $H_{(-\infty, \infty)}$ which, in suggestive notation, is
$H_{(-\infty, m]}$.
A little more clearly, it consists of all (finite sums of) elements $\cdots \rtimes x^{-1} \rtimes f^0 \rtimes x^1 \rtimes \cdots$ of $H_{(-\infty,
\infty)}$ with $f^i = \epsilon$ if $i > m$ is even and $x^i = 1$ if $i > m$ is odd. Similarly, we may define subalgebras 
$H_{[m, \infty)}$ of $H_{(-\infty, \infty)}$ for any $m \in \mathbb{Z}$. Finally, consider the subalgebra $H_{(-\infty, -1]} \otimes H_{[2, \infty)}$ 
of $H_{(-\infty, \infty)}$. Note that it consists of all (finite sums of) elements $\cdots \rtimes x^{-1} \rtimes f^0 \rtimes x^1 \rtimes \cdots$
of $H_{(-\infty, \infty)}$ with $f^0 = \epsilon, x^1 = 1$. The inclusion $H_{(-\infty,-1]} \otimes H_{[2,\infty)} \subset H_{(-\infty, \infty)}$
of (infinite-dimensional) algebras is of significant importance to us as it
will be used in $\S 3$ to construct the model for the quantum double inclusion of $R^H \subset R$.

The following results will be very useful. We refer to \cite[Theorem 2.1, Corollary 2.3(ii)]{BlaMnt1985} for the proof of Lemma \ref{matrixalg},
\cite[Lemma 4.5.3]{Jjo2008} or \cite[Proposition 3]{DeKdy2015} for the proof of Lemma 
\ref{commutants} and \cite[Lemma 4.2.3]{Jjo2008} for the proof of  Lemma \ref{anti1}.
\begin{lemma}\label{matrixalg}
 $H \rtimes H^* \rtimes H \rtimes \cdots$ ($2k$-terms) is isomorphic to the matrix algebra $M_{n^k}(\mathbb{C})$ where $n = dim ~H$.
\end{lemma}

\begin{lemma}\label{commutants}
For any $p\in {\mathbb Z}$, the subalgebras $H_{(-\infty,p]}$ and $H_{[p+2,\infty)}$ are mutual commutants in $H_{(-\infty, \infty)}$.
\end{lemma} 
Given a positive integer $l$ and $X \in A(H)_l$, say $X = x^1 \rtimes f^2 \rtimes \cdots \rtimes f^{l-1} \rtimes x^l$ or 
$x^1 \rtimes f^2 \rtimes \cdots \rtimes x^{l-1} \rtimes f^l$ according as $l$ is odd or even, let $X^{\prime}$ denote 
the element defined by 
\begin{align*}
 X^{\prime} = Sx^l \rtimes Sf^{l-1} \rtimes \cdots \rtimes Sf^2 \rtimes Sx^1 \ \mbox{or} \ Sf^l \rtimes Sx^{l-1} \rtimes \cdots \rtimes Sf^2 
 \rtimes Sx^1
\end{align*}
according as $l$ is odd or even. It is evident that $X^{\prime} \in A(H)_l$ or $A(H^*)_l$ according as $l$ is odd or even. Obviously, 
the map $X \mapsto X^{\prime}$ is a linear isomorphism 
of $A(H)_l$ onto $A(H)_l$ or $A(H^*)_l$ according as $l$ is odd or even. Similarly, if $X \in A(H^*)_l$, say $X = f \rtimes x
\rtimes g \rtimes \cdots$ and if $X^{\prime}$ denotes the element given by
\begin{align*}
X^{\prime} = \cdots \rtimes Sg \rtimes Sx \rtimes Sf, 
\end{align*}
then certainly $X^{\prime}
\in A(H^*)_l$ or $A(H)_l$ according as $l$ is odd or even and the map $X \mapsto X^{\prime}$ is a linear isomorphism 
of $A(H^*)_l$ onto $A(H^*)_l$ or $A(H)_l$ according as $l$ is odd or even.
\begin{lemma}\label{anti1}
 The map $X \mapsto X^{\prime}$ is a $*$-anti-isomorphism.
\end{lemma}
We now need to recall the Fourier transform map for $H$. The Fourier transform map $F_H : H \rightarrow H^*$ is defined by 
$F_H(a) = \delta \phi_1(a)\phi_2$ and
satisfies $F_{H^*} F_H = S$. We will usually omit the subscript of $F_H$ and $F_{H^*}$ and write
both as $F$ with the argument making it clear which is meant.

\section{Subfactor planar algebras and planar algebras of Kac algebras}

\subsection{Subfactor planar algebras}

The notion of planar algebras was introduced in \cite{Jns1999}. For the basics of (subfactor) planar algebras, we refer to \cite{Jns1999},
\cite{KdySnd2004} and \cite{KdyLndSnd2003}. We will use the older notion of planar algebras where $Col$, the set of colours, is given by 
$\{(0, \pm),1,2,\cdots\}$ (note that only $0$ has two variants, namely, $(0, +)$ and $(0,-)$).
This is equivalent to the 
newer notion of planar algebras where $Col = \{(k, \pm) : k \geq 0 \ \mbox{integer}\}$ and we refer to 
\cite[Proposition 1]{DeKdy2016} for the proof of this equivalence.
We will use the notation $T^{k_0}_{k_1, k_2, \cdots, k_b}$ to denote a tangle $T$ of colour $k_0$ (i.e., the colour of the external
box of $T$ is $k_0$) with $b$ internal
 boxes ($b$ may be zero also) such that the colour of the $i$-th internal box is $k_i$. Given a tangle $ T = T^{k_0}_{k_1, k_2, \cdots, k_b}$ and 
 a planar algebra $P$, $Z^P_T$ will always denote the associated linear map from $P_{k_1} \otimes P_{k_2} \otimes \cdots \otimes P_{k_b}$ to $P_{k_0}$ 
 induced by the tangle $T$.
 
\begin{figure}[!h]
\begin{center}
\psfrag{k}{\Huge $k$}
\psfrag{1}{\Huge $1$}
\resizebox{6.0cm}{!}{\includegraphics{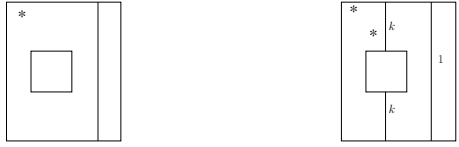}}
\end{center}
\caption{Inclusion tangles : $I_{0, +}^1, I_k^{k+1} (k \geq 1)$}
\label{fig:pic64}
\end{figure}

\begin{figure}[!h]
\begin{center}
\psfrag{a}{\Huge $D_1$}
\psfrag{b}{\Huge $D_2$}
\psfrag{k}{\Huge $k$}
\resizebox{5.0cm}{!}{\includegraphics{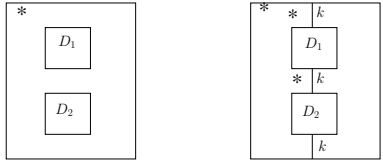}}
\end{center}
\caption{Multiplication tangle : $M_{(0, +), (0, +)}^{(0, +)}$, $M_{k, k}^k$}
\label{fig:pic65}
\end{figure}

\begin{figure}[!h]
\begin{center}
\psfrag{1}{\Huge $1$}s
\psfrag{2}{\Huge $k$}
\resizebox{6.0cm}{!}{\includegraphics{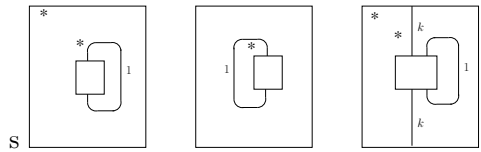}}
\end{center}
\caption{Conditional expectation Tangles : $E_1^{0, +}, E_1^{0, -}, E_{k+1}^k$}
\label{fig:pic66}
\end{figure}

\begin{figure}[!h]
\begin{center}
\psfrag{k}{\Huge $k$}
\psfrag{2}{\Huge $2$}
\psfrag{k-2}{\Huge $k-2$}
\resizebox{6.0cm}{!}{\includegraphics{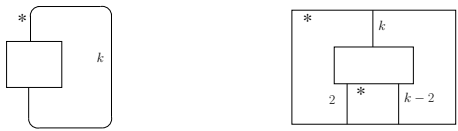}}
\end{center}
\caption{trace tangle : $tr_k^{(0, +)}$(left) and rotation tangle : $R_k^k (= R_k) (k \geq 2)$(right)}
\label{fig:pic67}
\end{figure}

\begin{figure}[!h]
\begin{center}
\psfrag{k}{\Huge $k$}
\resizebox{6.0cm}{!}{\includegraphics{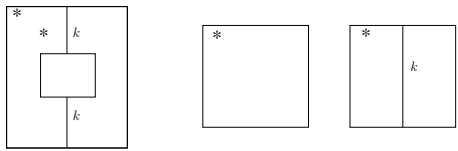}}
\end{center}
\caption{Identity tangles : $I_k^k$ (left) and Unit tangles : $1^{0,+}$ (middle), $1^k$ (right)}
\label{fig:pic70}
\end{figure}

\begin{figure}[!h]
\begin{center}
\psfrag{k}{\Huge $k$}
\psfrag{1}{\Huge $k$}
\psfrag{n-1}{\Huge $n-k$}
\resizebox{5.0cm}{!}{\includegraphics{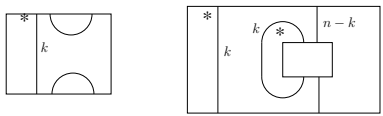}}
\end{center}
\caption{$\mathcal{E}^{k+2} (k \geq 0)$ (Jones projection tangles) (left) and $(Q(k))^n_n (1 \leq k \leq n)$(right)}
\label{fig:pic72}
\end{figure}

\begin{figure}[!h]
\begin{center}
\psfrag{1}{\Huge $1$}
\psfrag{2}{\Huge $2$}
\psfrag{3}{\Huge $3$}
\resizebox{6.0cm}{!}{\includegraphics{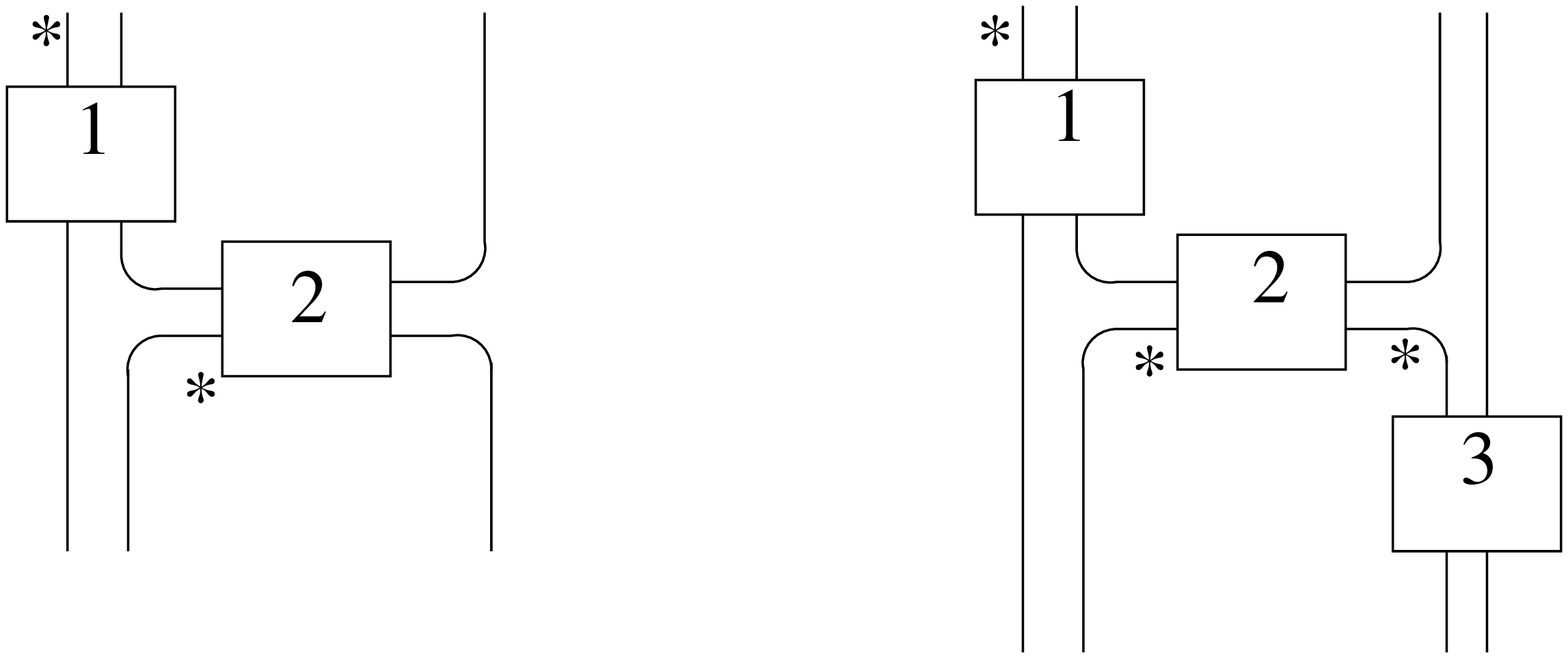}}
\end{center}
\caption{The tangles $T^3$(left) and $T^4$(right)}
\label{fig:pic59}
\end{figure}

 In Figures \ref{fig:pic64} - \ref{fig:pic59} we show and describe several tangles that will be useful to us in the sequel.
 Observe that Figure \ref{fig:pic59} shows some elements of a family of tangles. In Figure
\ref{fig:pic59} we have the tangles $T^n$ of colour $n$ for $n \geq 2$, with exactly $n-1$ internal $2$-boxes and no internal regions 
illustrated for $n = 3$ and $n = 4$. 

We will also find it useful to recall the notions of cabling and adjoint for tangles and for planar algebras.
Given any positive integer $m$ and a tangle $T$, say $T = T_{k_1, k_2, \cdots, k_b}^{k_0}$, 
the $m$-cabling of $T$, denoted by $T^{(m)}$, is the tangle obtained from $T$ by replacing each string of $T$ by a parallel cable of $m$-strings.
It is worth noting that the number of internal boxes of $T^{(m)}$ and $T$ are the same and that if $k_i(T^{(m)})$ denotes the colour of the 
$i$-th internal disc of $T^{(m)}$, then
\begin{align*}
        k_i(T^{(m)}) = 
        \begin{cases}
        mk_i, & \mbox{if} \ k_i > 0 \\
        (0, +), & \mbox{if} \ k_i = (0, +) \\
        (0, -), & \mbox{if} \ k_i = (0, -) \ \mbox{and} \  m \  \mbox{is odd}\\
        (0, +), & \mbox{if} \ k_i = (0, -) \ \mbox{and} \  m \ \mbox{is even}.\\
        \end{cases}
        \end{align*}
        Now given any planar algebra $P$, construct a new planar algebra $^{(m)}\!P$, called $m$-cabling of $P$, by setting
        \begin{align*}
        ^{(m)}\!P_k = 
        \begin{cases}
        P_{mk}, & \mbox{if} \ k > 0 \\
        P_{(0, +)}, & \mbox{if} \ k = (0, +) \\
        P_{(0, -)}, & \mbox{if} \ k = (0, -) \ \mbox{and} \  m \  \mbox{is odd}\\
        P_{(0, +)}, & \mbox{if} \ k = (0, -) \ \mbox{and} \  m \ \mbox{is even}\\
        \end{cases}
        \end{align*}
       and defining $Z_T^{^{(m)}\!P} = Z_{T^{(m)}}^P$ for any tangle $T$.  
Similarly,  given a planar algebra $P$, we construct a new planar algebra ${^*\!P}$, called the adjoint of $P$, where for 
      any $k \in Col, {(^*\!P)_k} = P_k$ as vector 
      spaces and given any tangle $T$, the action $Z_T^{^*\!P}$ of $T$ on $^*\!P$ is specified by $Z_{T^*}^P$ where $T^*$ is the tangle obtained by
      reflecting the tangle $T$ across any line in the plane.
 
 The following theorem on generating tangles will be useful.
\begin{theorem}\label{generating}\cite[Theorem 3.5]{KdySnd2004}
Let $\mathcal{T}$ denote the set of all tangles, and suppose $\mathcal{T}_0$ is a
subclass of $\mathcal{T}$ which satisfies:
\begin{itemize}
 \item[(a)] $\{1^{0, +}, 1^{0, -}\} \cup \{R_k : k\geq 2\} \cup \{E^k_{k+1}, M^k_{k, k}, I^{k+1}_k : k \in Col\} \subset \mathcal{T}_0$; and
 \item[(b)] $\mathcal{T}_0$ is closed under composition, when it makes sense.
\end{itemize}
Then, $\mathcal{T} = \mathcal{T}_0$. 
\end{theorem}

 Among planar algebras, the ones that we will be interested in are the subfactor
 planar algebras. If $P$ is a subfactor planar algebra of modulus $d$, then for each $k \geq 0$, we refer to the (faithful, positive, normalised) trace 
 $\tau : P_k \rightarrow \mathbb{C}$ defined for 
 $x \in P_k$ by 
 $\tau(x) = d^{-k} Z_{tr_k^{(0, +)}}(x)$ as the normalised pictorial trace on $P_k$. 
 
 The following fundamental theorem due to Jones \cite{Jns1999} relates subfactors and subfactor planar
algebras.
\begin{theorem}\label{jones}
 Let $N ( = M_0) \subset M ( = M_1) \subset^{e_2} M_2 \subset^{e_3} M_3 \subset \cdots$
be the tower of the basic construction associated to an extremal subfactor with $[M:N] = d^2 < \infty$, where, of course, 
$M_{n+2} = \left< M_{n+1}, e_{n+2} \right>$ ($n \geq 0$) is the result of basic construction applied to the initial inclusion $M_n \subset M_{n+1}$.
 Then there exists a unique subfactor planar algebra $P = P^{N \subset M}$ of
modulus $d$ satisfying the following conditions:
\begin{itemize}
\item[(i)] $P_k = N^{\prime} \cap M_k$ for all $k \geq 1$ - where this is regarded as an equality of $*$-algebras which is 
consistent with the inclusions on the two sides;
\item[(ii)] $Z_{\mathcal{E}^k}^P (1) = d e_k$ for all $k \geq 2$;
\item[(iii)] $Z_{(E^{\prime})^k_k}(x) = d E_{M^{\prime} \cap M_k}(x)$, for all $x \in N^{\prime} \cap M_k \ (k \geq 1)$ where $(E^{\prime})^n_n = Q(1)^n_n$;
\item[(iv)] $Z_{E^k_{k+1}}(x) = d E_{N^{\prime} \cap M_k}(x)$ for all $x \in N^{\prime} \cap M_{k+1}$
and this is required to hold for all $k$ in $Col$ where for $k = (0, \pm)$, the equation is interpreted as
\begin{align*}
 Z_{E_1^{(0, \pm)}}(x) = d tr_M(x), \forall x \in N^{\prime} \cap M.
\end{align*}

\end{itemize}

Conversely, any subfactor planar algebra $P$ with modulus $d$ arises from an extremal
subfactor of index $d^2$ in this fashion.
\end{theorem}
 \begin{remark}\label{a}
 It is a consequence of Theorem \ref{jones} that  
  for $1 \leq k \leq n$, $P_{k, n} \stackrel{\text{def}}{=} ran (Z_{Q(k)^n_n})$ is equal to $M^{\prime}_k \cap M_n$
  where $Q(k)^n_n$ is the tangle as shown on the right in Figure \ref{fig:pic72}.
 \end{remark}

 \subsection{Planar algebra associated to a Kac algebra}
 Suppose that $H$ acts outerly on the hyperfinite $II_1$ factor $M$. Let $P(H, \delta)$ (or, simply, $P(H)$) denote the subfactor planar algebra associated
to $M^H \subset M$ where $M^H$ is the fixed-point subalgebra of $M$. The following theorem (which is a reformulation of Theorem 5.1 of
\cite{KdyLndSnd2003}) gives a presentation of the planar algebra $P(H)$.
We now make a brief digression concerning notation. Given a label set $L = \coprod_{k \in col} L_k$, and a subset $R$ of the universal planar algebra $P(L)$ defined on the label set $L$, the 
notation $P(L, R)$ is used to denote the quotient of $P(L)$ by the planar ideal generated by the subset $R$ of $P(L)$.
\begin{theorem}\cite[Theorem 5.1]{KdyLndSnd2003}
 There is a $*$-isomorphism $P(L, R) \rightarrow P(H)$ of planar algebras where
\begin{align*}
        L_k = 
        \begin{cases}
        H, & \mbox{if} \ k = 2 \\
        \emptyset, & \mbox{otherwise}         
        \end{cases}
        \end{align*}
and $R$ being given by the set of relations in Figures \ref{fig:LnrMdl} - \ref{fig:XchNtp} (where (i) we write the relations
as identities - so the statement $a = b$ is interpreted as $a - b \in R;$ (ii) $\zeta \in k$  and  $a,b \in H;$ and
(iii) the external boxes of all tangles appearing in the relations are left 
undrawn and it is assumed that all external $*$-arcs are the leftmost arcs).
\begin{figure}[!h]
\begin{center}
\psfrag{zab}{\huge $\zeta a + b$}
\psfrag{eq}{\huge $=$}
\psfrag{a}{\huge $a$}
\psfrag{b}{\huge $b$}
\psfrag{z}{\huge $\zeta$}
\psfrag{+}{\huge $+$}
\psfrag{del}{\huge $\delta$}
\resizebox{11.0cm}{!}{\includegraphics{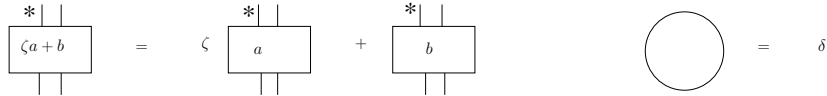}}
\end{center}
\caption{The L(inearity) and M(odulus) relations}
\label{fig:LnrMdl}
\end{figure}

\begin{figure}[!h]
\begin{center}
\psfrag{zab}{\huge $\zeta a + b$}
\psfrag{eq}{\huge $=$}
\psfrag{1h}{\huge $1_H$}
\psfrag{h}{\huge $h$}
\psfrag{z}{\huge $\zeta$}
\psfrag{+}{\huge $+$}
\psfrag{del}{\huge $\delta^{-1}$}
\resizebox{9.0cm}{!}{\includegraphics{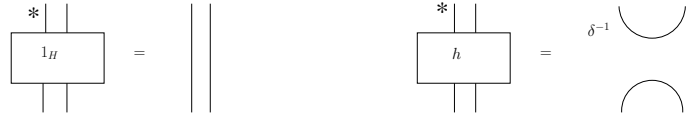}}
\end{center}
\caption{The U(nit) and I(ntegral) relations}
\label{fig:NitNtg}
\end{figure}

\begin{figure}[!h]
\begin{center}
\psfrag{epa}{\huge $\epsilon(a)$}
\psfrag{eq}{\huge $=$}
\psfrag{delinphia}{\huge $\delta \phi(a)$}
\psfrag{h}{\huge $h$}
\psfrag{a}{\huge $a$}
\psfrag{+}{\huge $+$}
\psfrag{del}{\huge $\delta^{-1}$}
\resizebox{9.0cm}{!}{\includegraphics{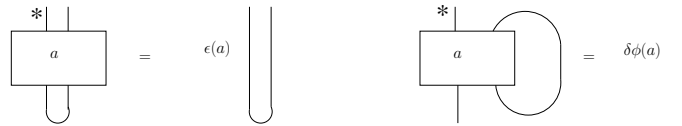}}
\end{center}
\caption{The C(ounit) and T(race) relations}
\label{fig:CntTrc}
\end{figure}

\begin{figure}[!h]
\begin{center}
\psfrag{epa}{\huge $\epsilon(a)$}
\psfrag{eq}{\huge $=$}
\psfrag{a1}{\huge $a_1$}
\psfrag{a2b}{\huge $a_2\,b$}
\psfrag{b}{\huge $b$}
\psfrag{a}{\huge $a$}
\psfrag{sa}{\huge $Sa$}
\psfrag{del}{\huge $\delta$}
\resizebox{10.0cm}{!}{\includegraphics{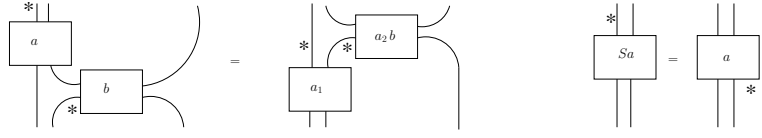}}
\end{center}
\caption{The E(xchange) and A(ntipode) relations}
\label{fig:XchNtp}
\end{figure}
\end{theorem}
In these figures, note that the shading is such that all the 2-boxes that occur are of colour $2$. Also note that the modulus relation is a pair
of relations - one for each choice of shading the circle.

Next, we recall (a reformulation of) a result from \cite{KdySnd2006}.
Let ${\mathcal T}(k, p) (p \leq k-1)$ denote the set of $k$ tangles (interpreted as $0$ for $k=0$) 
with $p$ internal boxes of colour $2$ and no `internal regions'. If $p = k-1$, we will simply write 
${\mathcal T}(k)$ instead of ${\mathcal T}(k, k-1)$. The result then asserts:
\begin{lemma}\label{iso1}
[reformulation of Lemma 16 of \cite{KdySnd2006}]
For each tangle $X \in {\mathcal T}(k, p)$, the map $Z_X^{P(H)}:
(P(H)_2)^{\otimes p} \rightarrow P(H)_k$ is an injective linear map and if $p = k-1$, then $Z_X^{P(H)}:
(P(H)_2)^{\otimes k-1} \rightarrow P(H)_k$ is a linear isomorphism.
\end{lemma}

The following lemma (a reformulation of \cite[Proposition 4.3.1]{Jjo2008}) establishes algebra isomorphisms between $P(H)_k$ and finite iterated 
crossed product algebras.
\begin{lemma}\label{iso}
  For each $k \geq 2$, the maps
  \begin{align*}
 A(H)_{k-1} \rightarrow P(H)_k \ \mbox{given by} \ x \rtimes f \rtimes y \rtimes \cdots (k-1 \ \mbox{terms}) \mapsto 
  Z^{P(H)}_{T^k}(x \otimes Ff \otimes y \otimes \cdots)
  \end{align*}
and 
\begin{align*}
  A(H^*)_{k-1} \rightarrow P(H^*)_k \ \mbox{given by} \ f \rtimes x \rtimes g \rtimes \cdots (k-1 \ \mbox{terms}) \mapsto 
  Z^{P(H^*)}_{T^k}(f \otimes Fx \otimes g \otimes \cdots)
  \end{align*}
 are $*$-algebra isomorphisms. 
\end{lemma}
We will use this identification of $A(H)_{k-1}$ (resp., $A(H^*)_{k-1}$) with $P(H)_k$ (resp., $P(H^*)_k$) very frequently
without mention. If $i \leq j$, let $tr_{H_{[i, j]}}$ denote the faithful, positive, tracial state on $H_{[i, j]}$ obtained by pulling back
the normalised pictorial trace on $P(H^i)_{j-i+2}$ using the algebra isomorphism of Lemma \ref{iso}.
It is easy to see that $tr_{H_{[i, j]}}$ indeed is the linear
functional on $H_{[i, j]}$ given by 
\begin{align*}
       \begin{cases}
        h^i \otimes \phi^{i+1} \otimes h^{i+2} \otimes \cdots (j-i+1 \mbox{-terms}), & \mbox{if} \  i \ \mbox{is even}\\
        \phi^i \otimes h^{i+1} \otimes \phi^{i+2} \otimes \cdots (j-i+1 \mbox{-terms}), & \mbox{if} \  i \  \mbox{is odd}.         
        \end{cases}
\end{align*}
Thus, for instance, if we assume $i$ to be odd, $j$ to be even 
and if $X \in H_{[i, j]}$, say,
$X = x^i \rtimes f^{i+1} \rtimes \cdots \rtimes x^{j-1} \rtimes f^j$,
then $tr_{H_{[i, j]}}(X) = \phi^i(x^i) f^{i+1}(h^{i+1}) \cdots \phi^{j-1}(x^{j-1}) f^j(h^j)$.

 \section{A model for the quantum double inclusion of $R^H \subset R$}
 The notion of quantum double inclusion associated to a finite-index and finite-depth subfactor has already been defined in the introduction of this article.
 The main goal of this section is to construct a model for the quantum double inclusion associated to the Kac algebra subfactor
 $R^{H} \subset R$ where $H$ is a finite-dimensional Kac algebra acting outerly on the hyperfinite $II_1$ factor $R$ and $R^{H}$ is the fixed point 
 subalgebra under this action. Our construction of the model for the quantum double inclusion of $R^H \subset R$ closely
 follows the construction of the model for the asymptotic inclusion of $R^H \subset R$ as given in \cite{Jjo2008}.
 
 We begin with recalling from \cite{KdySnd2008} the notion of finite pre-von Neumann algebras. 
By a finite pre-von Neumann algebra, we will mean a pair $(A, \tau)$ consisting of a complex $*$-algebra $A$ that is equipped with a normalised trace $\tau$ such that 
(i) the sesquilinear form defined by $\left< a, b \right> = \tau(b^* a)$ defines an inner-product on $A$ and such that (ii) for each $a \in A$, the
 left-multiplication map $\lambda_A(a) : A \longrightarrow A$ is bounded for the trace induced norm of $A$.
 By a compatible pair of finite pre-von Neumann algebras, we will mean
a pair $(A, \tau_A )$ and $(B, \tau_B)$ of finite pre-von Neumann algebras such that $A \subseteq B$ and ${\tau_B|}_{A} = \tau_A$.

If $A$ is a finite pre-von Neumann algebra with trace $\tau_A$, the symbol $L^2(A)$ will always denote
the Hilbert space completion of $A$ for the associated norm. Obviously, the left regular
representation $\lambda_A : A \rightarrow \mathcal{L}(L^2(A))$ is well-defined, i.e., for each $a \in A, \lambda_A(a) : A \rightarrow A$
extends to a bounded operator on $L^2(A)$. The notation $A^{\prime \prime}$ will always denote the von Neumann algebra 
$(\lambda_A(A))^{\prime \prime} \subset \mathcal{L}(L^2(A))$. The following lemma (a reformulation of \cite[Proposition 4.6(1)]{KdySnd2008})
 will be of great use in the sequel. 
\begin{lemma}\label{pre}\cite[Proposition 4.6(1)]{KdySnd2008} 
Let $(A, \tau_A )$ and $(B, \tau_B)$ be a compatible pair of finite pre-von Neumann algebras. The inclusion $A \subseteq B$
extends uniquely to a normal inclusion of $A^{\prime \prime}$ into $B^{\prime \prime}$ with image $(\lambda_B(A))^{\prime \prime}$. 
\end{lemma}

  Let $A_0 \subseteq A_1$ be a unital \textit{connected} inclusion of finite-dimensional $C^*$-algberas and let $A_0 \subseteq A_1 \subseteq A_2 \subseteq 
  \cdots$ be the Jones' basic construction tower of $A_0 \subseteq A_1$. For each $n \geq 1$, let $\tau_n$ denote the unique trace on $A_n$ which is
   a Markov trace (see \cite[$\S 2.7$]{GHJ1989} for the notion of Markov trace) for the inclusion $A_{n-1} \subseteq A_n$. Set $\tau_0 =$ trace on $A_0$ (= ${\tau_1|}_{A_0}$). 
   Clearly, $A = \cup_{n = 0}^{\infty}A_n$ comes equipped with a tracial state $\tau$
   (whose restriction to $A_n$ is $\tau_n$) making $A$ into a finite pre-von Neumann algebra; in fact this is the unique tracial state on $A$
   and consequently, $A^{\prime \prime}$ 
   turns out to be a (hyperfinite) $II_1$ factor.

   The next proposition shows that suitably compatible actions of a finite-dimensional
   Kac algebra on each $A_n$ may be extended to an action on $A^{\prime \prime}$. This is proved
   in \cite{Y1993} using the fundamental unitary or Kac-Takesaki operator associated to a Kac algebra. We give another proof.

   \begin{proposition}\label{extn}
   With notations as above, let $H$ be a finite-dimensional
  Kac algebra acting on each $A_n$ and let $\alpha^n$ denote the action of $H$ on $A_n$. Assume that $\tau_n \circ \alpha^n_h = \tau_n$ 
  for all $n \geq 0$ where, as usual, $h$ denotes the unique non-zero idempotent integral of $H$. Suppose further that for any $n \geq 0$ and any $x \in H$, the following diagram commutes.
   \begin{center}
\begin{tikzpicture}
\node at (0,0) {$A_{n+1}$};
 \draw [thick,->] (.5,0) -- (1.5,0);
 \node at (2,0) {$A_{n+1}$}; 
 \node at (1,.25){$\alpha_x^{n+1}$};
 \node at (0,-1.5) {$A_n$};
 \node at (2,-1.5) {$A_n$};
 \draw [thick,->] (.5,-1.5) -- (1.5,-1.5);
 \node at (1,-1.75) {$\alpha_x^n$}; 
 \node [rotate=90] at (0,-.75) {$\subseteq$};
 \node [rotate=90] at (2,-.75) {$\subseteq$};
 \end{tikzpicture}
 \end{center}  
  Clearly, commutativity of the preceding diagram yields an action $\alpha$ of $H$ on $A$. This action of $H$ on $A$ can be extended to $A^{\prime \prime}$.
  \end{proposition}
  Before we proceed to prove Proposition \ref{extn} we will need the following result whose proof is similar to that of \cite[Lemma 4.4(b)]{KdyLndSnd2003}
  and hence, we omit its proof.
  \begin{lemma}\label{Lnd}
   Let $A$ be a finite-dimensional $C^*$-algebra equipped with a faithful tracial state $\tau$ and suppose that $H$ is a finite-dimensional 
   Kac algebra acting on $A$ satisfying $\tau(h.a) = \tau(a)$ for all $a \in A$. Then there exists a unique morphism of $C^*$-algebras 
   $H \ni x \rightarrow L_x \in \mathcal{L}(L^2(A))$ such that $L_x(a) = x.a \ \forall a \in A$ where $L^2(A)$ denotes the Hilbert space 
   completion of the inner product space $A$ equipped with the inner product induced by $\tau$.
  \end{lemma}
  We now return to the proof of Proposition \ref{extn}. 
  \begin{proof}[Proof of Proposition \ref{extn}]
  Let $\lambda : A \rightarrow \mathcal{L}(L^2(A))$ denote the left regular representation. For any $a \in A$, let the notation $\lambda_a$ 
  stand for $\lambda(a)$. Let $a^{\prime \prime}$ be an element of $A^{\prime \prime}$. Then there exists a net $(a_j)$ in $A$ such that 
   $\lambda_{a_j}$ converges in strong operator topology (SOT) to $a^{\prime \prime}$. 
   Given $x \in H$, we show that the net $(\lambda_{\alpha_x(a_j)})$ converges in SOT. To this end it suffices to see that the net
   $(\lambda_{\alpha_x(a_j)})$ is bounded and converges pointwise on $A$ (i.e., on a dense subspace of $L^2(A)$). 
   Let $\Delta(x) = x_1 \otimes x_2 = \sum_i y^i \otimes z^i$ (finite sum). Then note that for any $a \in A$, 
   \begin{align}\label{eqn}
    \lambda_{\alpha_x(a_j)}(a) = \alpha_x(a_j)a = \alpha_{x_1}(a_j) \alpha_{x_2}(\alpha_{Sx_3}(a)) = \alpha_{x_1}(a_j \alpha_{Sx_2}(a)) = 
    \sum_i \alpha_{y^i}(a_j \alpha_{Sz^i}(a))
   \end{align}
   and consequently,   
   \begin{align*}
   \lVert \lambda_{\alpha_x(a_j)}(a) \rVert \ \leq \ \sum_i \lVert \alpha_{y^i}(a_j \alpha_{Sz^i}(a)) \rVert \ \leq \ \sum_i \ \lVert y^i \rVert
    \ \lVert a_j \rVert \ \lVert Sz^i \rVert \ \lVert a \rVert \ \leq \ M \lVert a \rVert (\sum_i \lVert y^i \rVert \ \lVert Sz^i \rVert)
    \end{align*}
 where $M = \underset{j}{Sup} \ \lVert a_j \rVert (< \infty)$ and the
      second inequality follows by an appeal to Lemma \ref{Lnd}. This shows that $\underset{j}{Sup} \ \lVert \lambda_{\alpha_x(a_j)} \rVert \leq M
      (\sum_i \lVert y^i \rVert \ \lVert Sz^i \rVert) < \infty$ and hence, $(\lambda_{\alpha_x(a_j)})$ is bounded. Note that
    \begin{align*}
     (\lambda_{\alpha_x(a_j)} - \lambda_{\alpha_x(a_k)})(a) = (\alpha_x(a_j - a_k))(a)
     = \alpha_{x_1}((a_j - a_k) \alpha_{Sx_2}(a))
    \end{align*} where the second equality is a consequence of \eqref{eqn}. Therefore, $(\lambda_{\alpha_x(a_j)}(a))$ is Cauchy and hence, 
    converges in $L^2(A)$. Thus we have shown that $(\lambda_{\alpha_x(a_j)})$ converges in SOT. We define 
    $\alpha_x(a^{\prime \prime}) (\in A^{\prime \prime})$ to be the SOT-limit of $(\lambda_{\alpha_x(a_j)})$. Thus, we have a linear map 
    $\alpha : H \rightarrow End_{\mathbb{C}}(A^{\prime \prime})$ carrying $x \in H$ to $\alpha_x$. It is not hard to see that $\alpha$ satisfies conditions (i)-(iv) of Definition
    \ref{def}. To verify condition (v), we note that given $a^{\prime \prime} \in A^{\prime \prime}$, write $a^{\prime \prime} = 
    b^{\prime \prime} + i c^{\prime \prime}$ with $b^{\prime \prime}, c^{\prime \prime}$ being self-adjoint elements of $A^{\prime \prime}$.
    By Kaplansky's density theorem, there exists nets $(b_j)_{j \in J}, (c_k)_{k \in K}$ of self-adjoint elements in $A$ such that 
    $(b_j)$ converges to $b^{\prime \prime}$ in SOT and $(c_k)$ converges to $c^{\prime \prime}$ in SOT. Consider the net $(d_{(j, k) \in I \times J})$
    such that $d_{(j, k)} = b_j +i c_k$. Clearly, $d_{(j, k)}$ converges in SOT to $a^{\prime \prime}$ as well as
    $d_{(j, k)}^*$ converges in SOT to $(a^{\prime \prime})^*$. Consequently, 
    \begin{align*}
    \alpha_{Sx^*}((a^{\prime \prime})^*) = 
    \mbox{SOT-lim} \ \lambda_{(\alpha_{Sx^*}(d_{(j, k)}^*))} = \mbox{SOT-lim} \ \lambda_{(\alpha_x(d_{(j, k)}))^*} = (\mbox{SOT-lim} \ \lambda_{(\alpha_x(d_{(j, k)}))})^* = 
    (\alpha_x(a^{\prime \prime}))^*.
    \end{align*}
    Thus condition (v) of Definition \ref{def} is satisfied and hence, $\alpha$ is an action of $H$ on $A^{\prime \prime}$.
     \end{proof}
     
     Recall that the algebra $H_{(-\infty, \infty)}$ was introduced in $\S 1$. Observe that $H_{(-\infty, \infty)}$ is the 
     union of the increasing chain of $C^*$-algebras given by $\mathbb{C} \subset H_{[-1, 1]} \subset H_{[-2, 2]} \subset H_{[-3, 3]} \subset \cdots$ 
     and it is evident that for any positive integer $k$, $tr_{H_{[-(k+1), k+1]}}$ restricts  
     $tr_{H_{[-k, k]}}$. Consequently, there is a consistently defined trace on $H_{(-\infty, \infty)}$ making it into a finite pre-von Neumann
     algebra. Let $H^{\prime \prime}_{(-\infty, \infty)}$ denote the von Neumann algebra $(H_{(-\infty, \infty)})^{\prime \prime}$.
     Later we will see that $H^{\prime \prime}_{(-\infty, \infty)}$ is indeed a hyperfinite $II_1$ factor.
     
      It is easy to see that for any integer $k$,
     $\mathbb{C} \subset H_k \subset H_{[k-1, k]} \subset H_{[k-2, k]} \subset H_{[k-3, k]} \subset \cdots$ is 
     the basic construction tower associated to the initial (connected) inclusion $\mathbb{C} \subset H_k$. Consequently, 
     it follows immediately from the discussion contained in the paragraph preceding Proposition \ref{extn} that 
     \begin{align*}
     H_{(-\infty, k]}^{\prime \prime} := (\cup_{i=0}^{\infty} H_{[k-i, k]})^{\prime \prime} = (H_{(-\infty, k]})^{\prime \prime}
     \end{align*} is a 
     hyperfinite $II_1$ factor .
     Further, there is an obvious action of $H^{k+1}$ ( = $H$ or $H^*$ according as $k$ is even or odd)
     on each $H_{[n, k]}$ ($n \leq k$)  and it is easy to check that all the hypotheses of Proposition \ref{extn} are satisfied. Therefore, by an appeal
     to Proposition \ref{extn} we conclude that this action of $H^{k+1}$ on $H_{(-\infty, k]}$ extends to $H_{(-\infty, k]}^{\prime \prime}$. 
     It is thus natural to ask the fixed point subalgebra of $H_{(-\infty, k]}^{\prime \prime}$ under this action of $H^{k+1}$ and 
     whether this action of $H^{k+1}$ on $H_{(-\infty, k]}^{\prime \prime}$ is outer.  
     Moreover, it is obvious that for any integer $k$, $H_{(-\infty, k]} \subset H_{(-\infty, \infty)}$ is a compatible pair of finite pre-von Neumann algebras
     and hence, by virtue of Lemma \ref{pre}, the inclusion $H_{(-\infty, k]} \subset H_{(-\infty, \infty)}$ extends uniquely to a unital normal
     inclusion $H_{(-\infty, k]}^{\prime \prime} \subseteq H^{\prime \prime}_{(-\infty, \infty)}$. It will be useful to know the relative commutant
     $(H_{(-\infty, k]}^{\prime \prime})^{\prime} \cap H^{\prime \prime}_{(-\infty, \infty)}$. The following Proposition, which is proved 
     in \cite{Jjo2008} (see \cite[Lemma 4.5.2, Lemma 4.5.3]{Jjo2008}), answers all these questions.
     For the sake of completeness we provide a proof of this.
     \begin{proposition}\label{loc}
   1.   For integers $k < l$,
     \begin{align*}
       (H_{(-\infty, k]}^{\prime \prime})^{\prime} \cap H_{(-\infty, l]}^{\prime \prime} = 
       \begin{cases}
        \mathbb{C}, & \mbox{if} \ l = k+1 \\
        H_{[k+2, l]}, & \mbox{if} \ l \geq k+2.          
        \end{cases}
        \end{align*}
         2. For any integer $k$,  $(H_{(-\infty, k]}^{\prime \prime})^{\prime} \cap H_{(-\infty, \infty)}^{\prime \prime} = H_{[k+2, \infty)}^{\prime \prime}$.
     \end{proposition}
     \begin{proof}
     1. Given integers $k < l$, it is easy to see that the square of finite-dimensional $C^*$-algebras as shown in Figure \ref{fig:comm} is a symmetric commuting 
     square. Further, all the inclusions are connected since the lower left corner is $\mathbb{C}$ and by an appeal to Lemma \ref{matrixalg} 
     we have that
     the lower right corner $H_{[k-1, k]}$ is a matrix algebra, the upper left corner $H_{[k, l]}$ is a matrix algebra if $l-k$ is odd while the 
     upper right corner $H_{[k-1, l]}$ is a matrix algebra if $l-k$ is even. 
      \begin{figure}[!h]
\begin{center}
\psfrag{a}{\LARGE $H_{[k, l]}$}
\psfrag{b}{\Large $\subset$}
\psfrag{c}{\LARGE $H_{[k-1, l]}$}
\psfrag{d}{\Large $\cup$}
\psfrag{e}{\Large $\cup$} 
\psfrag{f}{\LARGE $\mathbb{C}$}
\psfrag{g}{\Large $\subset$}
\psfrag{h}{\LARGE $H_{[k-1, k]}$}
\resizebox{3.0cm}{!}{\includegraphics{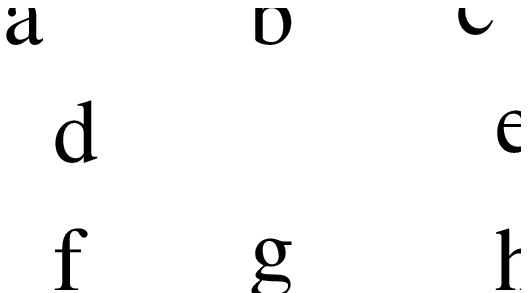}}
\end{center}
\caption{Commuting square}
\label{fig:comm}
\end{figure} 
We have already mentioned that for any integer $m$, $\mathbb{C} \subset H_m \subset H_{[m-1, m]} \subset H_{[m-2, m]} \subset H_{[m-3, m]} 
\subset \cdots$ is the basic construction 
tower of $\mathbb{C} \subset H_m$. From this it follows that $H_{[k, l]} \subset H_{[k-1, l]} \subset H_{[k-2, l]} \subset H_{[k-3, l]} 
\subset \cdots$ is the basic construction tower of $H_{[k, l]} \subset H_{[k-1, l]}$ with $(\cup_{i = 0}^{\infty}H_{[k-i, l]})^{\prime \prime} = 
H_{(-\infty, l]}^{\prime \prime}$. Again it follows, by an appeal to \cite[Proposition 4.3.6]{JnsSnd1997}, that $\mathbb{C} \subset 
H_{[k-1, k]} \subset H_{[k-3, k]} \subset H_{[k-5, k]} \subset \cdots$ 
is the basic construction tower of $\mathbb{C} \subset H_{[k-1, k]}$ with $(\cup_{i = 1}^{\infty}H_{[k-2i+1, k]})^{\prime \prime} = 
H_{(-\infty, k]}^{\prime \prime}$. 
Hence, by Ocneanu's compactness theorem (see \cite[Theorem 5.7.6]{JnsSnd1997}) we conclude that $(H_{(-\infty, k]}^{\prime \prime})^{\prime}
\cap H_{(-\infty, l]}^{\prime \prime}$ equals $(H_{[k-1, k]})^{\prime} \cap H_{[k, l]}$ which, by an appeal to Lemma \ref{commutants}, equals $\mathbb{C}$ or 
$H_{[k+2, l]}$ according as $ l = k+1$ or $l \geq k+2$.

2. Let $k$ be a fixed integer. For each $m \geq 2$, let $E_m$ denote the conditional expectation of $H_{(-\infty, \infty)}^{\prime \prime}$ onto
 $H_{(-\infty, k+m]}^{\prime \prime}$.  
 Given $x \in (H_{(-\infty, k]}^{\prime \prime})^{\prime} \cap H_{(-\infty, \infty)}^{\prime \prime}$, if we set 
 $x_m = E_m(x)$ ($m \geq 2$), then it is not hard to see that $x_m$ converges in weak operator topology to $x$. 
 Further note that for any $y \in (H_{(-\infty, k]}^{\prime \prime})$,
 $y x_m = y E_m(x) = E_m(y x) = E_m(x y) = E_m(x) y = x_m y$ so that 
 $x_m \in (H_{(-\infty, k]}^{\prime \prime})^{\prime} \cap H_{(-\infty, k+m]}^{\prime \prime}$ for all $m \geq 2$. Consequently,
 by an appeal to Lemma \ref{loc}(1) we conclude that $x_m \in H_{[k+2, k+m]} \subset H_{[k+2, \infty)}^{\prime \prime}$ for all $m \geq 2$.
 Therefore, $x \in H_{[k+2, \infty)}^{\prime \prime}$. The reverse inclusion is trivial. 
      \end{proof}
     It was shown in \cite{Jjo2008} (and is easy to see) that:
      \begin{lemma}\label{inv}
     $(H_{(-\infty, k]}^{\prime \prime})^{H^{k+1}} = H_{(-\infty, k-1]}^{\prime \prime}, \ H_{(-\infty, k]}^{\prime \prime} \rtimes H^{k+1} = 
     H_{(-\infty, k+1]}^{\prime \prime}$. 
     \end{lemma}
     It follows immediately from the preceding lemma together with an appeal to Lemma \ref{loc} that for any integer $k$,   
     $((H_{(-\infty, k]}^{\prime \prime})^{H^{k+1}})^{\prime} \cap H_{(-\infty, k]}^{\prime \prime} = H_{(-\infty, k-1]}^{\prime \prime} \cap  
     H_{(-\infty, k]}^{\prime \prime} = \mathbb{C}$. Thus we see that:
     \begin{lemma}\label{outer}
      The action of $H^{k+1}$ on $H_{(-\infty, k]}^{\prime \prime}$ is outer. 
     \end{lemma}
     Therefore, we conclude from Lemmas \ref{inv} and \ref{outer} that: 
     \begin{lemma}\label{model}
     For any integer $k$, $H_{(-\infty, k-1]}^{\prime \prime} \subset H_{(-\infty, k]}^{\prime \prime}$ is a model for $R^{H^{k+1}} \subset R$ 
     for some outer action of $H^{k+1}$ on the hyperfinite $II_1$ factor $R$ and $H_{(-\infty, k-1]}^{\prime \prime} \subset 
     H_{(-\infty, k]}^{\prime \prime} \subset H_{(-\infty, k+1]}^{\prime \prime} \subset H_{(-\infty, k+2]}^{\prime \prime} \subset \cdots$ is a 
     model for the basic construction tower of $R^{H^{k+1}} \subset R$.
     \end{lemma}
     Since, by Lemma \ref{model},
     $H_{(-\infty, k-1]}^{\prime \prime} \subset H_{(-\infty, k]}^{\prime \prime} \subset 
     H_{(-\infty, k+1]}^{\prime \prime} \subset H_{(-\infty, k+2]}^{\prime \prime} \subset \cdots$ is the basic construction tower of 
     $H_{(-\infty, k-1]}^{\prime \prime} \subset H_{(-\infty, k]}^{\prime \prime}$, 
      $\cup_{i=-1}^{\infty}H_{(-\infty, k+i]}^{\prime \prime}$ is a finite pre-von Neumann algebra and 
     $(\cup_{i=-1}^{\infty}H_{(-\infty, k+i]}^{\prime \prime})^{\prime \prime}$ is a $II_1$ factor which is hyperfinite since each of 
     $H_{(-\infty, k-1]}^{\prime \prime}, H_{(-\infty, k]}^{\prime \prime}$ is hyperfinite. It is easy to see that
     $(\cup_{i=-1}^{\infty}H_{(-\infty, k+i]}^{\prime \prime})^{\prime \prime} = H_{(-\infty, \infty)}^{\prime \prime}$. Thus we have shown that 
     $H_{(-\infty, \infty)}^{\prime \prime}$ is a hyperfinite $II_1$ factor.

     With $k = 0$ in Lemma \ref{model}, we have that $H_{(-\infty, -1]}^{\prime \prime} \subset H_{(-\infty, 0]}^{\prime \prime}$ is a model for 
     $R^H \subset R$ and thus, a model for the \textit{quantum double inclusion} of $R^H \subset R$ is given by 
     \begin{align*}
      H_{(-\infty, -1]}^{\prime \prime} \vee ((H_{(-\infty, 0]}^{\prime \prime})^{\prime} \cap H^{\prime \prime}_{(-\infty, \infty)}) \subseteq
      H^{\prime \prime}_{(-\infty, \infty)}. 
     \end{align*}
     By an appeal to Lemma \ref{loc}(2), one can easily see that 
     \begin{align*}
     (H_{(-\infty, 0]}^{\prime \prime})^{\prime} \cap H^{\prime \prime}_{(-\infty, \infty)} = H^{\prime \prime}_{[2, \infty)}
     \end{align*}
     and consequently, 
     \begin{align*}
      H_{(-\infty, -1]}^{\prime \prime} \vee 
     ((H_{(-\infty, 0]}^{\prime \prime})^{\prime} \cap H^{\prime \prime}_{(-\infty, \infty)})  = 
     H_{(-\infty, -1]}^{\prime \prime} \vee H^{\prime \prime}_{[2, \infty)} = 
     (H_{(-\infty, -1]} \otimes H_{[2, \infty)})^{\prime \prime}.
     \end{align*}
      \begin{definition}\label{def1}
      Set $\mathcal{N} = (H_{(-\infty, -1]} \otimes H_{[2, \infty)})^{\prime \prime}$ and 
      $\mathcal{M} = H_{(-\infty, \infty)}^{\prime \prime}$.
     \end{definition}
     We have thus shown that:
     \begin{proposition}\label{qdim}
      The subfactor $\mathcal{N} \subset \mathcal{M}$ is a model for the quantum double inclusion of $R^{H} \subset R$. 
     \end{proposition}

     \section{Basic construction tower of $\mathcal{N} \subset \mathcal{M}$ and relative commutants}

The purpose of this section is to 
find the basic construction tower associated to $\mathcal{N} \subset \mathcal{M}$ and also to compute the relative commutants.
\subsection{Some finite-dimensional basic constructions.}
This subsection is devoted to analysing the basic constructions associated to certain unital inclusions of finite-dimensional $C^*$-algebras.
We begin with recalling the following lemma (a reformulation of Lemma 5.3.1 of \cite{JnsSnd1997}) which provides an abstract characterisation of 
the basic construction associated to a unital inclusion of finite-dimensional $C^*$-algebras.
\begin{lemma}\cite[Lemma 5.3.1]{JnsSnd1997} \label{basic}
 Let $A \subseteq B \subseteq C$ be a unital inclusion of finite-dimensional $C^*$-algebras. Let $tr_B$ denote a faithful tracial state on $B$ and
 let $E_A$ denote the $tr_B$-preserving conditional expectation of $B$ onto $A$. Let $f \in C$ be a projection. Then $C$ is isomorphic to the 
 basic construction for  $A \subseteq B$ with $f$ as the Jones projection if the following conditions are satisfied:
 \begin{itemize}
 \item[(i)] $f$ commutes with every element of $A$ and $a \mapsto af$ is an injective map of $A$ into $C$,
 \item[(ii)] $f$ implements the trace-preserving conditional expectation of $B$ onto $A$ i.e., $fbf = E_A(b)f$ for all $b \in B$, and
 \item[(iii)] $BfB = C$.
 \end{itemize}
 \end{lemma}
 In the next lemma, we explicitly compute certain conditional expectation maps. 
 
 \begin{lemma}\label{exp}

   (i) Given a sequence of integers $k_1 \leq k_2 < k_3 \leq k_4 < \cdots < k_{2r-1} \leq k_{2r}$ in $[p,q]$, where $r$ is any positive integer,
   then the trace-preserving conditional expectation of $H_{[p, q]}$ onto
   $H_{[k_1, k_2]} \otimes H_{[k_3, k_4]} \otimes \cdots \otimes H_{[k_{2r-1}, k_{2r}]}$ is given by the map $tr_{H_{[p, k_1-1]}} \otimes 
   Id_{H_{[k_1, k_2]}} \otimes tr_{H_{[k_2+1, k_3-1]}} \otimes Id_{H_{[k_3, k_4]}} \otimes \cdots \otimes Id_{H_{[k_{2r-1}, k_{2r}]}} \otimes
   tr_{H_{[k_{2r}+1, q]}}$, with the obvious interpretations when, say, $k_{2i} + 1 > k_{2i+1} -1$.

(ii)    Given integers $l \geq 1$ and $s \geq 0$, let $\psi_{l, s}$ denote the embedding of $H_{[-l, 2+s]}$ inside 
   $H_{[-l, -1]} \otimes H_{[2, 6+s]}$ specified as follows:
   \begin{align*}
 \mbox{Let} \ X = x^{-l}/f^{-l} \rtimes \cdots \rtimes x^{2+s}/f^{2+s} \in H_{[-l, 2+s]}, 
\end{align*}
then $\psi_{l, s}(X) \in H_{[-l, -1]} \otimes H_{[2, 6+s]}$ is given by 
\begin{align*}
 (x^{-l}/f^{-l} \rtimes \cdots \rtimes f^{-2} \rtimes x^{-1}_1) \otimes 
 (\epsilon \rtimes x^{-1}_2 \rtimes f^0 \rtimes \cdots \rtimes x^{2+s}/f^{2+s}).
\end{align*}
Then the trace-preserving conditional expectation $E$ of $H_{[-l, -1]} \otimes H_{[2, 6+s]}$ onto 
$H_{[-l, 2+s]}$ is given by
       \begin{align*}       
       & E((x^{-l}/ f^{-l} \rtimes \cdots \rtimes x^{-1}) \otimes (f^2 \rtimes x^{3} \rtimes \cdots \rtimes x^{6+s}/f^{6+s}))\\
       & = \phi(Sx_2^{-1}x^{3}) f^2(h) \ x^{-l}/f^{-l} \rtimes \cdots \rtimes f^{-2} \rtimes
       x^{-1}_1 \rtimes f^{4} \rtimes \cdots \rtimes x^{6+s}/f^{6+s}. 
       \end{align*}
       \end{lemma}   
       
       \begin{proof}
       (i) Follows easily by direct computations and is left to the reader.

      (ii)  For notational convenience, we assume $l=1$ and $s=0$. The proof for the general case will follow
       in a similar fashion.
       It is trivial to verify that $E$ is trace-preserving. To see that $E$ is $H_{[-1, 2]}$ - $H_{[-1, 2]}$ linear, consider
       \begin{align*}
       X = x^{-1} \rtimes f^0 \rtimes x^1 \rtimes f^2 \in H_{[-1, 2]} 
       \end{align*}
       whose image in $H_{[-1, -1]} \otimes 
       H_{[2, 6]}$  under $\psi_{1, 0}$ , also denoted by the same symbol $X$, is given by 
       \begin{align*}
        x^{-1}_1 \otimes (\epsilon \rtimes x^{-1}_2 \rtimes f^0 \rtimes x^1 \rtimes f^2)
       \end{align*}
       and let
       \begin{align*}
        Y = y^{-1} \otimes (g^2 \rtimes y^3 \rtimes g^4 \rtimes y^5 \rtimes g^6) \in H_{[-1, -1]} \otimes H_{[2, 6]}.
       \end{align*}
       Then,
       \begin{align*}
       & X Y  =\  g^2_2(x^{-1}_2) \ f^0_1(y^3_2) \ g^4_2(x^1_1) \ f^2_1(y^5_2) \ 
       x^{-1}_1 y^{-1} \otimes (g^2_1 \rtimes x^{-1}_3 y^3_1 \rtimes f^0_2 g^4_1 \rtimes x^1_2 y^5_1 \rtimes f^2_2 g^6)\\
       \end{align*}
       and hence,
       \begin{align*}
         E(XY) & =  \phi(S(x^{-1}_2 y^{-1}_2) x^{-1}_4 y^{3}_1) \ g^2_1(h) \ g^2_2(x^{-1}_3) \ f^0_1(y^3_2) \ g^4_2(x^1_1) \ f^2_1(y^5_2) \
          x^{-1}_1 y^{-1}_1 \rtimes f^0_2 g^4_1 \rtimes x^1_2 y^5_1 \rtimes f^2_2 g^6\\
        & = \phi(Sy^{-1}_2 y^{3}_1) \ g^2(h) \ f^0_1(y^3_2) \ g^4_2(x^1_1) f^2_1(y^5_2) \ x^{-1} y^{-1}_1 \rtimes f^0_2 g^4_1 \rtimes 
        x^1_2 y^5_1 \rtimes f^2_2 g^6. 
        \end{align*}
        On the other hand, a little computation shows that  
        \begin{align*}
          \ X E(Y) & =  \phi(Sy^{-1}_2 y^3) \ g^2(h) \ (x^{-1} \rtimes f^0 \rtimes x^1 \rtimes f^2) \ ( y^{-1}_1 \rtimes g^4 \rtimes y^5 \rtimes g^6)\\
          & = \phi(Sy^{-1}_3 y^3) \ g^2(h) \ f^0_1(y^{-1}_2) \ g^4_2(x^1_1) \ f^2_1(y^5_2) \ x^{-1} y^{-1}_1 \rtimes f^0_2 g^4_1 \rtimes x^1_2 y^5_1 \rtimes 
          f^2_2 g^6 \\
          & = S \phi_1 (y^{-1}_3) \ \phi_2(y^3) \ f^0_1(y^{-1}_2) \ g^2(h) \ g^4_2(x^1_1) \ f^2_1(y^5_2) \ x^{-1} y^{-1}_1 \rtimes f^0_2 g^4_1 \rtimes x^1_2 y^5_1 \rtimes 
          f^2_2 g^6 \\
           & = (f^0_1 S \phi_1)(y^{-1}_2) \ \phi_2(y^3) \ g^2(h) \ g^4_2(x^1_1) \ f^2_1(y^5_2) \ x^{-1} y^{-1}_1 \rtimes f^0_2 g^4_1 \rtimes 
           x^1_2 y^5_1 \rtimes f^2_2 g^6  
        \end{align*}
        which, by an appeal to the formula $fS\phi_1 \otimes \phi_2 = S\phi_1 \otimes \phi_2 f$, is seen to be equal to
        \begin{align*}
         & (S \phi_1)(y^{-1}_2) \ (\phi_2 f^0_1) \ (y^3) \ g^2(h) \ g^4_2(x^1_1) \ f^2_1(y^5_2) \ x^{-1} y^{-1}_1 \rtimes f^0_2 g^4_1 \rtimes x^1_2 y^5_1 \rtimes 
          f^2_2 g^6 \\
         & =  \phi(Sy^{-1}_2 y^{3}_1) \ g^2(h) \ f^0_1(y^3_2) \ g^4_2(x^1_1) f^2_1(y^5_2) \ x^{-1} y^{-1}_1 \rtimes f^0_2 g^4_1 \rtimes 
        x^1_2 y^5_1 \rtimes f^2_2 g^6
        \end{align*}
        and this clearly equals $E(XY)$. Similarly, we can show that $E(YX) = E(Y)X$, completing the proof. 
        \end{proof}      
        
        Next, we apply Lemma \ref{basic} and Lemma \ref{exp} to explicitly describe certain basic constructions and their associated Jones projections.

\begin{proposition}\label{basic cons}
 The following are instances of basic constructions with the Jones projections being specified pictorially in appropriate planar algebras.
 \begin{itemize}
 \item[1.]
(i) For integers $k \leq p \leq q \leq l$, 
$H_{[p, q]} \subset H_{[k, l]} \subset H_{[2k - p, 2l - q]}$ is an instance of the basic construction with the Jones projection 
in $P(H^{2k-p})_{2l-2k-q+p+2}$ given by
 \begin{center}
  \begin{tikzpicture}
   \node at (0.2,0) {\small $ \delta^{-(p+l-k-q)}$};
   \draw [black,thick] (1.5,.5) to [out=275, in=265] (2,.5); 
   \draw [black,thick] (1.5,-.5) to [out=85, in=95] (2,-.5);
   \node [right] at (1.5,.2) {\tiny $p-k$};
   \node [right] at (1.5, -.18) {\tiny $p-k$};
   \draw [black,thick] (2.75,.5) -- (2.75,-.5);
   \node [right] at (2.75,0) {\tiny $q-p+2$};
   \draw [black,thick] (4.3,.5) to [out=275, in=265] (4.8,.5); 
  \draw [black,thick] (4.3,-.5) to [out=85, in=95] (4.8,-.5);
  \node [right] at (4.3,.18) {\tiny $l-q$}; 
  \node [right] at (4.3, -.18) {\tiny $l-q$};

  \end{tikzpicture}

 \end{center}
and $tr_{H_{[k, l]}}$ is a Markov trace of modulus $\delta^{2(p+l-k-q)}$
for the inclusion $ H_{[p, q]} \subset H_{[k, l]}$.

(ii) For any non-negative integer $k$, $\mathbb{C} \subset H_{[0, k]} \subset H_{[0, 2k+1]}$ is an instance of the basic construction with the Jones projection in 
$P(H^*)_{2k+3}$ given by 
 \begin{center}
      \begin{tikzpicture}
 \node at (.4,0) {\tiny $\delta^{-(k+1)}$};
 \draw [black,thick] (1.3,.5) -- (1.3,-.5);
 \node [right] at (.95,0) {\tiny $1$}; 
 \draw [black,thick] (1.7,.5) to [out=275, in=265] (2.2,.5); 
  \draw [black,thick] (1.7,-.5) to [out=85, in=95] (2.2,-.5);
  \node [below] at (2.2,.42) {\tiny $k+1$};
  \node [above] at (2.2,-.42) {\tiny $k+1$};
  \end{tikzpicture}
  \end{center}
  and $tr_{H_{[0, k]}}$ is a Markov trace of modulus $\delta^{2(k+1)}$ for the inclusion $\mathbb{C} \subset H_{[0, k]}$.

(iii) For any non-negative integer $k$, 
$\mathbb{C} \subset H_{[-k, 0]} \subset H_{[-2k-1, 0]}$  is an instance of the basic construction
with the Jones projection in $P(H)_{2k+3}$ given by 
 \begin{center}
  \begin{tikzpicture}
   \node at (.8,0) {\tiny $\delta^{-(k+1)}$};
 \draw [black,thick] (1.7,.5) to [out=275, in=265] (2.2,.5); 
  \draw [black,thick] (1.7,-.5) to [out=85, in=95] (2.2,-.5);
  \node [below] at (2.2,.45) {\tiny $k+1$};
  \node [above] at (2.2,-.45) {\tiny $k+1$};
  \draw [black,thick] (3,.5) -- (3, -.5);
  \node [right] at (3.001,0) {\tiny $1$}; 
\end{tikzpicture}
\end{center} 
and $tr_{H_{[-k, 0]}}$ is a Markov trace of modulus $\delta^{2(k+1)}$ for the inclusion $ \mathbb{C} \subset H_{[-k, 0]}$. 

\item[2.] If $l \geq 1, s \geq 0$  are integers, 
$H_{[-l, -1]} \otimes H_{[2, 2+s]} \subset H_{[-l, 2+s]} \subset H_{[-l, -1]} \otimes H_{[2, 6+s]}
(\subset H_{[-l, 6+s]} \cong P(H^l)_{s+l+8})$ is an instance of the basic construction with the Jones projection given by the following figure
\begin{center}
\begin{tikzpicture}
\node at (1.3,0) {\tiny $\delta^{-2}$};
 \draw [black,thick] (1.96,.5) -- (1.96,-.5);
 \node [right] at (2,0) {\tiny $l+2$}; 
 \draw [black,thick] (3,.5) to [out=275, in=265] (3.5,.5); 
  \draw [black,thick] (3,-.5) to [out=85, in=95] (3.5,-.5);
  \node [below] at (3.3,.4) {\tiny $2$};
  \node [above] at (3.3,-.4) {\tiny $2$};
  \draw [black,thick] (4,.5) -- (4, -.5);
  \node [right] at (4.002,0) {\tiny $s+2$};
 \end{tikzpicture}
 \end{center}
where the first inclusion is natural and the second inclusion is given by the map $\psi_{l, s}$ as defined in the statement of Lemma 
\ref{exp}(ii).
Furthermore, $tr_{H_{[-l, 2+s]}}$ is a Markov trace of modulus $\delta^{4}$ for the inclusion 
$H_{[-l, -1]} \otimes H_{[2, 2+s]} \subset H_{[-l, 2+s]}$.

 \item[3.]  If $l \geq 1, s \geq 0$ are integers, 
      $H_{[-l, 2+s]} \subset H_{[-l, -1]} \otimes H_{[2, 6+s]} \subset H_{[-l, 6+s]} (\cong P(H^l)_{s+l+8})$ is an instance of the 
      basic construction with the  Jones projection given by
      \begin{center}
      \begin{tikzpicture}
 \node at (.7,0) {\tiny $\delta^{-2}$};
 \draw [black,thick] (1.4,.5) -- (1.4,-.5);
 \node [right] at (1.42,0) {\tiny $l$}; 
 \draw [black,thick] (2,.5) to [out=275, in=265] (2.5,.5); 
  \draw [black,thick] (2,-.5) to [out=85, in=95] (2.5,-.5);
  \node [below] at (2.3,.4) {\tiny $2$};
  \node [above] at (2.3,-.4) {\tiny $2$};
  \draw [black,thick] (3,.5) -- (3, -.5);
  \node [right] at (3.001,0) {\tiny $s+4$};
 \end{tikzpicture}
      \end{center}
      where the first inclusion is given by the map $\psi_{l, s}$ as defined in the statement of Lemma \ref{exp}(ii) and the second inclusion is the natural 
      inclusion. Also, $tr_{H_{[-l, -1]} \otimes H_{[2, 6+s]}}$ is a Markov trace of modulus $\delta^{4}$ for the inclusion
      $H_{[-l, s+2]} \subset H_{[-l, -1]} \otimes H_{[2, 6+s]}$.
      \end{itemize}
      \end{proposition}

      \begin{proof}
      The strategy for the proof of Proposition \ref{basic cons} is to verify, in each case, conditions (i), (ii), and (iii) of
      Lemma \ref{basic}. We will frequently use Lemma \ref{iso1} without any mention in the proofs of all parts of Proposition \ref{basic cons}.
      \begin{itemize}
       \item [1.] 
       (i) For notational convenience, we assume that $k, p, q, l$ are all odd so that we may identify $H_{[2k-p, 2l-q]}$ with
       $P(H)_{2l-2k+p-q+2}$ and we regard $H_{[k,l]}, H_{[p, q]}$ as subalgebras of  $P(H)_{2l-2k+p-q+2}$. 
        Let $e_1$ denote the projection defined in the statement of Proposition \ref{basic cons}(1)(i).
        Given $Z \in H_{[k, l]}$, the natural inclusion of $H_{[k, l]}$ inside $H_{[2k-p, 2l-q]}$ shows that $Z$ is identified with the element of 
        $P(H)_{2l-2k+p-q+2}$ as shown on the left in Figure \ref{fig:ZW}. Similarly, any element $W \in H_{[p, q]}$ is identified with the element of 
        $P(H)_{2l-2k+p-q+2}$ as shown on the right in Figure \ref{fig:ZW}. 
         \begin{figure}[!h]
\begin{center}
\psfrag{a}{\Huge $Z_{T^{l-k+2}}^{P(H)}(Z)$}
\psfrag{b}{\Huge $Z_{T^{q-p+2}}^{P(H)}(W)$}
\psfrag{c}{\Huge $p-k$}
\psfrag{d}{\Huge $l-q$}
\psfrag{e}{\Huge $2p-2k$}
\psfrag{f}{\Huge $2l-2q$}
\resizebox{8.5cm}{!}{\includegraphics{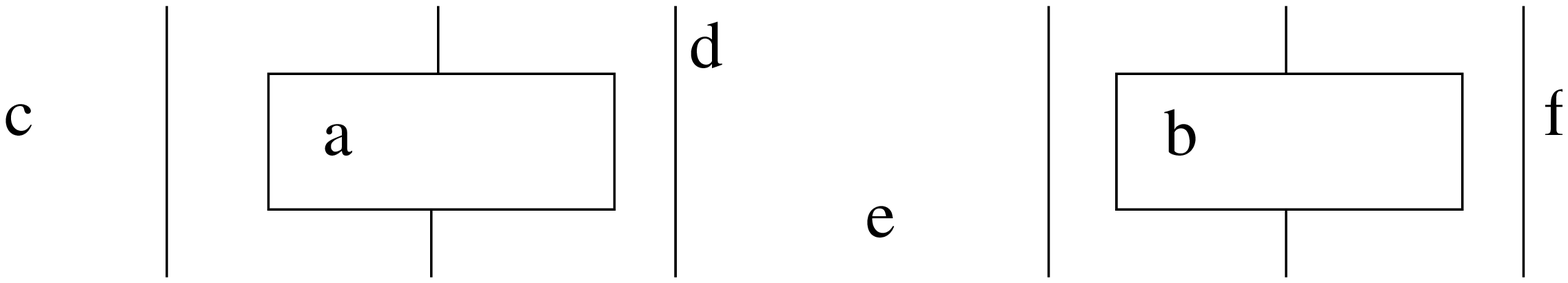}}
\end{center}
\caption{Elements: $Z$ (left) and $W$ (right)}
\label{fig:ZW}
\end{figure}
       Conditions (i) and (ii) of Lemma \ref{basic} are straightforward to verify. For condition (iii), one can easily observe that
       if we take 
       \begin{eqnarray*}
        X &=& x^k \rtimes f^{k+1} \rtimes \cdots \rtimes x^l,\\ Y &=& y^k \rtimes g^{k+1} \rtimes \cdots \rtimes g^{p-1} \underbrace{\rtimes 1 \rtimes 
        \epsilon \rtimes \cdots \rtimes 1}_\text{$q-p+1$ \mbox{terms}} \rtimes \ g^{q+1} \rtimes \cdots \rtimes y^l
       \end{eqnarray*}
       in $H_{[k, l]}$, then $Xe_1Y$ equals the element
        \begin{align*}
         \delta ^{-(p+l-k-q)} Z_U^{P(H)} & ( x^k \otimes Ff^{k+1} \otimes \cdots \otimes Ff^{l-1} \otimes x^l \otimes y^k \otimes Fg^{k+1}
         \otimes \cdots 
         \otimes Fg^{p-1}\\ & \otimes Fg^{q+1} \otimes y^{q+2} \otimes \cdots \otimes y^l)
       \end{align*}
       where  $U \in {\mathcal T}{(2l - 2k +p -q +2)}$ is the tangle as shown in Figure \ref{fig:pic26}.
        \begin{figure}[!h]
\begin{center}
\psfrag{a}{\Huge $1$}
\psfrag{b}{\Huge $2$}
\psfrag{c}{\huge $p-k-1$}
\psfrag{d}{\huge $p-k$}
\psfrag{e}{\huge $p-k+1$}
\psfrag{f}{\huge $p-k+2$}
\psfrag{g}{\huge $q-k+1$}
\psfrag{h}{\huge $q-k+2$}
\psfrag{i}{\huge $q-k+3$}
\psfrag{j}{\huge $l-k+1$}
\psfrag{k}{\huge $l-k+2$}
\psfrag{l}{\huge $l-k+3$}
\psfrag{m}{\huge $l+p-2k$}
\psfrag{n}{\huge $l+p-2k+1$}
\psfrag{o}{\huge $l+p-2k+2$}
\psfrag{p}{\huge $l+p-2k+3$}
\psfrag{q}{\huge $2l-2k+p-q+1$}
\resizebox{10.0cm}{!}{\includegraphics{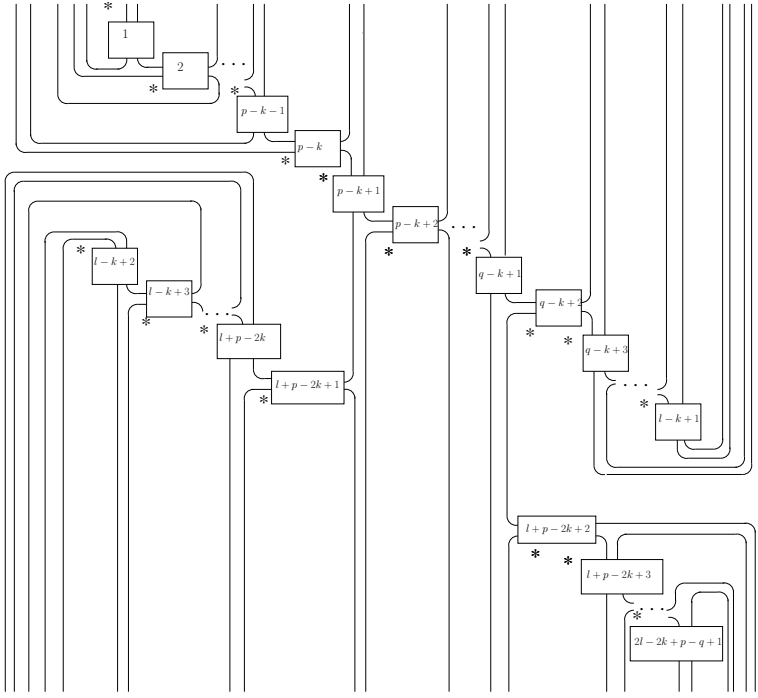}}
\end{center}
\caption{Tangle $U$}
\label{fig:pic26}
\end{figure}
       Thus we see that 
      $H_{[k, l]} e_1 H_{[k, l]}$ contains the image of the linear isomorphism $Z_U$ 
       and thus, $dim ~(H_{[k, l]} e_1 H_{[k, l]}) \geq$ rank of $Z_U = (dim ~H)^{2l - 2k +p -q +1}$ and whence the equality  
       $H_{[k, l]} e_1 H_{[k, l]} = H_{[2k - p, 2l - q]}$ follows. A routine computation in $P(H)_{2l-2k+p-q+2}$ 
       verifies that $tr_{H_{[k, l]}}$ is a 
       Markov trace of modulus $\delta^{2(p+l-k-q)}$ for the inclusion $H_{[p, q]} \subset H_{[k, l]}$ and we leave the verification
       to the reader. The proof for the case when $k, p, q, l$ are not all odd follows in a similar way and the reader should note that 
       if $2k-p$ is even, then $H_{[2k-p, 2l-q]}$ is identified with $P(H^*)_{2l-2k+p-q+2}$ and consequently, $H_{[k, l]}, H_{[p, q]}$ are 
       regarded as subalgebras of $P(H^*)_{2l-2k+p-q+2}$. 
       
       (ii) Identify $H_{[0, 2k+1]}$ with $P(H^*)_{2k+3}$. Then $H_{[0, k]}$ is identified with the subalgebra $P(H^*)_{k+2}$ of 
       $P(H^*)_{2k+3}$ while $\mathbb{C}$ is identified with the space of scalar multiples of the unit element of $P(H^*)_{2k+3}$.
       Let us use the symbol $e_2$ to denote the projection defined in the statement of Proposition \ref{basic cons}(1)(ii).
       Note that conditions (i) and (ii) of Lemma \ref{basic} are easy to verify. 
       A simple computation in $P(H^*)_{2k+3}$ shows that $H_{[0, k]} e_2 H_{[0, k]}$ 
       equals the image of the linear isomorphism $Z_T$ induced by the tangle $T \in {\mathcal T}{(2k+3)}$ as depicted in Figure 
       \ref{fig:pic27} and by comparing dimensions we conclude that $H_{[0, k]} e_2 H_{[0, k]} =$ Range of $Z_T = P(H^*)_{2k+3}$, verifying
       the condition (iii) of Lemma \ref{basic}.
       \begin{figure}[!h]
\begin{center}
\psfrag{a}{\huge $T^{k+2}$}
\psfrag{b}{\huge $k+1$}
\psfrag{c}{\huge $k+2$}
\psfrag{1}{\huge $1$}
\resizebox{2cm}{!}{\includegraphics{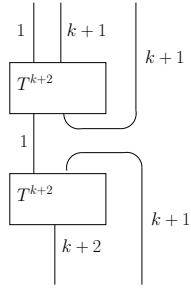}}
\end{center}
\caption{Tangle $T$}
\label{fig:pic27}
\end{figure} 
       
      (iii) Proceed along the same line of argument as in the proof of part (ii).
       
       \item[2.] Again, for notational convenience, rather than proving the result in its full generality, we shall just present the proof when
       $l=1, s = 0$. 
       The proof for the general case will follow in a similar fashion. Given 
       \begin{align*}
      X = x^{-1} \rtimes f^0 \rtimes x^1 \rtimes f^2 \in H_{[-1, 2]}, 
       \end{align*}
note that the image of $X$ under $\psi_{1, 0}$, also denoted $X$, can be identified with the element of $P(H)_9$ as shown on the left in  
       Figure \ref{fig:pic8}        
      \begin{figure}[!h]
\begin{center}
\psfrag{a}{\Huge $x^{-1}_1$}
\psfrag{b}{\Huge $x^{-1}_2$}
\psfrag{c}{\Huge $Ff^0$}
\psfrag{d}{\Huge $x^1$}
\psfrag{e}{\Huge $Ff^2$}
\psfrag{f}{\Huge $z^{-1}_1$}
\psfrag{g}{\Huge $z^{-1}_2$}
\psfrag{h}{\Huge $Fg^2$}
\resizebox{10cm}{!}{\includegraphics{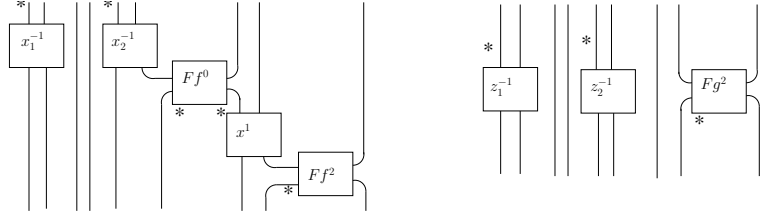}}
\end{center}
\caption{$X$ (left) and $Z \in H_{-1} \otimes H_2$ (right)}
\label{fig:pic8}
\end{figure}             
and consequently, $Z = z^{-1} \otimes g^2 \in H_{-1} \otimes H_2$ is identified with the element of $P(H)_9$ 
as shown on the right in Figure \ref{fig:pic8}.
Let $e$ denote the projection defined in the statement of Proposition \ref{basic cons}(2).

Note that the element $Ze$ is as shown on the left in Figure \ref{fig:D1} which, by an appeal to the Relation C, 
is easily seen to be equal to the element shown on the right in Figure \ref{fig:D1}. From this it follows immediately that $Z = 0$ whenever $Ze  = 0$,
verifying condition (i) of Lemma \ref{basic}.
\begin{figure}[!h]
\begin{center}
\psfrag{a}{\Huge $z^{-1}_1$}
\psfrag{b}{\Huge $z^{-1}_2$}
\psfrag{d}{\Huge $z^{-1}$}
\psfrag{c}{\Huge $Fg^2$}
\psfrag{e}{\Huge $=$}
\psfrag{f}{\Huge $\delta^{-2}$}
\resizebox{9cm}{!}{\includegraphics{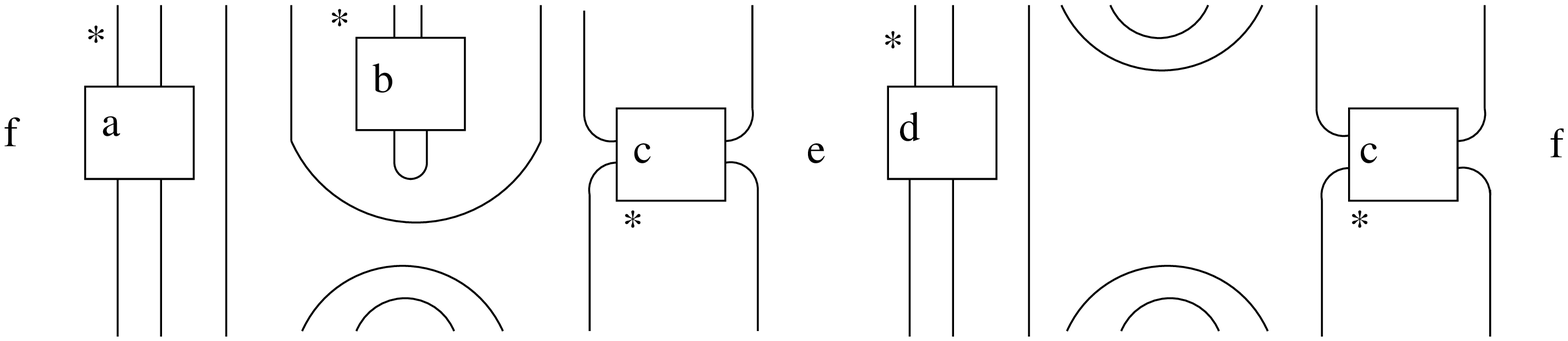}}
\end{center}
\caption{$Ze$}
\label{fig:D1}
\end{figure}

It follows from Lemma \ref{exp}(i) that $E(X) = f^0(h) \phi(x^1) \ x^{-1} \otimes f^2$ where $E$ is the trace-preserving conditional expectation 
of $H_{[-1, 2]}$ onto $H_{-1} \otimes H_2$. Observe that the element $eXe$ equals the element shown on the left in Figure \ref{fig:D2}. Now applying
the Relations C and T, we reduce the element on the left in Figure \ref{fig:D2} to that on the right in Figure \ref{fig:D2}. 
With $Z = E(X)$, the pictorial description of the element $Ze$ as shown on the right in Figure \ref{fig:D1} immediately yields that the 
element on the right in Figure \ref{fig:D2} is indeed $E(X) e$, verifying condition (ii) of Lemma \ref{basic}. 
\begin{figure}[!h]
\begin{center}
\psfrag{a}{\Huge $\delta^{-4}$}
\psfrag{b}{\Huge $x^{-1}_1$}
\psfrag{c}{\Huge $x^{-1}_2$}
\psfrag{d}{\Huge $Ff^0$}
\psfrag{e}{\Huge $=$}
\psfrag{f}{\Huge $x^1$}
\psfrag{g}{\Huge $Ff^2$}
\psfrag{h}{\Huge $x^{-1}$}
\psfrag{i}{\Huge $Ff^2$}
\psfrag{j}{\Huge $\delta^{-2} \ f^0(h) \ \phi(x^1)$}
\resizebox{10.5cm}{!}{\includegraphics{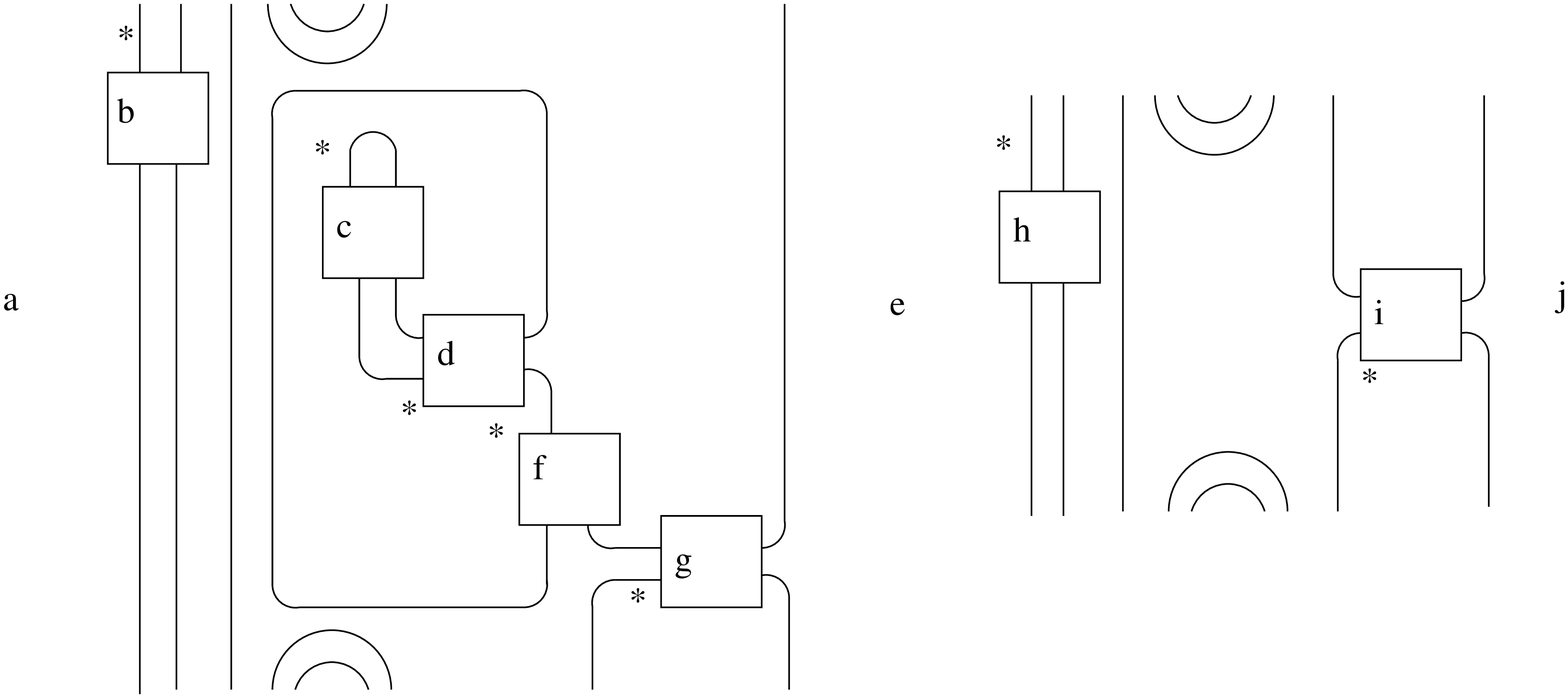}}
\end{center}
\caption{$eXe$}
\label{fig:D2}
\end{figure}

To verify condition (iii), we note that if we take elements 
  \begin{align*}         
       X = 1 \rtimes f^0 \rtimes x^1 \rtimes f^2, \ Y = y^{-1} \rtimes g^0 \rtimes y^1 \rtimes \epsilon
       \end{align*}
        in $H_{[-1, 2]}$, then $XeY$ equals the element
       \begin{align*}
      &  \delta^{-2} Z_S^{P(H)} (y^{-1} \otimes Ff^0 \otimes x^1 \otimes Ff^2 \otimes Fg^0 \otimes y^1)
       \end{align*}
       where $Z_S : H^{\otimes 6} \rightarrow
       P(H)_9$ is the injective linear map induced by the tangle $S \in \mathcal{T}(9, 6)$ as shown on the left in Figure \ref{fig:D3}.  
  \begin{figure}[!h]
\begin{center}
\psfrag{a}{\Huge $\delta^{-2}$}
\psfrag{b}{\Huge $x^{-1}_1$}
\psfrag{c}{\Huge $x^{-1}_2$}
\psfrag{d}{\Huge $Ff^0$}
\psfrag{e}{\Huge $x^1$}
\psfrag{f}{\Huge $Ff^2$}
\psfrag{1}{\Huge $1$}
\psfrag{2}{\Huge $2$}
\psfrag{3}{\Huge $3$}
\psfrag{4}{\Huge $4$}
\psfrag{5}{\Huge $5$}
\psfrag{6}{\Huge $6$}
\resizebox{10cm}{!}{\includegraphics{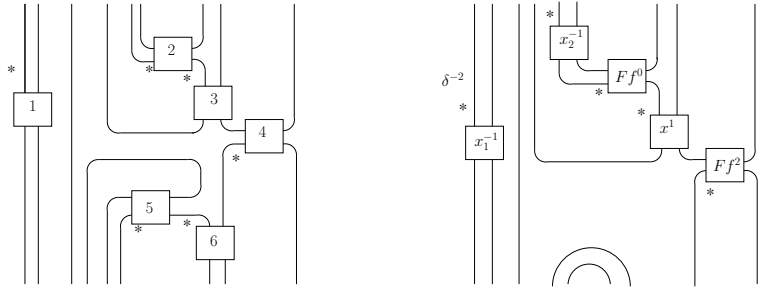}}
\end{center}
\caption{Tangle S (left) and $Xe$ (right)}
\label{fig:D3}
\end{figure}        
       It thus follows that $H_{[-1, 2]} e H_{[-1, 2]}$ contains the image of the map $Z_S$.            
       Now, by comparison of dimensions we see that 
       \begin{eqnarray*}
       dim ~(H_{[-1, 2]} e H_{[-1, 2]}) \geq \text{ rank of } Z_S = (dim ~H)^6 =  dim ~(H_{[-1, -1]} \otimes H_{[2, 6]})
       \end{eqnarray*}
       and obviously  $$dim ~(H_{[-1, -1]} \otimes H_{[2, 6]}) \geq 
       dim ~(H_{[-1, 2]} e H_{[-1, 2]})$$ as $H_{[-1, 2]} e H_{[-1, 2]}$ is contained in $H_{[-1, -1]} \otimes H_{[2, 6]}.$
       Thus, 
       condition (iii) of Lemma \ref{basic} is verified.
       
       Note that if $X \in H_{[-1, 2]}$ is as before, then $Xe$ equals the element shown on the right in Figure \ref{fig:D3}.
       It is then a routine computation to verify that $tr(Xe) = \delta^{-4} tr(X)$, proving that
$tr_{H_{[-1, 2]}}$ is a Markov trace of modulus $\delta^{4}$ for the inclusion 
$H_{-1} \otimes H_2 \subset H_{[-1, 2]}$.

       \item[3.] As before, we only present the proof when
       $l = 1$ and $s=0$, omitting the proof for the general case which is analogous. Let $e$ denote the projection defined in the statement
       of Proposition \ref{basic cons}(3).
       We identify as usual $H_{[-1, 6]}$ with $P(H)_9$. 
       
       Given 
       \begin{align*}
        X = x^{-1} \rtimes f^0 \rtimes x^1 \rtimes f^2 \in H_{[-1, 2]},
       \end{align*}
 its image in $H_{[-1, 6]}$ is given by 
 \begin{align*}
  x^{-1} \rtimes \epsilon \rtimes 1 \rtimes \epsilon \rtimes x^{-1}_2 \rtimes f^0 \rtimes x^1 \rtimes f^2.
 \end{align*}        
        The element $eX$ is shown on the left in Figure \ref{fig:D4}. An application of the Relation E shows that $eX$ equals the element on the right 
        in Figure \ref{fig:D4}. Similarly, by an appeal to the Relations A and E, one can easily see that the element $Xe$ equals the element on the 
        right in Figure \ref{fig:D4} so that $e X = X e$. Thus, we conclude that $e$ commutes with $X$. Further, it is evident from the pictorial representation
        of the element $Xe$ as shown on the right in Figure \ref{fig:D4} that the map $X \mapsto Xe$ of $H_{[-1, 2]}$ into $H_{[-1, 6]}$ is injective,
        verifying condition (i) of Lemma \ref{basic}.
        \begin{figure}[!h]
\begin{center}
\psfrag{a}{\Huge $\delta^{-2}$}
\psfrag{b}{\Huge $x^{-1}_1$}
\psfrag{c}{\Huge $x^{-1}_2$}
\psfrag{d}{\Huge $Ff^0$}
\psfrag{e}{\Huge $=$}
\psfrag{f}{\Huge $x^1$}
\psfrag{g}{\Huge $Ff^2$}
\psfrag{h}{\Huge $x^{-1}$}
\resizebox{10cm}{!}{\includegraphics{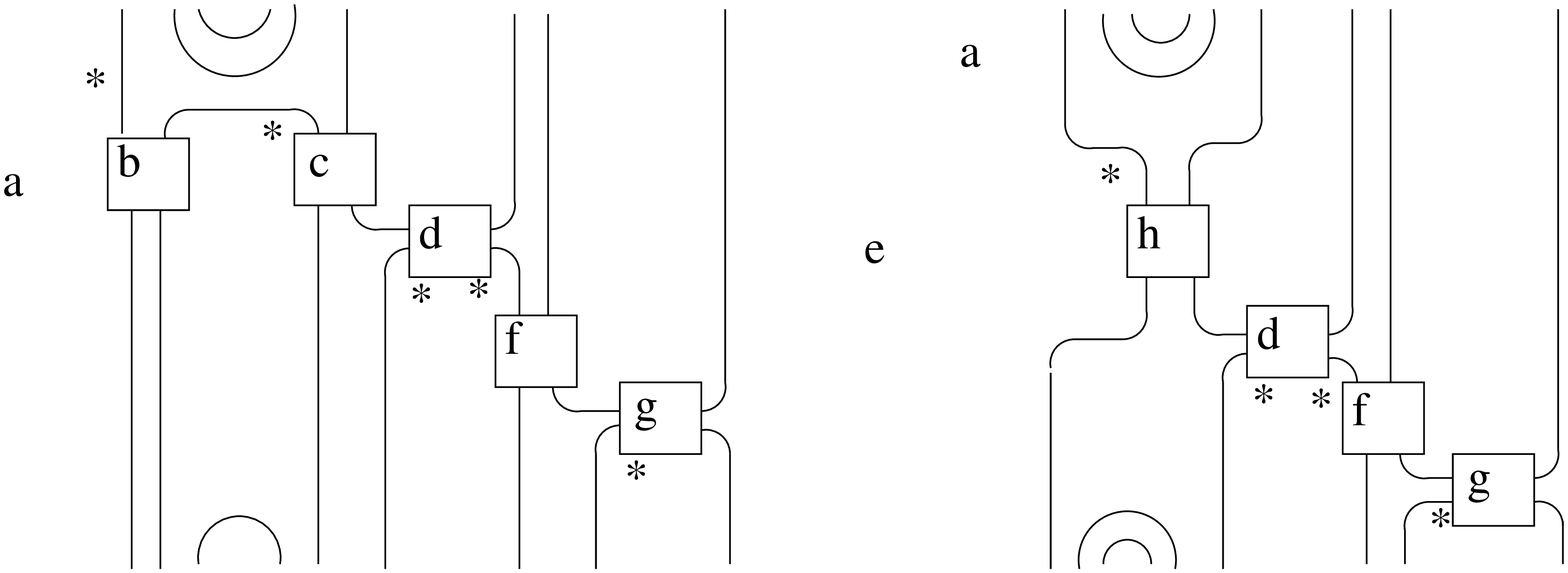}}
\end{center}
\caption{$eX = Xe$}
\label{fig:D4}
\end{figure}

         Given
       \begin{align*}            
        X = x^{-1} \otimes (f^2 \rtimes x^3 \rtimes f^4 \rtimes x^5 \rtimes f^6) \in H_{-1} \otimes H_{[2, 6]},
       \end{align*}
       we observe that $eXe$ equals the element in $P(H)_9$ shown on the left in Figure \ref{fig:D5}.
       \begin{figure}[!h]
\begin{center}
\psfrag{a}{\Huge $\delta^{-4}$}
\psfrag{b}{\Huge $x^{-1}$}
\psfrag{c}{\Huge $Ff^2$}
\psfrag{d}{\Huge $x^3$}
\psfrag{f}{\Huge $Ff^4$}
\psfrag{g}{\Huge $x^5$}
\psfrag{h}{\Huge $Ff^6$}
\psfrag{e}{\Huge $=$}
\psfrag{i}{\Huge $\delta^{-2} \ f^2(h) \ \phi(Sx^{-1}_3 x^3)$}
\psfrag{j}{\Huge $x^{-1}_1$} 
\psfrag{k}{\Huge $x^{-1}_2$}
\resizebox{12cm}{!}{\includegraphics{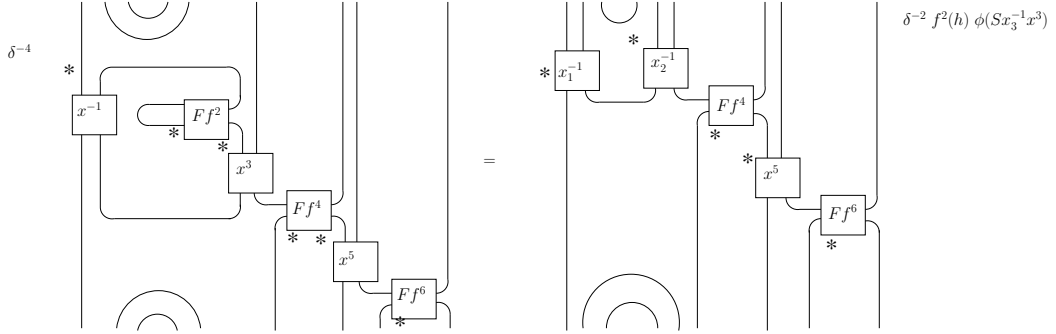}}
\end{center}
\caption{$eXe = E(X)e$}
\label{fig:D5}
\end{figure}
Firstly an application of the Relation C, then repeated applications of the Relations A and E and finally, an application of the Relation T
reduces the element on the left in Figure
\ref{fig:D5} to that on the right in Figure \ref{fig:D5}. If $E$ denotes the trace-preserving conditional expectation of 
$H_{[-1, -1]} \otimes H_{[2, 6]}$ onto $H_{[-1, 2]}$, then by Lemma \ref{exp}(ii), $E(X) = \phi(Sx^{-1}_2 x^3) \ f^2(h) \ x^{-1}_1 \rtimes f^4 
\rtimes x^5 \rtimes f^6$. Representing $E(X)e$ pictorially in $P(H)_9$, one can easily see that $E(X)e$ indeed equals the element on the right 
in Figure \ref{fig:D5}. Therefore, $eXe = E(X) e$, verifying condition (ii) of Lemma \ref{basic}.

Finally it just remains to verify that $(H_{-1} \otimes H_{[2, 6]}) \ e \ (H_{-1} \otimes H_{[2, 6]}) = H_{[-1, 6]}$.
Consider elements $X, Y$ in $H_{-1} \otimes H_{[2, 6]}$ given by 
        \begin{align*}
          X = x^{-1} \otimes (f^2 \rtimes x^3 \rtimes f^4 \rtimes x^5 \rtimes f^6), \ \mbox{and} \ 
         Y = 1 \otimes (g^2 \rtimes y^3 \rtimes \epsilon \rtimes 1 \rtimes \epsilon).         
        \end{align*}
        Then one can easily see that $XeY$ equals the element
        \begin{align*}
  Z_T^{P(H)}(x^{-1} \otimes Ff^2 \otimes x^3 \otimes Ff^4 \otimes x^5 \otimes Ff^6 \otimes Fg^2 \otimes y^3)
        \end{align*}        
        where $Z_T$ is the linear 
        isomorphism induced 
        by the tangle $T \in {\mathcal{T}}(9)$ as shown in Figure \ref{fig:D6}. Thus
        we see that $(H_{-1} \otimes H_{[2, 6]}) \ e \ (H_{-1} \otimes H_{[2, 6]})$ contains the image of
        $Z_T$. Then by comparing dimensions of spaces we have that
        $ H_{[-1, 6]} = (H_{-1} \otimes H_{[2, 6]}) \ e \ (H_{-1} \otimes H_{[2, 6]})$.         
 \begin{figure}[!h]
\begin{center}
\psfrag{1}{\Huge $1$}
\psfrag{2}{\Huge $2$}
\psfrag{3}{\Huge $3$}
\psfrag{4}{\Huge $4$}
\psfrag{5}{\Huge $5$}
\psfrag{6}{\Huge $6$}
\psfrag{7}{\Huge $7$}
\psfrag{8}{\Huge $8$}
\resizebox{4.0cm}{!}{\includegraphics{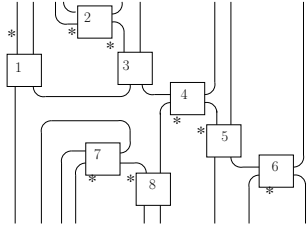}}
\end{center}
\caption{Tangle T}
\label{fig:D6}
\end{figure}

Finally, a routine computation shows that for any 
$X \in H_{-1} \otimes H_{[2, 6]}, \ tr(Xe) = \delta^{-4} tr(X)$,
 so that $tr_{H_{-1} \otimes H_{[2, 6]}}$ is a Markov trace of modulus $\delta^{4}$ for the inclusion $H_{[-1, 2]} \subset 
H_{-1} \otimes H_{[2, 6]}$, completing the proof.          
      \end{itemize}
    
      \end{proof}
       \subsection{Jones' basic construction tower of $\mathcal{N} \subset \mathcal{M}$ and relative commutants}
      The goal of this subsection is to explicitly determine the basic construction tower of $\mathcal{N} \subset \mathcal{M}$. 
      
      We begin with proving that certain squares of finite-dimensional $C^*$-algebras are symmetric commuting squares. 
\begin{lemma}\label{comm sq}
If $k < p < q < l$ are positive integers, then the square in Figure \ref{fig:pic35} is an instance of a symmetric commuting square with respect to $tr_{H_{[k, l]}}$ which is a Markov trace of modulus $\delta^{2(p-k+l-q+2)}$ for the inclusion $H_{[p+1, q-1]} \subset H_{[k, l]}$.
 \begin{figure}[!h]
\begin{center}
\psfrag{a}{\Large $H_{[p+1, q-1]}$}
\psfrag{b}{\Large $\subset$}
\psfrag{c}{\Large $H_{[k, l]}$}
\psfrag{d}{\Large $\cup$}
\psfrag{e}{\Large $\mathbb{C}$}
\psfrag{f}{\Large $\subset$}
\psfrag{g}{\Large $H_{[k, p]} \otimes H_{[q, l]}$}
\psfrag{h}{\Large $\cup$}
\resizebox{4cm}{!}{\includegraphics{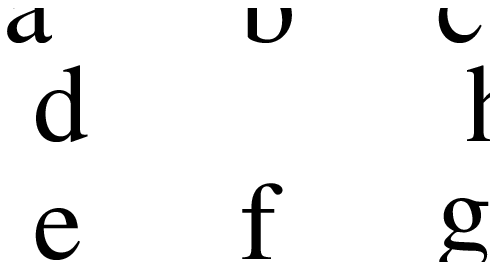}}
\end{center}
\caption{Commuting square}
\label{fig:pic35}
\end{figure}
   \end{lemma}
\begin{proof}
    By Lemma \ref{exp}(i) the square of finite-dimensional $C^*$-algebras in Figure \ref{fig:pic35} is a commuting square with respect to $tr_{H_{[k, l]}}$.
    In order to show that this
    is symmetric, we need to verify that $H_{[k, l]}$ is linearly spanned by $(H_{[k, p]} \otimes H_{[q, l]}) H_{[p+1, q-1]}$.
    Assume that $k, p, q, l$ are all odd so that we may identify $H_{[k, l]}$ with $P(H)_{l-k+2}$ and also identify 
    $H_{[k, p]} \otimes H_{[q, l]}, H_{[p+1, q-1]}$ as subalgebras of $P(H)_{l-k+2}$. A little computation in $P(H)$ shows that
    $(H_{[k, p]} \otimes H_{[q, l]}) H_{[p+1, q-1]}$ equals the image of the  linear isomorphism $Z_W$ where $W \in 
    {\mathcal T}(l-k+2)$ is as shown in Figure \ref{fig:pic34} which in turn implies, by comparing dimensions of spaces, that 
    $(H_{[k, p]} \otimes H_{[q, l]}) H_{[p+1, q-1]} = H_{[k, l]}$ and the desired result follows. Further, it is a direct consequence of the
    Proposition \ref{basic cons}(1)(i) that $tr_{H_{[k, l]}}$ is a Markov trace of modulus $\delta^{2(p-k+l-q+2)}$ 
    for the inclusion $H_{[p+1, q-1]} \subset H_{[k, l]}$. The general proof follows in a similar fashion with the difference that when $k$ is even,
    we identify $H_{[k, l]}$ with $P(H^*)_{l-k+2}$. 
    \begin{figure}[!h]
\begin{center}
\psfrag{a}{\Huge $1$}
\psfrag{b}{\Huge $2$}
\psfrag{c}{\Huge $q-k-1$}
\psfrag{d}{\Huge $q-k$}
\psfrag{e}{\Huge $q-k+1$}
\psfrag{f}{\Huge $q-k+2$}
\psfrag{g}{\Huge $l-k+1$}
\resizebox{6.0cm}{!}{\includegraphics{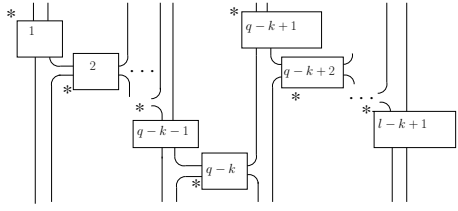}}
\end{center}
\caption{Tangle $W$}
\label{fig:pic34}
\end{figure}    
 \end{proof}

\begin{figure}[!h]
\begin{center}
\psfrag{a}{\Large $A_{0, 1}$}
\psfrag{b}{\Large $\subset$}
\psfrag{c}{\Large $A_{1, 1}$}
\psfrag{d}{\Large $\cup$}
\psfrag{e}{\Large $\cup$} 
\psfrag{f}{\Large $A_{0, 0}$}
\psfrag{g}{\Large $\subset$}
\psfrag{h}{\Large $A_{1, 0}$}
\resizebox{3.0cm}{!}{\includegraphics{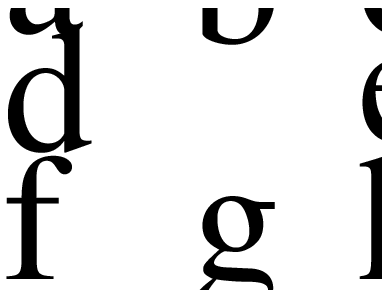}}
\end{center}
\caption{Commuting square}
\label{fig:pic37}
\end{figure}

      We set $A_{0, 0} = \mathbb{C}, \ A_{1, 0} = H_{-1} \otimes H_2, \ A_{0, 1} = H_{[0, 1]}$ and $ A_{1, 1} = H_{[-1, 2]}$.
      It is an immediate consequence of Lemma \ref{comm sq} that the square in Figure \ref{fig:pic37} is a symmetric commuting square with respect 
      to $tr_{A_{1, 1}}$ which is a Markov trace of modulus $\delta^4$ for the inclusion $A_{0, 1} \subset A_{1, 1}$. Further,
 all the inclusions are connected since the lower left corner is ${\mathbb C}$ while the upper right corner is a matrix algebra by 
 Lemma \ref{matrixalg}.
 
\begin{figure}[!h]
\begin{center}
\psfrag{a}{\Huge $\delta^{-2}$}
\psfrag{b}{\Huge $1$}
\psfrag{c}{\Huge $1$}
\psfrag{d}{\Huge $2k+3$}
\psfrag{e}{\Huge $1$} 
\psfrag{f}{\Huge $1$}
\resizebox{3.5cm}{!}{\includegraphics{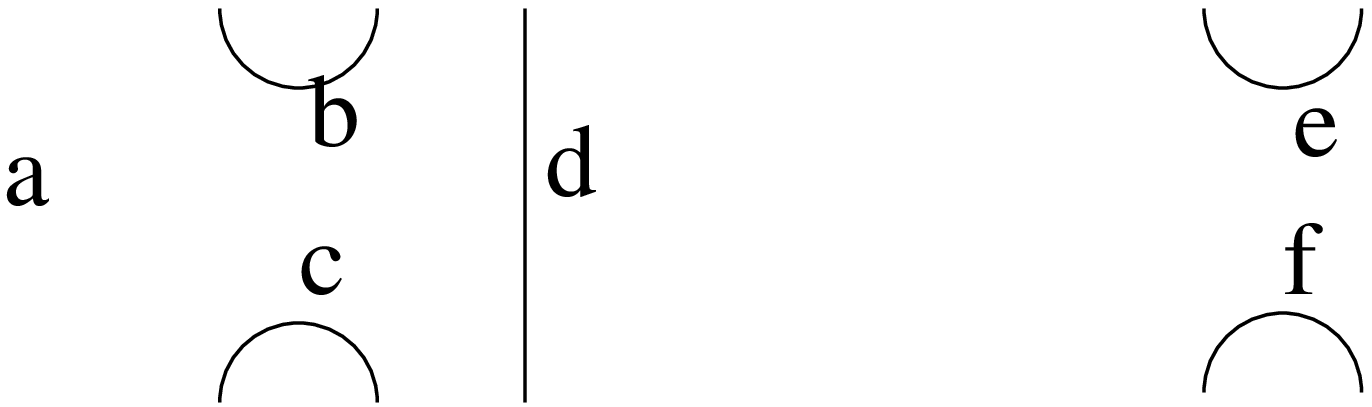}}
\end{center}
\caption{$e^{\prime}_{k+2}$}
\label{fig:pic100}
\end{figure}

 We also set $A_{k, 1} = H_{[-k, k + 1]}, \forall k \geq 2$. It is then a consequence of Proposition \ref{basic cons}(1)(i) that 
   $A_{0, 1} \subset A_{1, 1} \subset A_{2, 1} \subset A_{3, 1} \subset \cdots$ is the basic construction tower associated to
   the initial inclusion $A_{0, 1} \subset A_{1, 1}$ and for any $k \geq 0,$ if $e^{\prime}_{k+2}$ denotes the Jones projection lying in 
   $A_{k+2, 1}$ for the  basic construction of $A_{k, 1} \subset A_{k+1, 1}$, then  $e^{\prime}_{k+2}$ is given by Figure \ref{fig:pic100}.

   Further, we define inductively 
   \begin{align*}
    A_{k+2, 0} = < A_{k+1, 0}, e^{\prime}_{k+2} >
   \end{align*}
   for each $k \geq0$.
   It is well known that $A_{0, 0} \subset A_{1, 0} \overset{e^{\prime}_2} {\subset} A_{2, 0} \overset{e^{\prime}_3} {\subset} 
   A_{3, 0} \overset{e^{\prime}_4} \subset \cdots$ 
   is the basic construction tower of $A_{0, 0} \subset A_{1, 0}$.
   The following lemma explicitly describes the $C^*$-algebras $A_{k, 0}$ for $k > 0$. 
   \begin{lemma}   
     $A_{k, 0} = H_{[-k, -1]} \otimes H_{[2, k+1]}, k > 0$.  
   \end{lemma}
   \begin{proof}
   
     
    Note first that since for any $k \geq 0, A_{k, 0} \subset A_{k+1, 0} \subset A_{k+2, 0}$ is an instance of the basic construction with
    the Jones projection 
    $e^{\prime}_{k+2} \in A_{k+2, 0}$, we must have that $A_{k+2, 0} = A_{k+1, 0} \ e^{\prime}_{k+2} \ A_{k+1, 0}$.
    Our proof proceeds by induction on $k$. The case $k = 1$ is obvious from the definition of $A_{1, 0}$.    
    Now suppose that there is a positive integer
$r \geq 2$ such that the statement holds for all positive integers $k < r$ so that in particular, 
$A_{r-1, 0} = H_{[-(r-1), -1]} \otimes H_{[2, r]}$. Regard $H_{[-r, -1]} \otimes H_{[2, r+1]}$ as a subalgebra of
$H_{[-r, r+1]} \cong P(H^r)_{2r+3}$. Noting that $e^{\prime}_r \in H_{[-r, -1]} \otimes H_{[2, r+1]}$, we conclude that 
$A_{r, 0} = A_{r-1, 0} \ e^{\prime}_r \ A_{r-1, 0} \subseteq H_{[-r, -1]} \otimes H_{[2, r+1]}$. To prove the reverse inclusion, 
consider the elements in $A_{r-1, 0}$ given  by
\begin{align*}
 X = (f^{-(r-1)} / x^{-(r-1)} \rtimes \cdots \rtimes f^{-2} \rtimes x^{-1}) \otimes (f^2 \rtimes x^3 \rtimes \cdots \rtimes x^r / f^r), \ 
 Y = g^{-(r-1)} / y^{-(r-1)} \otimes y^r / g^r.
\end{align*}
A little computation in $P(H^r)_{2r+3}$ shows that $H_{[-r, -1]} \otimes H_{[2, r+1]}$ is linearly spanned by elements of the form
$X \ e^{\prime}_r \ Y$, proving the reverse inclusion. Hence the proof follows.
\end{proof}

As  $A_{0, 1} \subset A_{1, 1} \subset A_{2, 1} \subset \cdots$ (resp., $A_{0, 0} \subset A_{1, 0} \subset A_{2, 0} \subset \cdots$) is the basic
construction tower of $A_{0, 1} \subset A_{1, 1}$ (resp., $A_{0, 0} \subset A_{1, 0}$), we conclude that 
$(\cup_{k = 0}^{\infty} A_{k, 0})^{\prime \prime}$ as 
well as $(\cup_{k = 0}^{\infty} A_{k, 1})^{\prime \prime}$ are hyperfinite $II_1$ factors. 
   Now note that $\cup_{k = 0}^{\infty} A_{k, 0} = H_{(-\infty, -1]} \otimes H_{[2, \infty)}, \cup_{k = 0}^{\infty} A_{k, 1} = 
   H_{(-\infty, \infty)}$ and hence, it follows from Definition \ref{def} that $(\cup_{k = 0}^{\infty} A_{k, 0})^{\prime \prime} = \mathcal{N},
   (\cup_{k = 0}^{\infty} A_{k, 1})^{\prime \prime} = \mathcal{M}$.      
   Thus, we have proved that:
   \begin{lemma}\label{hyper}
    $\mathcal{N}$ and $\mathcal{M}$ are hyperfinite $II_1$ factors.
   \end{lemma}
 The following lemma shows that $\mathcal{N} \subset \mathcal{M}$ is of finite index  
   equal to $\delta^4$.
   \begin{lemma}\label{index}
    $[\mathcal{M} : \mathcal{N}] = \delta^4$.
   \end{lemma}
   \begin{proof}
    It is well known that (see \cite[Corollary 5.7.4]{JnsSnd1997}) $[\mathcal{M} : \mathcal{N}]$ equals the square of the norm of the inclusion matrix for 
    $A_{0, 0} \subset A_{0, 1}$ which further equals the modulus of the Markov trace $tr_{A_{0, 1}}$ for the inclusion 
    $A_{0, 0} ( = \mathbb{C}) \subset A_{0, 1} ( = H_{[0, 1]})$ which, again, by an application of Proposition 
    \ref{basic cons}(1)(ii),
    equals $\delta^4$. 
   \end{proof}  
  
  For each $k \geq 0$ and $n \geq 2$, we now define a finite-dimensional $C^*$-algebra, denoted $A_{k, n}$, as follows.   
   \begin{align*}
        A_{k, n} = 
        \begin{cases}
        H_{[-k, -1]} \otimes H_{[2, \ 2n + k + 1]}, & \mbox{if} \ n \geq 0 \hspace{1mm} \mbox{is even,} \ k > 0  \\
        H_{[-k, 2n + k - 1]}, & \mbox{if} \ n > 0 \hspace{1mm} \mbox{is odd,} \ k > 0  \\    
        H_{[-2(n-1), 1]}, & \mbox{if} \ n>0, k=0  \\
        \end{cases}
        \end{align*}
        For $k, n \geq 0$, we regard $A_{k, n}$ as a unital $C^*$-subalgebra of both $A_{k+1, n}$ and $A_{k, n+1}$ via the embeddings
        as described below.
        \begin{itemize}
         \item For any $n \geq 0, k > 0, A_{k, n}$ sits inside $A_{k+1, n}$ in the natural way.
         \item If $n \geq 0$ is even and $k > 0$, then $A_{k, n}$ sits inside $A_{k, n+1}$ in the natural way.
         \item If $n > 0$ is odd and $k > 0$, the embedding of $A_{k, n}$ inside $A_{k, n+1}$ is given by $\psi_{l, s}$ as 
         defined in the statement of Lemma \ref{exp}(ii) with $l = k, s = 2n+k-3$.
         \item If $n$ is odd, then $A_{0, n}$ is identified with the subalgebra $H_{[0, 2n-1]}$ of $A_{1, n} = H_{[-1, 2n]}$. 
         \item If $n > 0$ is even, then $A_{0, n}$ is identified with the subalgebra $H_{[2, 2n+1]}$ of
         $A_{1, n} = H_{-1} \otimes H_{[2, 2n+2]}$.         
         \item Embedding of $A_{0, n}$ inside $A_{0, n+1}$ is natural for all $n \geq 0$.
        \end{itemize}
        Thus, we have a grid $\{ A_{k, n} : k, n \geq 0 \}$ of finite-dimensional $C^*$-algebras. 
        The following remark contains several useful facts concerning the grid $\{ A_{k, n} : k, n \geq 0 \}$.
        \begin{remark}\label{B}
         \begin{itemize}
          \item[(i)] We have already seen (as an application of Lemma \ref{comm sq}) that the square of finite-dimensional $C^*$-algebras as shown 
          in Figure \ref{fig:pic37} is a symmetric commuting square with respect to $tr A_{1,1}$ which is a Markov trace of modulus $\delta^4$
          for the inclusion $A_{0,1} \subset A_{1,1}$ and all the inclusions are connected. Further, by Lemma \ref{index}, 
          $[\mathcal{M} : \mathcal{N}] = \delta^4$ where $\mathcal{N} = (\cup_{k=0}^{\infty} A_{k, 0})^{\prime \prime}, \mathcal{M} =
          (\cup_{k=0}^{\infty} A_{k, 1})^{\prime \prime}$.
           \item[(ii)]It follows from the embedding prescriptions that the following diagram (see Figure \ref{fig:cd}) commutes for all $k, n \geq 0$. 
         \begin{figure}[!h]
\begin{center}
\psfrag{a}{\Huge $A_{k, n+1}$}
\psfrag{b}{\Huge $\subset$}
\psfrag{c}{\Huge $A_{k+1, n+1}$}
\psfrag{x}{\Huge $A_{k, n}$}
\psfrag{y}{\Huge $\subset$} 
\psfrag{z}{\Huge $A_{k+1, n}$}
\psfrag{p}{\Huge $\cup$}
\psfrag{q}{\Huge $\cup$}
\resizebox{2.5cm}{!}{\includegraphics{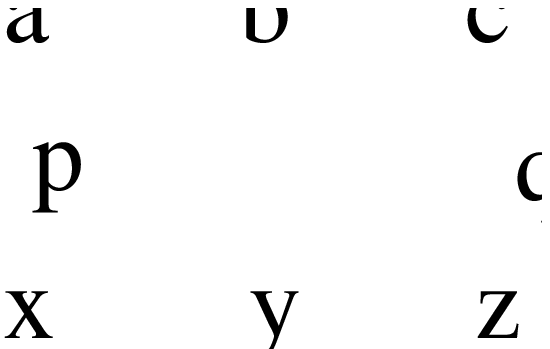}}
\end{center}
\caption{Commuative diagram}
\label{fig:cd}
\end{figure}
     \item[(iii)]It is a direct consequence of Proposition
        \ref{basic cons} that
        for any $k, n \geq 0, A_{k, n} \subset A_{k, n+1} \subset A_{k, n+2}$ is an instance of the basic construction
         and further, $tr_{A_{k, n+1}}$ is a Markov trace of modulus $\delta^4$ for the inclusion 
        $A_{k, n} \subset A_{k, n+1}$. Also, if $e_{k, n+2}$ ($k \geq 0, n \geq 0$) denotes the Jones projection lying in $A_{k, n+2}$, the 
        result of basic construction for $A_{k, n} \subset A_{k, n+1}$, then $e_{k, n+2}$ is shown in Figure \ref{fig:D7}.
         \begin{figure}[!h]
\begin{center}
\psfrag{a}{\Huge $\delta^{-2}$}
\psfrag{2}{\Huge $2$}
\psfrag{c}{\Huge $k+2$}
\psfrag{b}{\Huge $k$}
\psfrag{d}{\Huge $2n+k+1$}
\psfrag{e}{\Huge $2n+k+1$}
\resizebox{6cm}{!}{\includegraphics{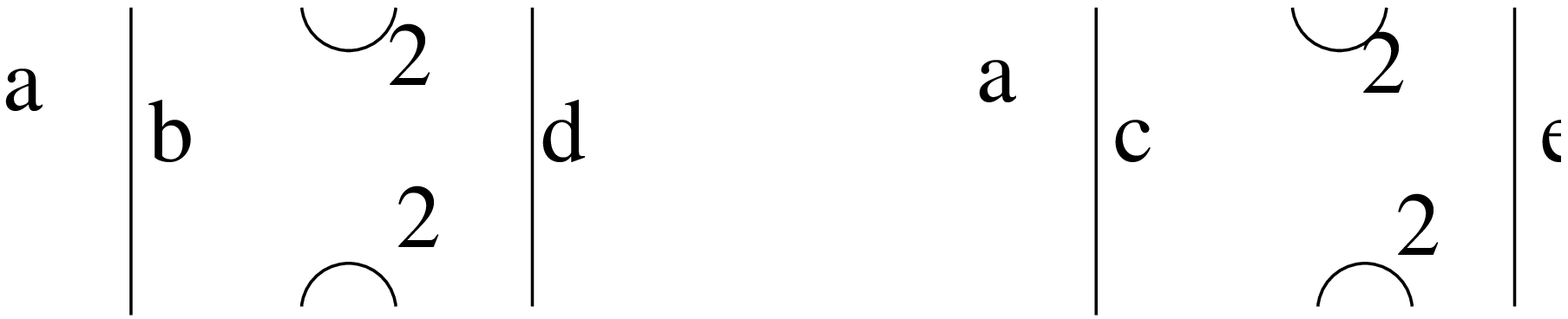}}
\end{center}
\caption{$e_{k, n+2}$ ($k > 0, n$ odd or $k = 0, n \geq 0$) (left) and $e_{k, n+2}$ ($k > 0, n$ even) (right)}
\label{fig:D7}
\end{figure} 
\item[(iv)] For any $k \geq 0, n \geq 2$, the embedding of $A_{k, n}$ inside $A_{k+1, n}$ carries $e_{k, n}$ to $e_{k+1, n}$. 
         \end{itemize}
        \end{remark}
        Consider the tower of finite pre-von Neumann algebras $\cup_{k = 0}^{\infty}A_{k, 0} \subset \cup_{k = 0}^{\infty}A_{k, 1}
        \subset \cup_{k = 0}^{\infty}A_{k, 2} \subset \cdots$. Observe that for any $m \geq 0, \  
        \cup_{k = 0}^{\infty}A_{k, m} \subset \cup_{k = 0}^{\infty}A_{k, m+1}$ is a compatible pair. Note also that  
        $\cup_{k = 0}^{\infty} A_{k, m} = H_{(-\infty, -1]} \otimes H_{[2, \infty)} \ \mbox{or} \ 
    H_{(-\infty, \infty)}$ according as $m$ is even or odd.
    For each $m \geq 0$,  we define $\mathcal{M}_m:= (\cup_{k = 0}^{\infty}A_{k, m})^{\prime \prime}$. Then $\mathcal{M}_m = (H_{(-\infty, -1]} 
    \otimes H_{[2, \infty)})^{\prime \prime} \ \mbox{or} \ H^{\prime \prime}_{(-\infty, \infty)}$ according as $m$ is even or odd.
    In view of the facts concerning the grid $\{ A_{k, n} : k, n \geq 0 \}$ as mentioned in Remark \ref{B} and \cite[Proposition 5.7.5]{JnsSnd1997},
    one can conclude that:
    \begin{proposition}\label{towerbasic}
    $\mathcal{M}_0 (=\mathcal{N}) \subset \mathcal{M}_1 (=\mathcal{M}) \subset \mathcal{M}_2 \subset \mathcal{M}_3 \subset \cdots$ is the basic construction tower of $\mathcal{N}
    \subset \mathcal{M}$. 
    \end{proposition}

 \subsection{Computation of the relative commutants.}
     We now proceed to compute the relative commutants. By virtue of Ocneanu's compactness theorem (see \cite[Theorem 5.7.6]{JnsSnd1997}), the relative commutant 
$\mathcal{N}^{\prime} \cap \mathcal{M}_k$ ($k > 0$) is given 
by
\begin{align*}
 \mathcal{N}^{\prime} \cap \mathcal{M}_k = A_{0,k} \cap (A_{1,0})^{\prime}, \ k \geq 1.
\end{align*}
 The following proposition describes the spaces $ A_{0, k} \cap (A_{1, 0})^{\prime}, k \geq 1$. 
      
\begin{figure}[!h]
\begin{center}
\psfrag{1}{\huge $2$}
\psfrag{3}{\huge $2m$}
\psfrag{a}{\huge $\delta^{-2}$}
\psfrag{x}{\huge $\mathcal{N}^{\prime} \cap \mathcal{M}_{m+1}$}
\psfrag{y}{\huge $\mathcal{N}^{\prime} \cap \mathcal{M}_{m}$}
\psfrag{z}{\Huge $\tilde{Q}_{m+1}$}
\psfrag{w}{\Huge $\tilde{Q}_m$}
\psfrag{b}{\Huge $\longrightarrow$}
\psfrag{c}{\Huge $\cup$}
\resizebox{8.5cm}{!}{\includegraphics{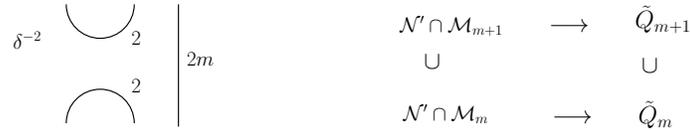}}
\end{center}
\caption{$\tilde{e}_{m+2}$ (left) and a commutative diagram (right)}
\label{fig:pic73}
\end{figure}

\begin{figure}[!h]
\begin{center}
\psfrag{a}{\Huge $\tilde{Q}_{m+1}$}
\psfrag{b}{\Huge $Q_{m+1}$}
\psfrag{c}{\Huge $\tilde{Q}_m$}
\psfrag{d}{\Huge $Q_m$}
\psfrag{e}{\Huge $\overset{\gamma_{m+1}}{\longrightarrow}$}
\psfrag{f}{\Huge $\cup$}
\psfrag{g}{\Huge $\overset{\gamma_m}{\longrightarrow}$}
\resizebox{3.0cm}{!}{\includegraphics{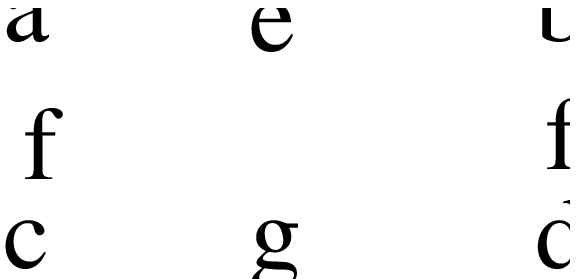}}
\end{center}
\caption{Commutative diagram}
\label{fig:D8}
\end{figure}

\begin{figure}[!h]
\begin{center}
\psfrag{1}{\Huge $2$}
\psfrag{3}{\Huge $2(m-2)$}
\psfrag{a}{\Huge $\delta^{-2}$}
\resizebox{4cm}{!}{\includegraphics{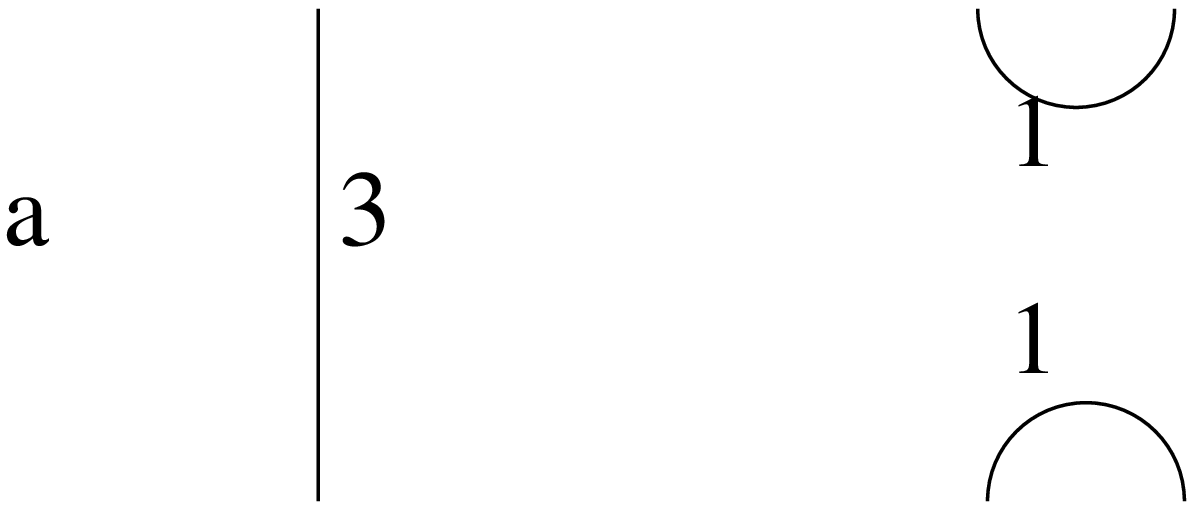}}
\end{center}
\caption{$e_m$}
\label{fig:pic77}
\end{figure}

\begin{figure}[!h]
\begin{center}
\psfrag{a}{\Huge $\mathcal{N}^{\prime} \cap \mathcal{M}_{m+1}$}
\psfrag{b}{\Huge $Q^{op}_{m+1}$}
\psfrag{c}{\Huge $\mathcal{N}^{\prime} \cap \mathcal{M}_m$}
\psfrag{d}{\Huge $Q^{op}_m$}
\psfrag{e}{\Huge $\overset{\Psi_{m+1}}{\longrightarrow}$}
\psfrag{f}{\Huge $\cup$}
\psfrag{g}{\Huge $\overset{\Psi_m}{\longrightarrow}$}
\resizebox{3.5cm}{!}{\includegraphics{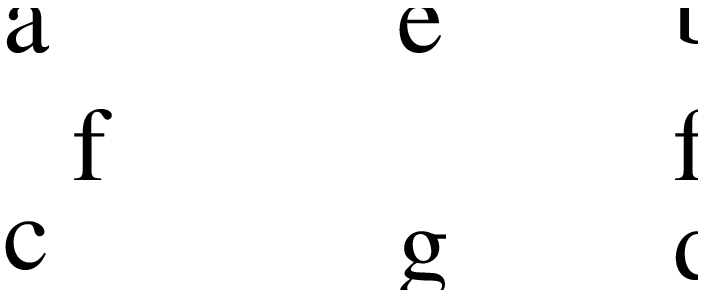}}
\end{center}
\caption{Commutative diagram}
\label{fig:D9}
\end{figure}

      \begin{proposition}\label{rel comm}      
      Let $m \geq 1$ be an integer. Then, 
     \begin{itemize}
       \item[(a)]$A_{0, 2m} \cap (A_{1, 0})^{\prime} \cong \tilde{Q}_{2m}:= \{ X \in H_{[2, 4m]} : X  \ \mbox{commutes with} \ \Delta_{m-1}(x) \in 
       \otimes_{i = 1}^m H_{4i-1}, \forall x \in H \}$ and       
       \item[(b)] $A_{0, 2m-1} \cap (A_{1, 0})^{\prime} \cong \tilde{Q}_{2m-1}:= \{ X \in H_{[0, 4m-4]} : X  \ \mbox{commutes with} \ 
       \Delta_{m-1}(x) \in \otimes_{i = 0}^{m-1} H_{4i-1}, \forall x \in H \}$.
      \end{itemize}
  \end{proposition}
  \begin{proof}
   \begin{itemize}
 \item[(a)]We observe from the embedding prescription of $A_{1, 0}$ inside $A_{1, 2m}$ that given 
 $Y = x \otimes f \in A_{1, 0}$, its image in $A_{1, 2m}$ is given by
   \begin{align*}
    \Delta_m(x) \otimes f = \otimes_{i=1}^{m+1} x_i \otimes f \in \otimes_{i=0}^m H_{4i-1} \otimes H_{4m+2}. 
   \end{align*}
   Recall that $A_{0, 2m}$ is identified with the subalgebra $H_{[2, 4m+1]}$ of
         $A_{1, 2m}$. 
         Thus, given $X \in A_{0, 2m}$, its image in $A_{1, 2m}$ is given by $1 \otimes (X \rtimes \epsilon)$. 
   Consequently, $A_{0, 2m} \cap (A_{1, 0})^{\prime}$ equals
   \begin{align*}
   & \{X \in H_{[2, 4m+1]} : 1 \otimes (X \rtimes \epsilon) \in H_{-1} \otimes H_{[2, 4m+2]} \ \mbox{commutes with}\\
  & \Delta_m(x) \otimes f \in \otimes_{i = 0}^m H_{4i-1} 
   \otimes H_{4m+2}, \forall f \otimes x \in H^* \otimes H \}.
   \end{align*}
   Thus given $X \in A_{0, 2m} \cap (A_{1, 0})^{\prime}$, since $X (\in H_{[2, 4m+1]})$ commutes with $H_{4m+2}$, 
   it commutes with every element of $H_{[4m+2, \infty)}$ and consequently, by an appeal to Lemma \ref{commutants}
   we conclude that $X$ indeed lies in $H_{[2, 4m]}$. 
  Thus, $A_{0, 2m} \cap (A_{1, 0})^{\prime}$ can be identified with
   \begin{align*}
   & \{X \in H_{[2, 4m]} : 1 \otimes X \in H_{-1} \otimes \in H_{[2, 4m]} \ \mbox{commutes with} \ 
   \Delta_m(x) \in \otimes_{i = 0}^{m} H_{4i-1}, \forall x \in H \}
   \end{align*} 
   and a little thought should convince the reader that this space equals 
   \begin{align*}
    \{ X \in H_{[2, 4m]} : X  \ \mbox{commutes with} \ \Delta_{m-1}(x) \in 
       \otimes_{i = 1}^m H_{4i-1}, \forall x \in H \}
   \end{align*}
 \item[(b)] It follows from the embedding formula of $A_{1, 0}$ inside $A_{1, 2m}$ that given 
 $Y = x \otimes f \in A_{1, 0}$, its image in $A_{1, 2m-1}$ is given by 
   \begin{align*}
    \Delta_{m-1}(x) \otimes f \in \otimes_{i = 0}^{m-1} H_{4i-1} \otimes H_{4m-2}.
   \end{align*}
   Recall that $A_{0, 2m-1}$ is identified with the subalgebra $H_{[0, 4m-3]}$ of $A_{1, 2m-1}$.
   Thus, given $X \in A_{0, 2m-1}$, its image in $A_{1, 2m-1}$ is given by $1 \rtimes X \rtimes \epsilon$.
   Consequently, $A_{0, 2m-1} \cap (A_{1, 0})^{\prime}$ equals
   \begin{align*}
   & \{X \in H_{[0, 4m-3]} : 1 \rtimes X \rtimes \epsilon \in H_{[-1, 4m-2]} \ \mbox{commutes with}\\  
  & \Delta_{m-1}(x) \otimes f \in \otimes_{i = 0}^{m-1} H_{4i-1} \otimes H_{4m-2}, \forall x \otimes f \in H \otimes H^*\}. 
   \end{align*}
   Similar kind of argument as in the proof of part (a) shows that this space can be identified with 
   \begin{align*}
    \{ X \in H_{[0, 4m-4]} : X  \ \mbox{commutes with} \ 
       \Delta_{m-1}(x) \in \otimes_{i = 0}^{m-1} H_{4i-1}, \forall x \in H \}.
   \end{align*}
   \end{itemize}
  \end{proof}    
     It follows from Remark \ref{B}(iii) that the Jones projection lying in $\mathcal{N}^{\prime} \cap \mathcal{M}_{m+2} = 
   A_{0, m+2} \cap (A_{1, 0})^{\prime}$ ($m \geq 0$)
   is given by $e_{0, m+2}$ (see Figure \ref{fig:D7}), which, under the identification of 
   $A_{0, m+2} \cap (A_{1, 0})^{\prime}$ with $\tilde{Q}_{m+2}$ as given by Proposition \ref{rel comm}, is easily seen to be identified with the 
   projection $\tilde{e}_{m+2}$ in $\tilde{Q}_{m+2}$ as shown on the left in Figure \ref{fig:pic73}.
\begin{remark}\label{embtil}
  It is worth knowing the embedding of $\tilde{Q}_m$ inside $\tilde{Q}_{m+1}$ ($m \geq 1$).
  It follows easily from the embedding formulae of $A_{1, m}$ inside $A_{1, m+1}$ and Proposition \ref{rel comm} that given $X \in \tilde{Q}_m$,
  it sits inside $\tilde{Q}_{m+1}$ as $\epsilon \rtimes 1 \rtimes X$ and the diagram on the right in Figure \ref{fig:pic73} commutes
 where each horizontal arrow indicates the $*$-isomorphism.
 \end{remark}
 Consider the $*$-anti-isomorphism of $H_{[2, 4m]}$ (resp., $H_{[0, 4m-4]}$) onto $A(H^*)_{4m-1}$ (resp., $A(H^*)_{4m-3}$) as given by 
 Lemma \ref{anti1}. Let $Q_{2m}$ (resp., $Q_{2m-1}$) denote the image of $\tilde{Q}_{2m} (\subset H_{[2, 4m]})$ (resp., $\tilde{Q}_{2m-1}
 (\subset H_{[0, 4m-4]}$)) inside
 $A(H^*)_{4m-1}$ (resp., $A(H^*)_{4m-3}$) under this $*$-anti-isomorphism. One can easily see that 
 \begin{align*}
  Q_{2m} = \{ &X \in A(H^*)_{4m-1}: X \ \mbox{commutes with} \ \epsilon \rtimes x_1 \rtimes \epsilon \rtimes 1 \rtimes \epsilon \rtimes x_2 \rtimes
  \epsilon \rtimes 1 \rtimes \epsilon \rtimes x_3 \rtimes \cdots\\
  &\rtimes x_{m-1} \rtimes \epsilon \rtimes 1 \rtimes \epsilon \rtimes x_m \rtimes \epsilon \ \mbox{where} \ \otimes_{i=1}^m x_i = 
  \Delta_{m-1}(x), \forall x \in H \},\\
  Q_{2m-1} = \{ &X \in A(H^*)_{4m-3}: X \rtimes 1 \ \mbox{commutes with} \ \epsilon \rtimes x_1 \rtimes \epsilon \rtimes 1 \rtimes \epsilon \rtimes x_2 \rtimes
  \epsilon \rtimes 1 \rtimes \epsilon \rtimes x_3 \rtimes \cdots\\
  &\rtimes x_{m-1} \rtimes \epsilon \rtimes 1 \rtimes \epsilon \rtimes x_m \ \mbox{where} \ \otimes_{i=1}^m x_i = \Delta_{m-1}(x), 
  \forall x \in H \}.
 \end{align*}
 Let $\gamma_m$ denote the aforementioned $*$-anti-isomorphism of 
 $\tilde{Q}_m$ onto $Q_m$ ($m \geq 1$). Now, keeping in mind the embedding of $\tilde{Q}_m$ inside $\tilde{Q}_{m+1}$ as mentioned in Remark 
 \ref{embtil} and noting the obvious fact that $Q_m$ sits inside $Q_{m+1}$ naturally, we conclude that the following diagram 
 (see Figure \ref{fig:D8}) commutes. 
%
%
%
%
%
%
Further, if $e_m \in Q_m (m \geq 2)$ denotes the projection which is the image of $\tilde{e}_m \in \tilde{Q}_m$ under 
$\gamma_m$, 
  it is then not hard to see that $e_m$ is given by Figure \ref{fig:pic77}.
Let $Q_m^{op}$ denote the opposite algebra of $Q_m$. That is, $Q_m^{op} = Q_m$ as vector spaces, only the multiplication of $Q_m^{op}$ is 
opposite to that of $Q_m$. 
  Obviously, the identity map of $Q_m$ onto $Q_m^{op}$ is $*$-anti-isomorphism. 
  For each $m \geq 1$, let $\Psi_m: \mathcal{N}^{\prime} \cap \mathcal{M}_m \rightarrow$ $Q^{op}_m$ denote the following composite map: 
    \begin{align*}
     \mathcal{N}^{\prime} \cap \mathcal{M}_m \xrightarrow{*-\mbox{isom}} \tilde{Q}_m \xrightarrow{\gamma_m (*- \mbox{anti-isom})} 
     Q_m \xrightarrow{\mbox{Identity}} Q^{op}_m. 
    \end{align*}
    Obviously $\Psi_m$ is a $*$-isomorphism for each $m \geq 1$ and it carries $e_{0, m}$ to $e_m$ for $m \geq 2$. The commutative 
    diagram on the right in Figure \ref{fig:pic73} and the commutative diagram of Figure \ref{fig:D8} together imply commutativity of the following diagram
    (see Figure \ref{fig:D9}).

Once again applying Ocneanu's compactness theorem we obtain that 
\begin{align*}
 \mathcal{M}^{\prime} \cap \mathcal{M}_m = A_{0, m} \cap A_{1, 1}^{\prime}, m \geq 2.
\end{align*}
Proceeding along the same line of argument as in the proof of Proposition \ref{rel comm}, one can
show that: 
\begin{lemma}\label{second}
The $*$-isomorphism $\Psi_m$ of $\mathcal{N}^{\prime} \cap \mathcal{M}_m$ onto $Q^{op}_m$ carries $\mathcal{M}^{\prime} \cap \mathcal{M}_m$ 
onto the subspace of $Q^{op}_m$ given by
 \begin{align*}
  \{X = f^0 \rtimes x^1 \rtimes f^2 \rtimes \cdots \rtimes f^{2m-2} \in Q^{op}_m (\subset A(H^*)_{2m-1}^{op}) : f^0 = \epsilon, x^1 = 1 \}.
   \end{align*}
\end{lemma}

We conclude this section with the following lemma.
\begin{lemma}\label{irr}
  $\mathcal{N} \subset \mathcal{M}$ is irreducible.
 \end{lemma}
 \begin{proof}
  An appeal to Lemma \ref{commutants} immediately shows that the space $Q_1 = \{f \in H^* : f \rtimes 1 \ \mbox{commutes with} \ \epsilon 
  \rtimes x, \forall x \in H\}$ is trivial 
  so that $\mathcal{N} \subset \mathcal{M}$ is irreducible. 
 \end{proof}
 
 \section{On planar algebra of $\mathcal{N} \subset \mathcal{M}$}
 In this section we give an explicit description of the planar algebra associated to $\mathcal{N} \subset \mathcal{M}$ and
 it turns out to be an interesting planar subalgebra of $^{*(2)}\!P(H^*)$. 
 
\begin{figure}[!h]
\begin{center}
\psfrag{a}{\huge $1$}
\psfrag{b}{\Huge $Fx_1$}
\psfrag{c}{\huge $2$}
\psfrag{e}{\Huge $Fx_k$}
\psfrag{f}{\huge $1$}
\psfrag{h}{\Huge $Fx_{k+1}$}
\psfrag{z}{\huge $2$}
\psfrag{u}{\Huge $Fx_{k-1}$}
\psfrag{w}{\Huge $Fx_k$}
\psfrag{k}{\Huge $Fx_{2k}$}
\psfrag{j}{\huge $2$}
\psfrag{l}{\huge $1$}
\psfrag{m}{\huge $2$}
\psfrag{n}{\Huge $Fx_{2k-1}$}
\psfrag{r}{\huge $2$}
\psfrag{s}{\huge $2$}
\psfrag{o}{\Huge $x_k$}
\psfrag{q}{\Huge $Fx_1$}
\psfrag{r}{\huge $2$}
\psfrag{s}{\huge $2$}
\psfrag{p}{\huge $1$}
\psfrag{v}{\huge $2$}
\psfrag{Y}{\Huge $X$}
\psfrag{X}{\Huge $X$}
\psfrag{y}{\Huge $Fx_{k+1}$}
\psfrag{2}{\huge $1$}
\psfrag{x}{\huge $2$}
\psfrag{g}{\Large $1$}
\psfrag{1}{\Large $1$}
\psfrag{2}{\huge $1$}
\psfrag{3}{\huge $2$}
\psfrag{i}{\huge $2$}
\resizebox{11cm}{!}{\includegraphics{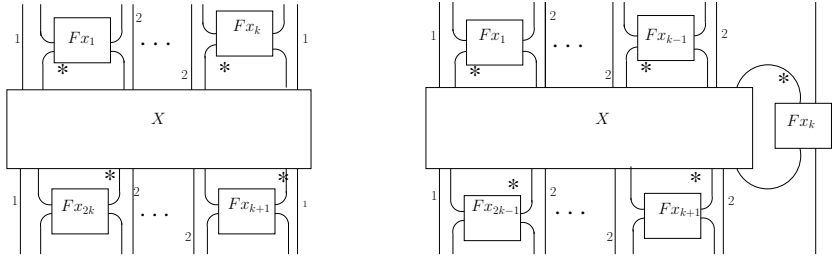}}
\end{center}
\caption{$\alpha^{2k}_x(X)$ (Left) and $\alpha^{2k-1}_x(X)$ (Right), where $k \geq 1$}
\label{fig:D10}
\end{figure}

For each $m \geq 1$, consider the linear map $\alpha^m : H \rightarrow End(P(H^*)_{2m})$ defined for $x \in H$ and
$X \in P(H^*)_{2m}$ by Figure \ref{fig:D10} where the notation $\alpha^m_x$ stands for 
$\alpha^m(x)$. 

With the help of the maps $\alpha^m$ defined above we give an alternative description of the spaces $Q_m$ ($m \geq 1$). 
\begin{proposition}\label{descrip}
    For any $m \geq 1$,
    \begin{align*}
    Q_m = \{X \in P(H^*)_{2m} : \alpha^m_h(X) = X \}.    
    \end{align*}
    \end{proposition}
    
     Before we proceed to prove Proposition \ref{descrip}, we pause for a simple Hopf algebraic lemma.
    \begin{lemma}\label{Hopf}
    Let $k \geq 1$ be an integer.

     (a) If $X \in A(H^*)_{4k-3}$, the following are equivalent :   
     \begin{itemize}
    \item[(i)] $X \rtimes 1$ commutes with $\epsilon \rtimes x_1 \rtimes \epsilon \rtimes 1 \rtimes \epsilon \rtimes x_2 \rtimes
    \epsilon \rtimes 1 \rtimes \epsilon \rtimes x_3 \rtimes \cdots \rtimes \epsilon \rtimes x_k, \forall x \in H$,
    
    \item[ii)] $(\epsilon \rtimes h_1 \rtimes \epsilon \rtimes 1 \rtimes \epsilon \rtimes h_2 \rtimes
    \epsilon \rtimes 1 \rtimes \epsilon \rtimes h_3 \rtimes \cdots \rtimes \epsilon \rtimes h_k) \ (X \rtimes 1) \  
    (\epsilon \rtimes Sh_{2k} \rtimes \epsilon \rtimes 1 \rtimes \epsilon \rtimes Sh_{2k-1} \rtimes
    \epsilon \rtimes 1 \rtimes \epsilon \rtimes Sh_{2k-2} \rtimes \cdots \rtimes \epsilon \rtimes Sh_{k+1}) = X \rtimes 1$.
    \end{itemize}   
    
    \item[(b)] If $X \in A(H^*)_{4k-1}$, the following are equivalent :
    \begin{itemize}
     \item[(i)] $X$ commutes with $\epsilon \rtimes x_1 \rtimes \epsilon \rtimes 1 \rtimes \epsilon \rtimes x_2 \rtimes
    \epsilon \rtimes 1 \rtimes \epsilon \rtimes x_3 \rtimes \cdots \rtimes \epsilon \rtimes x_k \rtimes \epsilon, \forall x \in H$, 
    
    \item[(ii)] $(\epsilon \rtimes h_1 \rtimes \epsilon \rtimes 1 \rtimes \epsilon \rtimes h_2 \rtimes
    \epsilon \rtimes 1 \rtimes \epsilon \rtimes h_3 \rtimes \cdots \rtimes \epsilon \rtimes h_k \rtimes \epsilon) \ X \  
    (\epsilon \rtimes Sh_{2k} \rtimes \epsilon \rtimes 1 \rtimes \epsilon \rtimes Sh_{2k-1} \rtimes
    \epsilon \rtimes 1 \rtimes \epsilon \rtimes Sh_{2k-2} \rtimes \cdots \rtimes \epsilon \rtimes Sh_{k+1} \rtimes \epsilon) = X$.
    \end{itemize}
%
%
%
    \end{lemma}

    \begin{proof}
    \begin{itemize}
     \item[(a)]  
     $(i) \Rightarrow (ii):$ obvious.
     
     $(ii) \Rightarrow (i):$
     Note that 
     \begin{align*}
     &(\epsilon \rtimes x_1 \rtimes \epsilon \rtimes 1 \rtimes \epsilon \rtimes x_2 \rtimes
    \epsilon \rtimes 1 \rtimes \epsilon \rtimes x_3 \rtimes \cdots \rtimes \epsilon \rtimes x_k) (X \rtimes 1)\\
    &= (\epsilon \rtimes x_1 \rtimes \epsilon \rtimes 1 \rtimes \epsilon \rtimes x_2 \rtimes
    \epsilon \rtimes 1 \rtimes \epsilon \rtimes x_3 \rtimes \cdots \rtimes \epsilon \rtimes x_k)\\
    &(\epsilon \rtimes h_1 \rtimes \epsilon \rtimes 1 \rtimes \epsilon \rtimes h_2 \rtimes
    \epsilon \rtimes 1 \rtimes \epsilon \rtimes h_3 \rtimes \cdots \rtimes \epsilon \rtimes h_k) \ (X \rtimes 1)\\  
    &(\epsilon \rtimes Sh_{2k} \rtimes \epsilon \rtimes 1 \rtimes \epsilon \rtimes Sh_{2k-1} \rtimes
    \epsilon \rtimes 1 \rtimes \epsilon \rtimes Sh_{2k-2} \rtimes \cdots \rtimes \epsilon \rtimes Sh_{k+1})\\
    &= (\epsilon \rtimes x_1h_1 \rtimes \epsilon \rtimes 1 \rtimes \epsilon \rtimes x_2h_2 \rtimes
    \epsilon \rtimes 1 \rtimes \epsilon \rtimes x_3h_3 \rtimes \cdots \rtimes \epsilon \rtimes x_kh_k) \ (X \rtimes 1)\\  
    &(\epsilon \rtimes Sh_{2k} \rtimes \epsilon \rtimes 1 \rtimes \epsilon \rtimes Sh_{2k-1} \rtimes
    \epsilon \rtimes 1 \rtimes \epsilon \rtimes Sh_{2k-2} \rtimes \cdots \rtimes \epsilon \rtimes Sh_{k+1}).
      \end{align*}
      Now using the formula $xh_1 \otimes Sh_2 = h_1 \otimes Sh_2x$ (and hence, $\Delta_{k-1}(xh_1) \otimes \Delta_{k-1}(Sh_2) = 
      \Delta_{k-1}(h_1) \otimes \Delta_{k-1}(Sh_2x))$, the desired result follows easily. 
     \item[(b)] Proof is similar to that of part (a).
      \end{itemize} 
      \end{proof}
      We are now ready to prove Proposition \ref{descrip}.
      \begin{proof}[Proof of Proposition \ref{descrip}]
      We prove the result for $m$ odd say, $m = 2k-1 (k \geq 1)$, leaving the case when $m$ is even for the reader. This is an 
      immediate consequence of Lemma \ref{Hopf}(a) that the space $Q_{2k-1}$ can equivalently be described as
      \begin{align*}
      Q_{2k-1} = &\{ X \in A(H^*)_{4k-3}: (\epsilon \rtimes h_1 \rtimes \epsilon \rtimes 1 \rtimes \epsilon \rtimes h_2 \rtimes
    \epsilon \rtimes 1 \rtimes \epsilon \rtimes h_3 \rtimes \cdots \rtimes \epsilon \rtimes h_k) \ (X \rtimes 1)\\  
   & (\epsilon \rtimes Sh_{2k} \rtimes \epsilon \rtimes 1 \rtimes \epsilon \rtimes Sh_{2k-1} \rtimes
    \epsilon \rtimes 1 \rtimes \epsilon \rtimes Sh_{2k-2} \rtimes \cdots \rtimes \epsilon \rtimes Sh_{k+1}) = X \rtimes 1 \}. 
      \end{align*}
      Interpreting this equivalent description of $Q_{2k-1}$ pictorially in $P(H^*)$,
      we note that $Q_{2k-1}$ consists of
      precisely those elements $X \in P(H^*)_{4k-2}$ such that the pictorial equation of Figure \ref{fig:eqv} holds. 
       Now applying the conditional expectation tangle $E^{4k-2}_{4k-1}$ and then using once the multiplication relation in $P(H^*)$ 
      we reduce the element on the left in Figure \ref{fig:eqv}
      to that on the left in Figure \ref{fig:eqv1}. On the other hand an application of the conditional expectation tangle $E^{4k-2}_{4k-1}$
      and then an appeal to the modulus relation reduces the element on the right in Figure \ref{fig:eqv} to that on the right in Figure 
      \ref{fig:eqv1}. Now note that for any $x \in H$,
      \begin{align*}
       Fx_1 Fx_2 = \delta^2 \phi_2(x_1) \tilde{\phi}_2(x_2) \phi_1 \tilde{\phi}_1 = 
       \delta^2 (\phi_2 \tilde{\phi}_2)(x) \phi_1 \tilde{\phi}_1 = \delta F(x).
      \end{align*}
Using this Hopf algebraic identity, one can easily see that the element on the left in Figure \ref{fig:eqv1} just equals 
$\delta \alpha^{2k-1}_h(X)$ and consequently, we obtain from the pictorial equation of Figure \ref{fig:eqv1} that $\alpha^{2k-1}_h(X) = X$,
showing that $Q_{2k-1} \subseteq \{X \in P(H^*)_{4k-2}: \alpha^{2k-1}_h(X) = X \}$. To see the reverse inclusion, 
let $X \in P(H^*)_{4k-2}$ be such that $\alpha^{2k-1}_h(X) = X$.      
\begin{figure}[!h]
\begin{center}
\psfrag{a}{\huge $Fh_1$}
\psfrag{b}{\huge $Fh_2$}
\psfrag{c}{\huge $Fh_k$}
\psfrag{d}{\huge $Fh_{k+1}$}
\psfrag{e}{\huge $=$}
\psfrag{f}{\huge $Fh_{2k-1}$}
\psfrag{g}{\huge $Fh_{2k}$}
\psfrag{X}{\huge $X$}
\psfrag{o}{\huge $\cdots$}
\psfrag{X}{\huge $X$}
\resizebox{9cm}{!}{\includegraphics{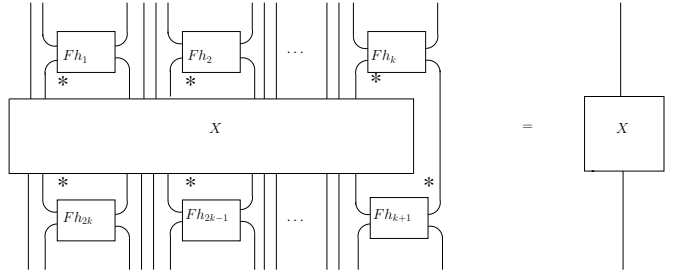}}
\end{center}
\caption{A characterisation of $X \in Q_{2k-1}$}
\label{fig:eqv}
\end{figure}
\begin{figure}[!h]
\begin{center}
\psfrag{a}{\huge $Fh_1$}
\psfrag{b}{\huge $Fh_2$}
\psfrag{e}{\huge $=$}
\psfrag{f}{\huge $Fh_{2k-1}$}
\psfrag{g}{\huge $Fh_{2k}$}
\psfrag{X}{\huge $X$}
\psfrag{o}{\huge $\cdots$}
\psfrag{j}{\huge $\delta$}
\psfrag{k}{\huge $Fh_k \ Fh_{k+1}$}
\psfrag{h}{\huge $4k-2$}
\resizebox{9.5cm}{!}{\includegraphics{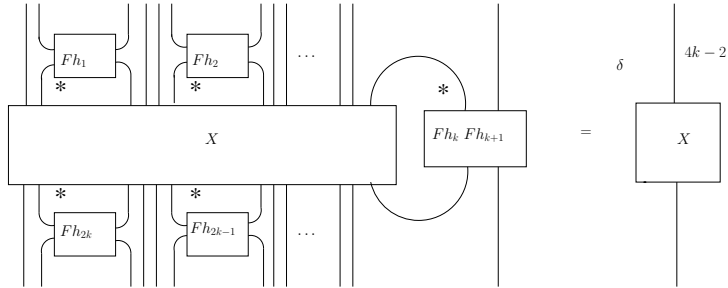}}
\end{center}
\caption{Another characterisation of $X \in Q_{2k-1}$}
\label{fig:eqv1}
\end{figure}
In order to verify that $X \in Q_{2k-1}$, it suffices to show, by Lemma \ref{Hopf}(a), that 
\begin{align}\label{lhs}
 &(\epsilon \rtimes \tilde{h}_1 \rtimes \epsilon \rtimes 1 \rtimes \epsilon \rtimes \tilde{h}_2 \rtimes
    \epsilon \rtimes 1 \rtimes \epsilon \rtimes \tilde{h}_3 \rtimes \cdots \rtimes \epsilon \rtimes \tilde{h}_k) \ (X \rtimes 1)\nonumber \\  
 &(\epsilon \rtimes S\tilde{h}_{2k} \rtimes \epsilon \rtimes 1 \rtimes \epsilon \rtimes S\tilde{h}_{2k-1} \rtimes
    \epsilon \rtimes 1 \rtimes \epsilon \rtimes S\tilde{h}_{2k-2} \rtimes \cdots \rtimes \epsilon \rtimes S\tilde{h}_{k+1}) = X \rtimes 1,  
\end{align}
where $\tilde{h}$, of course, denotes another copy of the unique non-zero idempotent integral of $H$.
Keeping in mind that $\alpha^{2k-1}_h(X) = X$ and applying multiplication rule in $P(H^*)$, it is easy to see that the element of 
$P(H^*)_{4k-1}$ as depicted in Figure \ref{fig:eqv2} represents the left hand side of \eqref{lhs}. 
\begin{figure}[!h]
\begin{center}
\psfrag{a}{\huge $F(\tilde{h}_1 h_1)$}
\psfrag{b}{\huge $F(\tilde{h}_{k-1}h_{k-1})$}
\psfrag{c}{\huge $F(\tilde{h}_{k+2} h_{k+1})$}
\psfrag{d}{\huge $F(\tilde{h}_{2k} h_{2k-1})$}
\psfrag{e}{\huge $F\tilde{h}_k$}
\psfrag{f}{\huge $F\tilde{h}_{k+1}$}
\psfrag{g}{\huge $Fh_k$}
\psfrag{X}{\huge $X$}
\psfrag{o}{\huge $\cdots$}
\resizebox{10cm}{!}{\includegraphics{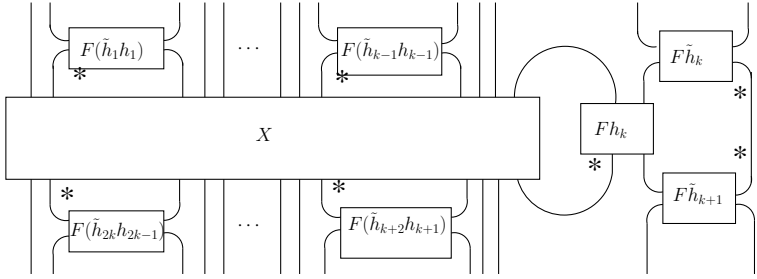}}
\end{center}
\caption{Pictorial description of left hand side of \eqref{lhs}}
\label{fig:eqv2}
\end{figure}
\begin{figure}[!h]
\begin{center}
\psfrag{a}{\huge $Fx_1$}
\psfrag{b}{\huge $Fx_2$}
\psfrag{c}{\huge $Fy$}
\psfrag{d}{\huge $F(xy)$}
\psfrag{e}{\huge $=$}
\resizebox{4.5cm}{!}{\includegraphics{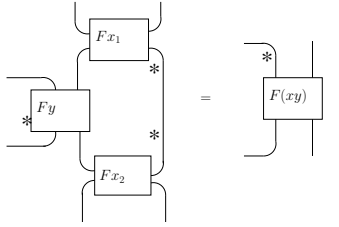}}
\end{center}
\caption{A relation in $P(H)^*$}
\label{fig:eqv3}
\end{figure}
Now, using the Hopf algebraic identity $\phi_1 \otimes \phi_2 f = \phi_1 Sf \otimes \phi_2$ and the antipode, integral and exchange relations
in $P(H^*)$, the equation of Figure \ref{fig:eqv3} can be verified to hold in $P(H^*)$ for all $x, y \in H$ and we leave this pleasant 
verification to the reader.
Using the equation of Figure \ref{fig:eqv3} and the fact that $\tilde{h}_1 h_1 \otimes \tilde{h}_2 h_2 \otimes \cdots \otimes \tilde{h}_{2k-1}
h_{2k-1} = \Delta_{2k-2}(\tilde{h}h) = \Delta_{2k-2}(h)$, one can easily see that the element in Figure \ref{fig:eqv2} indeed equals 
$\alpha^{2k-1}_h(X) \rtimes 1$ which, again by virtue of the fact that $\alpha^{2k-1}_h(X) = X$, equals $X \rtimes 1$ and the proof is complete.
 \end{proof}
     Thus, we have a family $\{Q_m : m \geq 1 \}$ of vector spaces where for each $m \geq 1, Q_m$ 
    is a subspace of  $P(H^*)_{2m} = ^{(2)}\!\!\!P(H^*)_{m}$. Setting $Q_{0, \pm} = \mathbb{C}$, we note that 
    $Q := \{Q_m: m \in Col \}$ is a subspace of $^{(2)}\!P(H^*)$. The following proposition shows that $Q$ is indeed a planar subalgebra of 
    $^{(2)}\!P(H^*)$.
   \begin{proposition}\label{planar}    
     $Q$ is a planar subalgebra of  $^{(2)}\!P(H^*)$. 
    \end{proposition}
    \begin{proof}
     By an appeal to Theorem \ref{generating}, it suffices to prove that $Q$ is closed under the action of the following  set of tangles
     \begin{align*}
      \{1^{0,+}, 1^{0, -}\} \cup \{R_k^k:k \geq 2\} \cup \{M_{k, k}^k, E^k_{k+1}, I_k^{k+1}: k \in Col \}.
     \end{align*}
    It is obvious that $Q$ is closed under the action of the tangles $1^{0, \pm}$ and $M_{k, k}^k, I_k^{k+1} (k \in Col)$.

     To see that $Q$ is closed under the action of the rotation tangle $R_k^k \ (k \geq 2)$, we note that for any 
     $X \in Q_k \ (k \geq 2)$, we have 
     \begin{align*}
      Z^{^{(2)}\!P(H^*)}_{R_k^k}(X) =   Z^{^{(2)}\!P(H^*)}_{R_k^k}(\alpha^k_h(X)) = \alpha^k_h(  Z^{^{(2)}\!P(H^*)}_{R_k^k}(X))  
     \end{align*}
     where the first equality follows from the fact that $\alpha^k_h(X) = X$ and to see the second equality we need to use the Hopf algebra
     identity $h_1 \otimes h_2 \otimes \cdots \otimes h_l = h_2 \otimes h_3 \otimes 
     \cdots \otimes h_l \otimes h_1$ ($l \geq 2$) which basically follows from $h_1 \otimes h_2 = h_2 \otimes h_1$ (which essentially 
     expresses traciality of $h$).

     Verifying that $Q$ is closed under the action of $E^k_{k+1} \ (k \geq 1)$ amounts to verification of the following identity
     \begin{align}\label{coexp}
      Z^{^{(2)}\!P(H^*)}_{E^k_{k+1}}(X) =  Z^{^{(2)}\!P(H^*)}_{E^k_{k+1}}(\alpha^{k+1}_h(X)) = \alpha^{k}_h(Z^{^{(2)}\!P(H^*)}_{E^k_{k+1}}(X))
     \end{align}
     for $X \in Q_{k+1}$. We observe that the first equality of \eqref{coexp} is obvious from the fact $\alpha^{k+1}_h(X) = X$. If $k$ is even, the second 
     equality of \eqref{coexp} follows easily by pictorially representing in $P(H^*)$ both sides of the equality and then using the Relation T. 
     When $k$ is odd, verification
     of the second equality needs more effort. We first claim that the equation of Figure \ref{fig:eqv4} holds in $P(H^*)$ for all $x \in H$.
     To see this note that for any $x \in H$, $Fx_1 \otimes Fx_2 = \delta^2 (\phi_2 \tilde{\phi}_2)(x) \phi_1 \otimes \tilde{\phi}_1$
     which, by an appeal to the formula $\phi_1 \otimes \phi_2 f = \phi_1 Sf \otimes \phi_2$, equals $\delta^2 \ \phi_2(x) \ \phi_1 S 
     \tilde{\phi}_2 \otimes \tilde{\phi}_1$. Using this expression for $Fx_1 \otimes Fx_2$, we can reduce the element on the left in Figure \ref{fig:eqv4} to that
     on the left in Figure \ref{fig:eqv5}. Using the multiplication and antipode relations in $P(H^*)$, one can easily see that the equation
     of Figure \ref{fig:eqv5} holds. Finally, an application of the exchange relation and then the integral relation reduces the element on the right in 
     Figure \ref{fig:eqv5} to that on the right in \ref{fig:eqv4}, establishing our claim. When $k$ is odd, to verify the second equality of 
     \eqref{coexp}, we just need to present pictorially elements on both sides of the equality in $P(H^*)$ and then apply the equation
     of Figure \ref{fig:eqv4}. This completes the proof of the proposition.     
     \begin{figure}[!h]
\begin{center}
\psfrag{a}{\Huge $Fx_1$}
\psfrag{b}{\Huge $Fx_2$}
\psfrag{c}{\Huge $Fx$}
\psfrag{e}{\Huge $=$}
\resizebox{7.5cm}{!}{\includegraphics{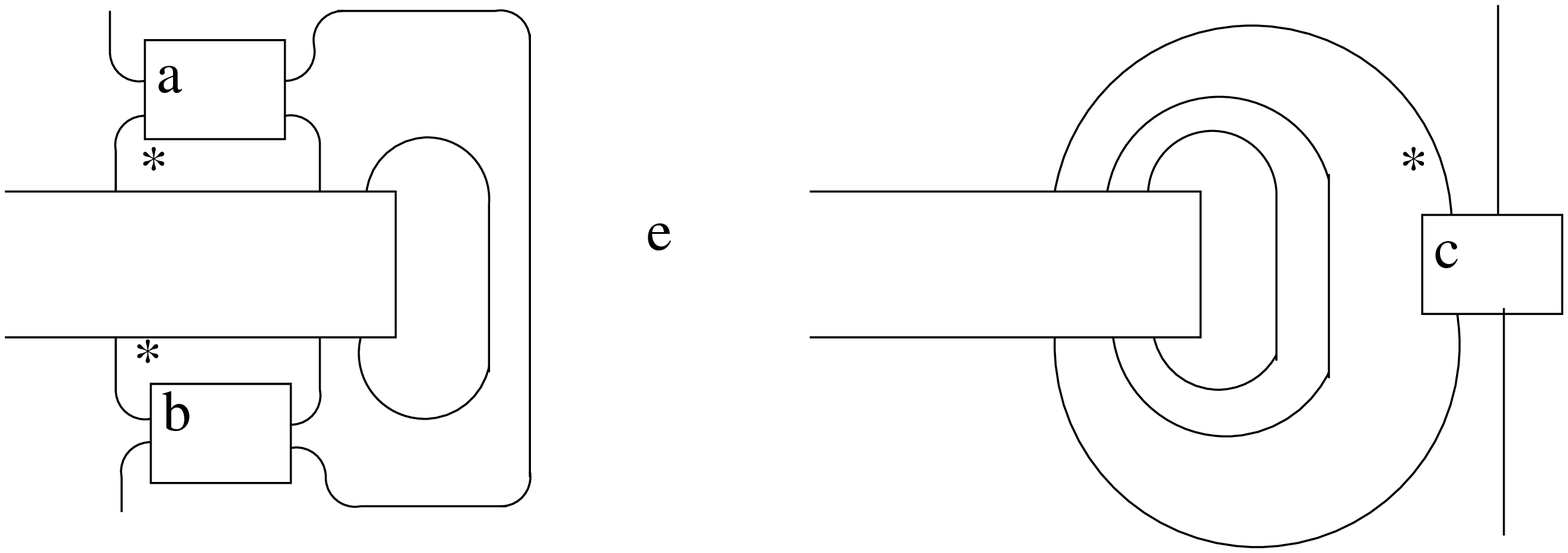}}
\end{center}
\caption{An equation in $P(H^*)$}
\label{fig:eqv4}
\end{figure}
 \begin{figure}[!h]
\begin{center}
\psfrag{a}{\Huge $\phi_1$}
\psfrag{b}{\Huge $\tilde{\phi}_2$}
\psfrag{c}{\Huge $\tilde{\phi}_1$}
\psfrag{e}{\Huge $=$}
\psfrag{f}{\Huge $\phi_1 S \tilde{\phi}_2$}
\psfrag{g}{\Huge $\tilde{\phi}_1$}
\psfrag{x}{\Huge $\delta^2 \ \phi_2(x)$}
\psfrag{y}{\Huge $\delta^2 \ \phi_2(x)$}
\resizebox{9cm}{!}{\includegraphics{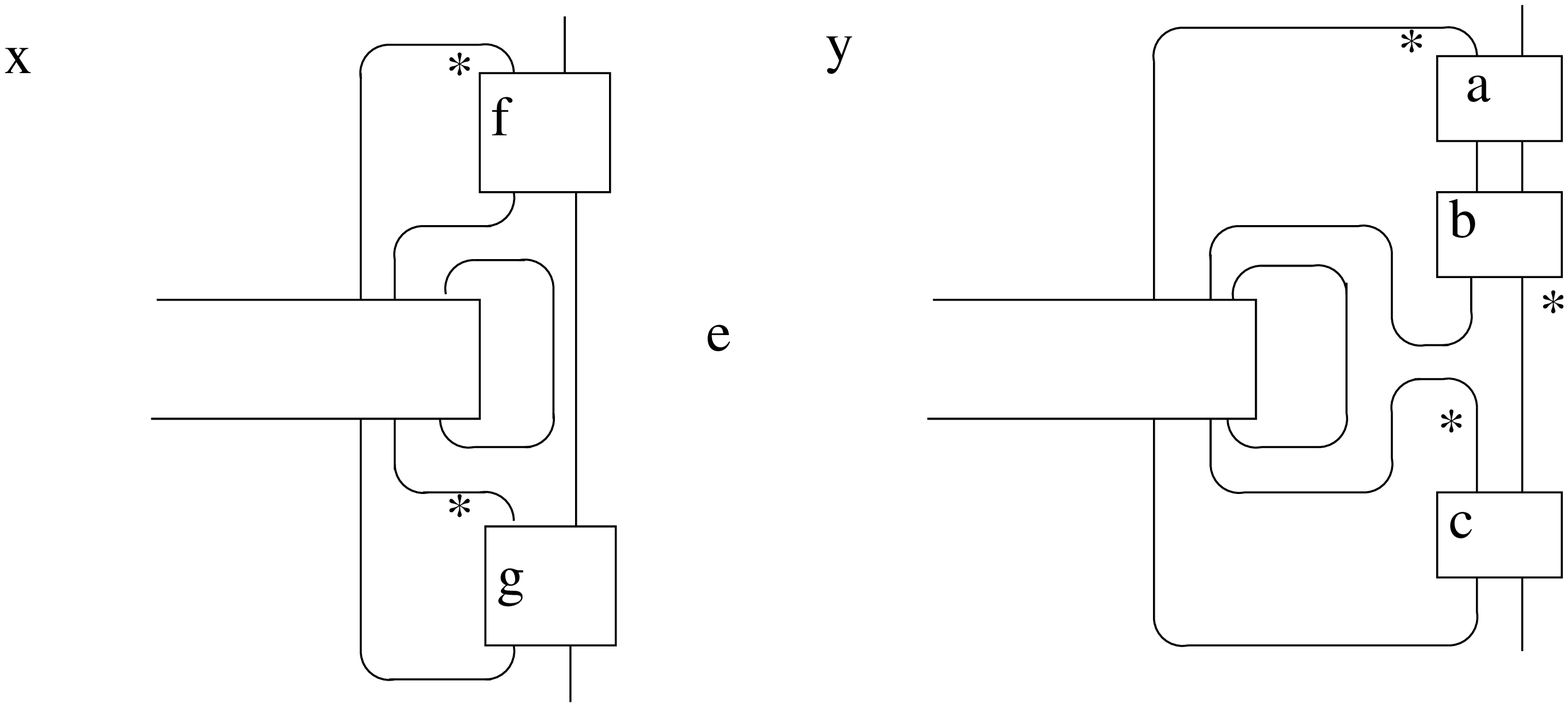}}
\end{center}
\caption{An equation in $P(H^*)$}
\label{fig:eqv5}
\end{figure}
    \end{proof}
    As an immediate corollary of Proposition \ref{planar} we obtain that:
    \begin{corollary} \label{planar1}
     $^*\!Q$ is a planar subalgebra of  $^{*(2)}\!P(H^*)$.
    \end{corollary}

      The next proposition shows that $P^{\mathcal{N} \subset \mathcal{M}}$, the planar algebra associated to $\mathcal{N} \subset \mathcal{M}$,
      is given by the adjoint of the planar algebra $Q$.    
    \begin{proposition}\label{pl alg}
     $P^{\mathcal{N} \subset \mathcal{M}} = {^*\!Q}$.
    \end{proposition}
    
    \begin{proof}
    In order to establish that $P^{\mathcal{N} \subset \mathcal{M}} = {^*\!Q}$, we need to verify all the conditions of Theorem \ref{jones}.
    We note first that the subfactor $\mathcal{N} \subset \mathcal{M}$ is extremal since it is irreducible by Lemma \ref{irr}.    
    
    Obviously, $^*\!Q$ has modulus $\delta^2$. 
    In view of the commutative diagram in Figure \ref{fig:D9}, we note that for each $m \geq 1$, the identification of 
    $^*\!Q_m ( = Q^{op}_m)$ with $\mathcal{N}^{\prime} \cap \mathcal{M}_m$ as $*$-algebras (via $\Psi^{m}$) respects inclusion, verifying the condition (i) of 
    Theorem \ref{jones}. Verification of the condition (ii) follows from the trivial observation that 
    $\delta^{2} \ e_m = Z_{\mathcal{E}^m}^{^*\!Q}(1)$ (see Figure \ref{fig:pic77} for the definition of $e_m$) 
    for each $m \geq 2$. Since $^*\!Q$ is a planar subalgebra of  $^{*(2)}\!P(H^*)$, condition (iv) of Theorem \ref{jones} is
    automatically satisfied. We next observe, by virtue of Theorem \ref{jones} and Remark \ref{a}, that the linear map 
    $Z^{^{*(2)}\!P(H^*)}_{(E^{\prime})^k_k}$ induced by the tangle
    $(E^{\prime})^k_k \ (k \geq 1)$ is such that $\delta^{-2} Z^{^{*(2)}P(H^*)}_{(E^{\prime})^k_k}$
     equals the conditional expectation of $^{*(2)}\!P(H^*)_k (\cong A(H^*)_{2k-1}^{op})$ onto the subspace $^{*(2)}\!P(H^*)_{1, k} 
     (\cong \{f^0 \rtimes x^1 \rtimes \cdots \rtimes f^{2m-2} \in A(H^*)_{2k-1}^{op} : f^0 = \epsilon, x^1 = 1\} )$. 
     Now, by Corollary \ref{planar1}, $^*\!Q$ is a planar subalgebra of  
     $^{*(2)}\!P(H^*)$ and hence, $Z^{^*\!Q}_{(E^{\prime})^k_k}$
     carries $^*\!Q_k (= \mathcal{N}^{\prime} \cap \mathcal{M}_k)$ into $^*\!Q_k \cap {^{*(2)}\!P(H^*)_{1, k}}$. It follows immediately
     from Lemma \ref{second} that $^*\!Q_k \cap {^{*(2)}\!P(H^*)_{1, k}}$ indeed equals 
     $\mathcal{M}^{\prime} \cap \mathcal{M}_k$ so that $\delta^{-2} Z^{\!^*Q}_{(E^{\prime})^k_k}$ induces 
     the conditional expectation of $\mathcal{N}^{\prime} \cap \mathcal{M}_k$ onto $\mathcal{M}^{\prime} \cap \mathcal{M}_k$, verifying the
     condition (iii) of Theorem \ref{jones}. This completes the proof of the proposition.

    \end{proof}
    
    We collect the results of the previous statements into a single main theorem.
    \begin{theorem}\label{maintheorem}
      $^*\!Q$ is a planar subalgebra of $^{*(2)}\!P(H^*)$ and $^*\!Q = P^{\mathcal{N} \subset \mathcal{M}}$. 
      For each $k \geq 1, {^*\!Q_k}$ consists of all $X \in 
      P(H^*)_{2k}$ such that the element shown in Figure \ref{fig:charac20} equals $X$.
     \begin{figure}[!h]

\begin{center}

\psfrag{1}{\huge $Fh_1$}

\psfrag{2}{\huge $Fh_2$}

\psfrag{3}{\huge $Fh_3$}

\psfrag{km1}{\huge $Fh_{k-1}$}

\psfrag{k}{\huge $Fh_k$}





%

\psfrag{X}{\huge $X$}

\resizebox{6cm}{!}{\includegraphics{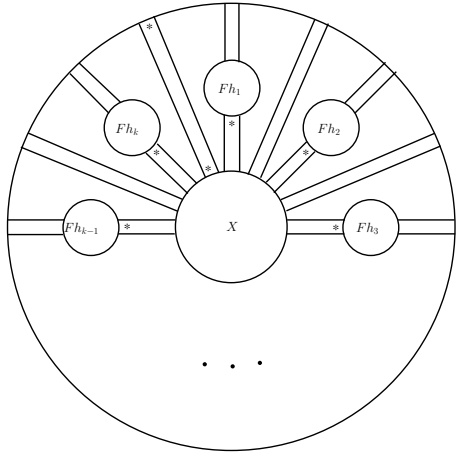}}

\end{center}

\caption{Characterisation of the image}

\label{fig:charac20}

\end{figure}

    \end{theorem}
    \begin{proof}
     It follows immediately from Proposition \ref{pl alg} after observing that $\alpha^k_h(X)$ 
     in Figure \ref{fig:D10} is equivalent to the element in Figure \ref{fig:charac20}.  
    \end{proof}
    
\section{Main result}
In this section we show that $\mathcal{M} \cong \mathcal{N} \rtimes D(H)^{cop}$ by showing that the planar algebra of $\mathcal{N} \subset 
\mathcal{M}$ is isomorphic to the planar algebra of $R \subset R \rtimes D(H)^{cop}$ for an outer action of $D(H)^{cop}$ on the hyperfinite
$II_1$ factor $R$ where $D(H)$ is the Drinfeld  double of $H$. We refer to \cite{Mjd2002} for the Drinfeld double construction.

\begin{figure}[!h]

\begin{center}

\psfrag{1}{\huge $Fh_1$}

\psfrag{2}{\huge $Fh_2$}

\psfrag{3}{\huge $Fh_3$}

\psfrag{km1}{\huge $Fh_{k-1}$}

\psfrag{k}{\huge $Fh_k$}

\psfrag{2k-4}{\huge $2k-4$}

\psfrag{2k-3}{\huge $2k-3$}

\psfrag{2k-2}{\huge $2k-2$}

\psfrag{2k-1}{\huge $2k-1$}

\psfrag{X}{\huge $X$}

\resizebox{11.5cm}{!}{\includegraphics{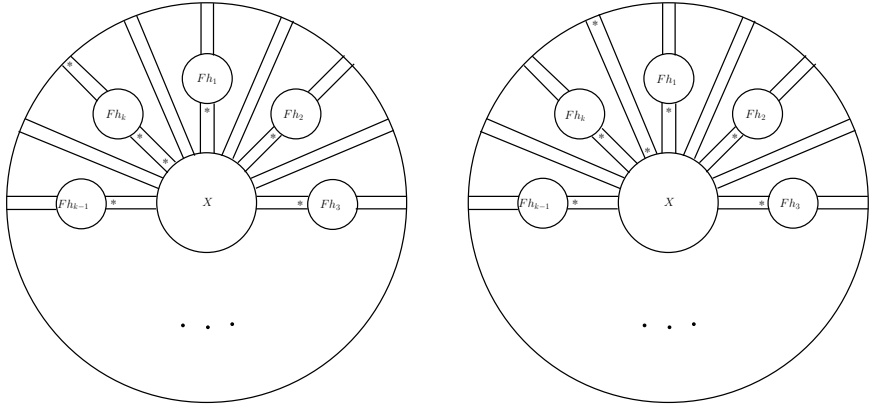}}

\end{center}

\caption{Characterisation of the image}

\label{fig:charac22}

\end{figure}

In \cite{DeKdy2016}, the authors produced an explicit embedding of the planar algebra of the Drinfeld double of a finite-dimensional, semisimple and 
cosemisimple Hopf algebra (and hence, in particular, a Kac algebra) $H$ into $\!^{(2)}P(H^*)$ and characterised the image. The following theorem, which is a 
reformulation of \cite[Theorem 10]{DeKdy2016}, gives an explicit characterisation of the image $Q$ of the planar algebra of $D(H)$ 
inside $\!^{(2)}P(H^*)$. It is worth mentioning that this 
statement uses the newer version of planar algebras with spaces indexed by $\{ \mathbb{N} \cup 0\} \times \{+, -\}$. We refer to $\S 2.1$ for 
the older and newer notions of planar algebras. 
\begin{theorem}\label{qt}
\cite[Theorem 10]{DeKdy2016}\\
 $Q$ is characterised as follows: $Q_{k,+}$ (resp. $Q_{k,-}$) is the set of all $X \in P(H^*)_{2k,+}$  such the element on the left 
(resp. right) in Figure \ref{fig:charac22} equals $X$ where $h \in H$ is the unique non-zero idempotent integral.
\end{theorem}

Let $Q_1$ and $Q_2$ be planar subalgebras of $\!^{(2)}P(H^*)$ in the older sense defined by setting $(Q_1)_{0, \pm} = (Q_2)_{0, \pm} = \mathbb{C}$
and for any positive integer $k, (Q_1)_{k} = Q_{{k, +}}, \ (Q_2)_{k} = Q_{k, -}$. Then note that $Q_1$ is is isomorphic to the planar algebra 
associated to $D(H)$ in the older sense and corresponds to the subfactor 
$R^{D(H)} \subset R$ whereas $Q_2$ corresponds to its dual i.e., to $R \subset R \rtimes D(H)$ for some outer action of $D(H)$ on $R$.
Using the fact - see \cite[Remark 4.18]{KdySnd2004} - that for any Kac algebra $K$, $^*\!P(K) \cong P(K^{op})$, we conclude that 
$^*\!Q_1$ corresponds to the subfactor $R^{D(H)^{op}} \subset R$ and consequently, $^*\!Q_2$ corresponds to 
$R \subset R \rtimes D(H)^{op} \cong R \subset R \rtimes D(H)^{cop}$. It is now an immediate consequence of these observations together 
with Theorem \ref{qt} and the description of the planar algebra of $\mathcal{N} \subset \mathcal{M}$ as given by Theorem \ref{maintheorem}  that 
$P^{\mathcal{N} \subset \mathcal{M}} = {^*\!Q_2}$. Finally, an application of Jones' theorem - see Theorem \ref{jones} - immediately yields 
our main result. 
\begin{theorem}\label{result}
The quantum double inclusion of $R^H \subset R$ is isomorphic to 
$R \subset R \rtimes D(H)^{op} (\cong R \subset R \rtimes D(H)^{cop})$ for some outer action of $D(H)^{op}
(\cong D(H)^{cop})$ on the hyperfinite $II_1$ factor $R$.
\end{theorem}

 \section*{acknowledgement}
 The author sincerely thanks Prof. Vijay Kodiyalam, IMSc, for a careful reading of the manuscript and suggestions for improvement and also for  
 his support, constructive comments and several fruitful discussions during the whole project. He also thanks Prof. David Evans for his
 support in making it possible to attend the inspiring programme ``Operator Algebras: Subfactors and their Applications'' held at INI, 2017 
 while this work was going on. The author is supported by the 
 National Board of Higher Mathematics (NBHM), India.


\begin{thebibliography}{99}
 \bibitem{BlaMnt1985}
 Robert J. Blattner and Susan Montgomery,
 {\em A duality theorem for Hopf module algebras},
 Journal of Algebra, 95, No. 1 (1985), 153-172. 
 
 
 \bibitem{DeKdy2015}
 Sandipan De and Vijay Kodiyalam, 
 {\em Note on infinite iterated crossed products of
 Hopf algebras and the Drinfeld double},
 Journal of Pure and Applied Algebra, 219, No. 12 (2015), 5305-5313.
 
 \bibitem{DeKdy2016}
Sandipan De and Vijay Kodiyalam, 
{\em Planar algebras, cabling and the Drinfeld double},
Quantum Topology, 9, No. 1 (2018), 141-165.

 
 \bibitem{EK1993}
 D. E. Evans and Y. Kawahigashi,
 {\em Quantum Symmetries on operator
algebras}, Oxford Univ. Press, 1998.

\bibitem{EK1995}
D. E. Evans and Y. Kawahigashi, 
{\em On Ocneanu's theory of asymptotic inclusions for subfactors, topological
quantum field theories and quantum doubles}, 
Internat. J. Math., 6, No. 2 (1995), 205-228.
 
 \bibitem{GHJ1989}
F.M. Goodman, P. de la Harpe and V.F.R Jones,
{\em Coxeter graphs and towers of algebras},
MSRI Publ., 14, Springer, New York, 1989.

\bibitem{Izumi2000}
M. Izumi,
{\em The structure of sectors associated with Longo-Rehren inclusions. I. General theory},
Comm. Math. Phys., 213, No. 1 (2000), 127-179.

\bibitem{Izumi2001}
M. Izumi,
{\em The structure of sectors associated with Longo-Rehren inclusions. II. Examples}, Rev. Math. Phys.,
13, No. 5 (2001), 603-674.

\bibitem{Jjo2008}
S. Jijo,
{\em Planar algebra associated to the Asymptotic inclusion of a Kac
algebra subfactor},
Ph.D. thesis, Homi Bhabha National Institute, 2008.

\bibitem{Jns1999}
V. F. R. Jones, 
{\em Planar algebras I}, To appear in New Zealand J. Math.
arXiv:math/9909027.


\bibitem{JnsSnd1997}
V.F.R. Jones and V.S. Sunder, 
{\em Introduction to Subfactors},
Cambridge University Press, 1997.



\bibitem{KdyLndSnd2003}
Vijay Kodiyalam, Zeph Landau, V.S. Sunder,
{\em The planar algebra associated to a Kac algebra},
Proc. Indian Acad. Sci. Math. Sci. 113, No. 1 (2003), 15-51.



\bibitem{KdySnd2004}
Vijay Kodiyalam and V. S. Sunder, 
{\em On Jones' planar algebras}, Journal of Knot
Theory and its Ramifications, 13, No. 2 (2004), 219-247.

\bibitem{KdySnd2006}
Vijay Kodiyalam and V. S. Sunder, 
{\em The planar algebra of a semisimple and
cosemisimple Hopf algebra},
Proc. Indian Acad. Sci. Math. 116,
No. 4 (2006), 443-458.


\bibitem{KdySnd2008}
Vijay Kodiyalam and V.S. Sunder,
{\em From subfactor planar algebras to subfactors},
Internat. J. Math. 20, No. 10  (2009), 1207-1231.


\bibitem{Mjd2002}
 Shahn Majid, 
 {\em A Quantum Groups Primer}, London Mathematical Society Lecture
 Note Series 292, Cambridge University Press, Cambridge, 2002.
 
 \bibitem{Mug2003}
 Michael \"{M}uger,
 {\em From subfactors to categories and topology II. The quantum double of tensor categories and subfactors}, 
 J. Pure Appl. Algebra, 180 (2003), 159-219.
 
 \bibitem{Oc1988}
 A. Ocneanu,
 {\em Quantized groups, string algebras and Galois theory for algebras, in Operator algebras and applications}, Vol. 2, vol. 136 of 
 London Math. Soc. Lecture Note Ser., Cambridge Univ. Press, Cambridge, 1988,
119-172.

\bibitem{Y1993}
T. Yamanouchi,
 {\em Construction of an outer action of a finite-dimensional Kac algebra on the AFD factor of type $II_1$},
 Internat. J. Math. 4, No. 06 (1993), 1007-1045.
  \end{thebibliography}
\end{document}